\documentclass[final]{siamart1116}
\pdfoutput=1

\usepackage{mathtools}
\usepackage{amssymb}
\usepackage{stmaryrd}
\usepackage{color}
\usepackage{pdflscape}
\usepackage{epsfig}
\usepackage{amsmath}
\usepackage{amsfonts}

\usepackage{graphicx}
\usepackage{epstopdf}
\usepackage{multirow}
\usepackage{hhline}
\usepackage{txfonts}
\usepackage{mathrsfs}
\usepackage{color}
\usepackage{cite}
\usepackage{algorithmic}

\usepackage{sidecap}
\usepackage{subcaption}
\usepackage{mathrsfs}
\usepackage{enumitem}

\newtheorem{thm}{Theorem}[section]

\newtheorem{rem}[thm]{Remark}
\numberwithin{equation}{section}
\title{Vibration Analysis of Geometrically Nonlinear and Fractional 
Viscoelastic Cantilever Beams} 

\author{
	Pegah Varghaei
	\footnote{ D\lowercase{epartment of} M\lowercase{echanical} E\lowercase{ngineering} \& D\lowercase{epartment of} C\lowercase{omputational} M\lowercase{athematics}, S\lowercase{cience}, \lowercase{and}, E\lowercase{ngineering}
		,	
		M\lowercase{ichigan} S\lowercase{tate} U\lowercase{niversity}, 428 S S\lowercase{haw} L\lowercase{ane}, E\lowercase{ast} L\lowercase{ansing}, MI 48824, USA,  C\lowercase{orresponding author; zayern@msu.edu}}
	, Ehsan Kharazmi
	\footnote{ D\lowercase{ivision of} A\lowercase{pplied} M\lowercase{athematics},	
		B\lowercase{rown} U\lowercase{niversity}, 182 G\lowercase{eorge} S\lowercase{treet}, P\lowercase{rovidence}, RI 02912, USA}
	, Jorge Suzuki
	\footnote{ D\lowercase{epartment of} M\lowercase{echanical} E\lowercase{ngineering} \& D\lowercase{epartment of} C\lowercase{omputational} M\lowercase{athematics}, S\lowercase{cience}, \lowercase{and}, E\lowercase{ngineering}
		,	
		M\lowercase{ichigan} S\lowercase{tate} U\lowercase{niversity}, 428 S S\lowercase{haw} L\lowercase{ane}, E\lowercase{ast} L\lowercase{ansing}, MI 48824, USA}
	, Mohsen Zayernouri
	\footnote{D\lowercase{epartment of} M\lowercase{echanical} E\lowercase{ngineering} \& D\lowercase{epartment of} S\lowercase{tatistics and} P\lowercase{robability},
		M\lowercase{ichigan} S\lowercase{tate} U\lowercase{niversity}, 428 S S\lowercase{haw} L\lowercase{ane}, E\lowercase{ast} L\lowercase{ansing}, MI 48824, USA,  C\lowercase{orresponding author; zayern@msu.edu}}
}

\begin{document}

\maketitle

\begin{abstract}
We investigate the nonlinear vibration of a fractional viscoelastic cantilever 
beam, subject to base excitation, where the viscoelasticity takes the general 
form of a distributed-order fractional model, and the beam curvature introduces 
geometric nonlinearity into the governing equation. We utilize the extended 
Hamilton's principle to derive the governing equation of motion for specific 
material distribution functions that lead to fractional Kelvin-Voigt 
viscoelastic model. By spectral decomposition in space, the resulting governing 
fractional PDE reduces to nonlinear time-fractional ODEs. We use direct 
numerical integration in the decoupled system, in which we observe the 
anomalous power-law decay rate of amplitude in the linearized model. We further 
develop a semi-analytical scheme to solve the nonlinear equations, using method 
of multiple scales as a perturbation technique. We replace the expensive 
numerical time integration with a cubic algebraic equation to solve for 
frequency response of the system. We observe the super sensitivity of response 
amplitude to the fractional model parameters at free vibration, and bifurcation 
in steady-state amplitude at primary resonance.

%Moreover, we create a vibration-based feature library of system response for different values of fractional derivative order to develop a machine learning tool that returns estimation of system parameter in an inverse problem setting. We use this framework to assess the heath of considered beam by assuming a multi-stage failure threshold in the model parameters.   

%We employ fractional constitutive equations in modeling anomalous behavior of material deterioration in the progressive mechanical damage process under undesired dynamic loading. The micro-scale changes of the internal structure of material during the damage process is usually observed via changes in macro-scale material properties. Such alteration and its inherent associated uncertainty propagates into the mathematical model parameters and therefore, affects the time response characteristics under prescribed loading. We develop a mathematical framework to study the state of damage propagation and thus, monitor and estimate the vibration-based health of structure. We collect proper learning data sets and their associated uncertainties by simulating the stochastic models for a vast range of parameters distribution. Then, we train the estimator algorithms via stochastic supervised learning, where we derive the cost functions, using non-linear regression and neural network. Finally, we test the performance of developed scheme with some synthetic vibration data.

	%
\end{abstract}

\pagestyle{myheadings}
\thispagestyle{plain}

%
%%%%%%%%%%%%%%%%%%%%%%%%%%%
\section{Introduction}
\label{Sec: Introduction}
%%%%%%%%%%%%%%%%%%%%%%%%%%%
%
%\noindent \textbf{Damage and Cantilever beams:}

%Structures under periodic loading experience different types of damage and fatigue. 
%Discuss different types of damage
%focus on vibration-based identification of damage propagation
As mechanical structures undergo vibration due to many unwanted external 
forcing sources, they experience cyclic stresses resulting in structural 
damage. Such damage process is in fact initiated in the micro-scale, through 
the development of micro-cracks/bond breakage, changing the micro-structure and 
therefore the material properties. As the destructive loads keep being applied 
to the system, damage propagates to macro-scales and eventually leads to 
failure. This change of material properties due to damage makes dynamics of the 
system to vary in response to similar excitation. For instance, the stiffness 
of a damaged cantilever beam made of a ductile material is less than an intact 
one\cite{chen1996ductile}. By applying a similar base excitation, it is 
observed that the frequency response of the two structures are distinct. 

Structural health monitoring (SHM) lowers the rates of catastrophic failure, 
increases the safety and reduces the overall cost of system maintenance of a 
mechanical structure by early damage detection. Due to their easy 
implementation, response-based damage detections are one of the most common 
practice in SHM.  This damage detection method uses the information  such as 
natural frequency and mode shape that dynamic response of structure provides to 
predict the existence, location and progress of damage in a 
structure\cite{wang2007improved}.
Damage identification of structural systems such as cantilever beam can be done by tracking the changes in it's vibration characteristics, since damage decreases stiffness of the structure.

Several experimental research have been done to identify damage in a beam using different methods \cite{khatir2018crack,zhou2017cosine,zhou2017structural}. 
A Cantilever beam is a structural element with various engineering applications
which is subjected to dynamical loading. SHM using vibration analysis have five 
levels: detection, localization, classification, assessment, and prediction. 
There are plenty of research focus on most of these levels. SHM for 
vibration-based phenomenon requires high accuracy data to examine the health of 
system. This data can be obtained from modal parameters information based on 
the structural response measurement \cite{khatir2018crack}. 
Here, we try to build the bridge between system response due to specific known 
inputs, and the health of the system in a sense to predict the life span of the 
system. We accomplish this by estimating system properties from its response to 
inputs, and then, finding a relation between damage and those properties. 

Most methods to identify the vibration-based damage use linear models \cite{doebling1996damage,staszewski2004health}. Since structural response is usually nonlinear, capturing the nonlinearity via nonlinear models is promising. Investigating the vibration of geometrically non-linear viscoelastic cantilever beam requires developing compliant models that well-describe the behavior of viscoelastic materials and lend themselves to (nonlinear) numerical methods where the exact solutions are not analytically available.\\

\noindent \textbf{Viscoelasticity.} Many experimental observations in the 
literature show viscoelastic behavior of material in different 
environmental/boundary conditions, meaning that they do not behave purely 
elastic and there exists some internal dissipation mechanism. Viscoelastic 
materials have both a notion of stored and dissipative energy component. This 
common characteristic reveals the material stiffness and damping properties as 
a function of temperature and frequency \cite{torvik1984appearance}. The 
classical models such as Maxwell and Kelvin-Voigt are considered as a 
combination of spring and dash-pot \cite{pipkin2012lectures, 
	christensen2012theory}. Although these models provide accurate results for 
exponential viscoelastic behavior with a finite number of relaxation times, 
they merely truncate the power-law physics of a broad class of anomalous 
materials, providing satisfactory representations only for short observation 
time \cite{Jaishankar2013}. Nutting and 
Gemant separately show in their work that power law function with real order 
exponent is more descriptive for creep or relaxation 
\cite{nutting1921new,gemant1936method}. Later Bagley and Torvik 
\cite{bagley1983theoretical}, show the agreement of fractional calculus models 
of viscoelastic materials with the molecular theories describing the behavior 
of viscoelastic materials. Considering the power law, we can use Riemann 
Liouville fractional integral for constitutive law relating deformation-stress, 
and Caputo fractional derivative for stress-strain constitutive law 
\cite{di2013fractional, suzuki2016fractional}. Several works investigated such 
modeling for bio engineering 
\cite{naghibolhosseini2015estimation, naghibolhosseini2018fractional, 
	perdikaris2014fractional}, visco-elasto-plastic modeling for power law 
dependent stress-strain \cite{Suzuki2017Thesis} and more 
\cite{shitikova2017interaction,rossikhin1997applications}. In the present 
paper, we consider a general case of fractional constitutive laws through 
distributed order differential equations (DODEs), where employing specific 
material distributions functions, we recover a fractional 
Kelvin-Voigt viscoelastic constitutive model to study time dependent frequency 
response of the system.\\

\noindent \textbf{Numerical methods for fractional differential equations.} 
The advancement of numerical methods to solve fractional ordinary differential 
equations (FODEs) and fractional partial differential equations (FPDEs), 
increases the number of research using fractional modeling. Spectral methods 
for spatial discretization \cite{azrar2002semi,zayernouri2014exponentially, 
	zayernouri2014spectral, zayernouri2015fractionalcollocation, 
	samiee2019unified1, samiee2019unified2} and distributed order differential 
equations\cite{samiee2018petrov} are two classes of these numerical methods. 
One of the numerical scheme to solve nonlinear FPDEs is finite element 
method\cite{rao2017finite}, but the downside of this scheme is being 
computationally expensive. Even using high performance computing, direct 
integration of non-linear equations is not the most suitable method. Therefore, 
the problem should be reduced by representing the unknowns as a linear 
combination of several well defined functions assuming that the unknowns (in 
our case displacement) can be represented by a linear combination of several 
well-defined functions. This method reduce the computational cost and have a 
reasonable accuracy. The technique of modal analysis in linear vibration is an 
appropriate and classic procedure. The system of differential equations 
obtained from eigenmodes can be solved separately. For non-linear vibration 
Galerkin method due to its higher-order accuracy\cite{li2009space}, is a 
suitable choice to approximate the governing partial differential equation by a 
system of ordinary differential equations obtained by representing the solution 
as a sum of independent temporal functions, each of which satisfies the spatial 
boundary conditions and is multiplied by a time dependent 
coefficient\cite{qian1996anomalous}. This is a standard method for deriving the 
ordinary differential equations used in studies of chaotic vibrations of 
nonlinear elastic systems\cite{qian1996anomalous}. In the present work, our 
focus is to construct a spectral Galerkin approximation in space to solve the 
fractional partial differential equation of motion. Spectral method with mode 
shapes of fractional beam is used as basis to discretize the space. To 
approximate the time derivative, we use L1 scheme. Among a class of numerical 
methods for time-fractional differential equations \cite{ 
	lubich1986discretized, zayernouri2015tempered, 
	suzuki2018automated, zayernouri2016fractionalAdams}, in this work we employ 
the direct L1 scheme developed by Lin and Xu \cite{lin2007finite}.

This work is organized as follows. In section \ref{Sec: Mathematical 
	Formulation}, we derive the governing equation for the nonlinear in-plane 
vibration of a visco-elastic cantilever beam. Afterwards, we use the extended 
Hamilton's principle to derive the equation of motion when the only external 
force is the base excitation. Then, we obtain the weak formulation of the 
problem and use assumed modes in space to reduce the problem to system of 
differential equations in time. In section  \ref{Sec: Direct Numerical Time 
	Integration}, we obtain the corresponding linearized equation of motion. We 
perform the perturbation analysis in section \ref{Sec: Perturbatin Analysis} to 
solve the resulting nonlinear FODE. We report the super sensitivity of 
response amplitude to the fractional model parameters at free vibration, and 
bifurcation in steady-state amplitude at primary resonance and finally conclude 
the paper with a summary.

\section{Mathematical Formulation}
\label{Sec: Mathematical Formulation}
%%%%%%%%%%%%%%%%%%%%%%%%%%%
%
We formulate the mathematical model that describes the behavior of the considered physical system, and discuss the main assumptions and theorems, used to derive the equation of motion. We employ spectral decomposition in space to discretize the problem and further use the perturbation method to solve the resulting nonlinear equations.

%
%%%%%%%%%%%%%%%%%%%%%%%%%%%
%\subsection{In-Plane Vibration of a Visco-Elastic Euler-Bernoulli Cantilever Beam with Nonlinear Geometry}
\subsection{Nonlinear In-Plane Vibration of a Visco-Elastic Cantilever Beam}
\label{Sec: Fractional NL Beam}
%%%%%%%%%%%%%%%%%%%%%%%%%%%
%

\begin{figure}[h]
	\centering
	\includegraphics[width=1\linewidth]{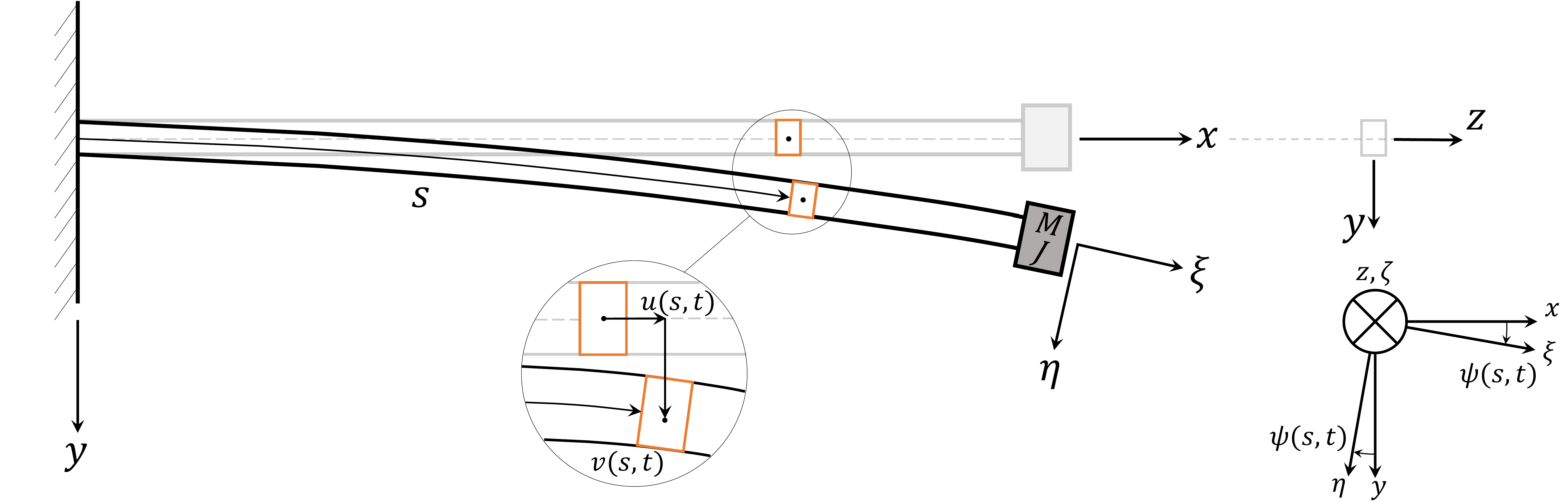}
	\caption{In-plane lateral deformation of a slender isotropic cantilever 
		beam. The terms $u(s,t)$ and $v(s,t)$ denote, respectively, the axial and 
		lateral 
		displacements, and $\psi(s,t)$ is the rotation angle about $z$ axis.}
	\label{Fig: Cantilever Beam}
\end{figure}

\begin{SCfigure}
	\centering
	\includegraphics[width=0.5\linewidth]{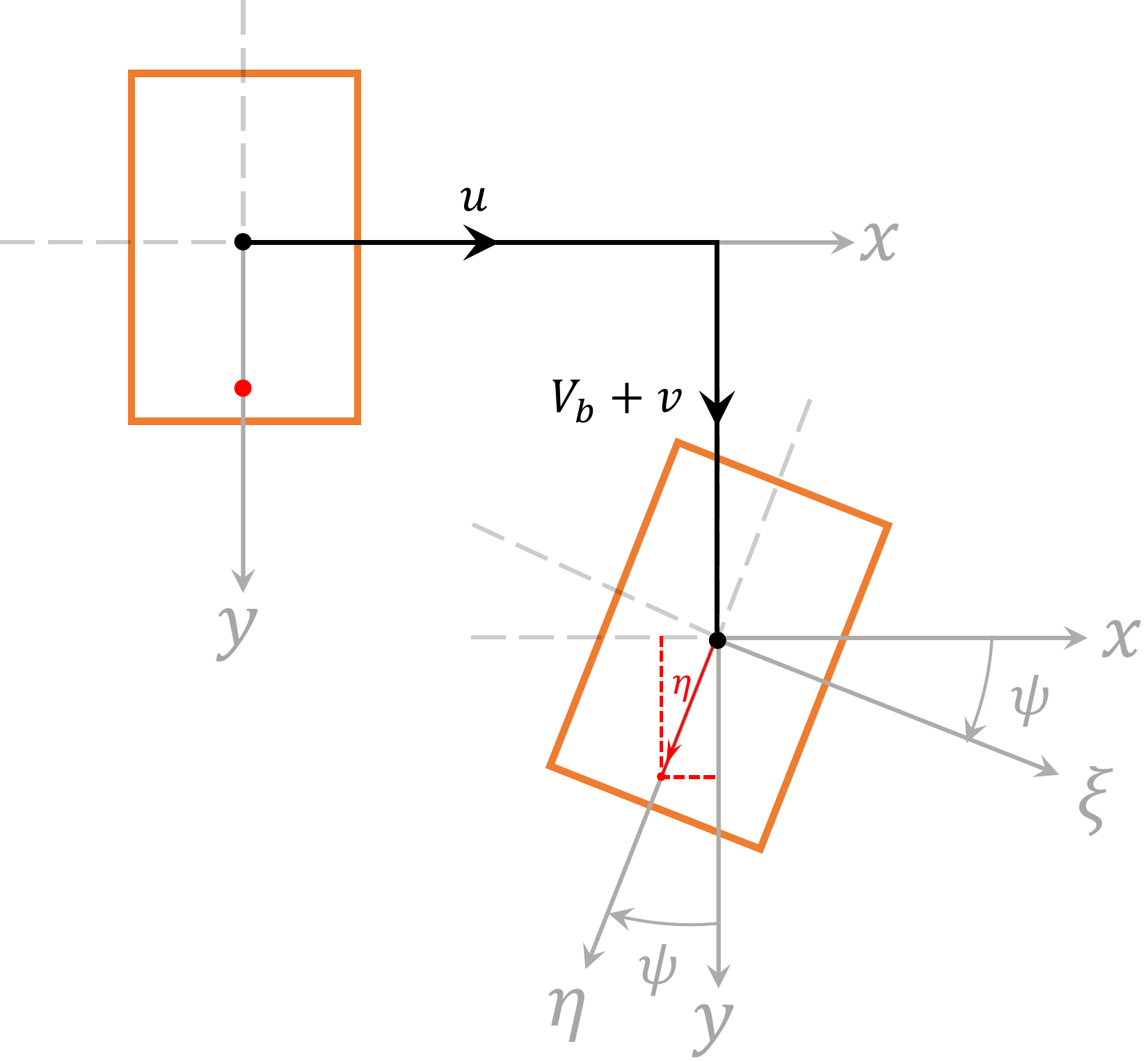}
	\caption{Detailed in-plane lateral deformation of a slender isotropic 
		cantilever beam. The total deformation of an arbitrary point (shown in red 
		color) as the beam undergoes a deformation process is illustrated. This 
		deformation 
		is comprised of the axial displacement of the beam $u$, the lateral 
		displacement of beam in addition to the base motion $v+V_b$, and the 
		displacement due to rotation $\psi$.}
	\label{Fig: Cantilever Beam deformation}
\end{SCfigure}

We consider the nonlinear response of a slender isotropic visco-elastic 
cantilever beam with lumped mass $M$ at the tip, subject to harmonic transverse 
base excitation, $V_b$ \textit{(see Figures \ref{Fig: Cantilever Beam} and 
	\ref{Fig: Cantilever Beam deformation})}. We use the nonlinear Euler-Bernoulli 
beam theory to obtain the governing equations, where the geometric 
nonlinearities in a cantilever beam with symmetric cross section is included in 
the equations of motion. We consider the following list of kinematic and 
geometric assumptions for the beam and derive the corresponding governing 
equations. 
\begin{itemize}[leftmargin=8mm]
	\item The beam is idealized as an inextensional one, i.e., stretching of 
	the neutral axis is insignificant. The effects of warping and shear 
	deformation are ignored. Therefore, the strains acting in the cross section 
	are only due to the bending kinematics.
	
	\item The beam is slender with symmetrical cross section undergoes purely planar flexural vibration.
	
	\item The length $L$, cross section area $A$, mass per unit length $\rho$, 
	density $\varrho$, mass $M$ and rotatory inertia $J$ of the lumped mass at 
	the tip of beam are constant.
	
	\item The axial displacement along length of beam and the lateral 
	displacement are respectively denoted by $u(s,t)$ and $v(s,t)$.
	
	\item We consider the in-plane transverse vibration of the beam and reduce the problem to 1-dimension.
	
\end{itemize}
Figures \ref{Fig: Cantilever Beam} and \ref{Fig: Cantilever Beam deformation} show the lateral deformation of the considered cantilever beam .
%We consider the in plane lateral vibration of a 1-dimensional Euler-Bernoulli cantilever beam, shown in Fig. \ref{Fig: Cantilever Beam}. The displacement along length of beam and lateral displacement are denoted by $u(s,t)$ and $v(s,t)$, respectively. 
As the beam deforms, we let the inertial coordinate system $(x,y,z)$ rotate 
about the $z$ axis with rotation angle $\psi(s,t)$ to the coordinate system 
$(\xi,\eta,\zeta)$, where
\begin{align*}
\begin{bmatrix} \textbf{e}_{\xi} \\ \textbf{e}_{\eta} \\ \textbf{e}_{\zeta} \end{bmatrix}
= \begin{pmatrix}cos(\psi) & sin(\psi) & 0 \\ -sin(\psi) & cos(\psi) & 0 \\ 0 & 0 & 1\end{pmatrix}
\begin{bmatrix} \textbf{e}_{x} \\ \textbf{e}_{y} \\ \textbf{e}_{z} \end{bmatrix},
\end{align*}
and $\textbf{e}_{i}$ is the unit vector along the $i^{\text{\,th}}$ coordinate. 
The angular velocity and curvature of the beam at any point along the length of 
the beam $s$ and any time $t$ can be written, respectively, as
\begin{align}
\label{Eq: angular velocity curvature}
\boldsymbol{\omega}(s,t) = \frac {\partial{\psi}}{\partial{t}} \, \textbf{e}_{z}, \quad 
\boldsymbol{\rho}(s,t) =  \frac {\partial{\psi}}{\partial{s}} \, \textbf{e}_{z}
\end{align}
%++++++++++++++++++++++++++++++++++++++++++
%++++++++++++++++++++++++++++++++++++++++++
%

\noindent The total displacement and velocity of an arbitrary point along the $y$ axis takes the form:
\begin{align}
\label{Eq: displacement}
\textbf{r} & = \left (u - \eta \,  \sin(\psi) \right) \, \textbf{e}_x + \left (v + V_b + \eta \, \cos(\psi)\right) \, \textbf{e}_y ,
\\ \label{Eq: velocity}
\frac {\partial{\textbf{r}}}{\partial{t}} & = \left( \frac {\partial{u}}{\partial{t}} - \eta \,  \frac {\partial{\psi}}{\partial{t}} \, \cos(\psi) \right) \, \textbf{e}_x + \left ( \frac {\partial{v}}{\partial{t}}  +  \frac {\partial{V_b}}{\partial{t}} - \eta \,  \frac {\partial{\psi}}{\partial{t}} \, \sin(\psi)\right) \, \textbf{e}_y .
\end{align}
%

%
%++++++++++++++++++++++++++++++++++++++++++
%++++++++++++++++++++++++++++++++++++++++++
%

\begin{figure}[t]
	\centering
	\includegraphics[width=0.75\linewidth]{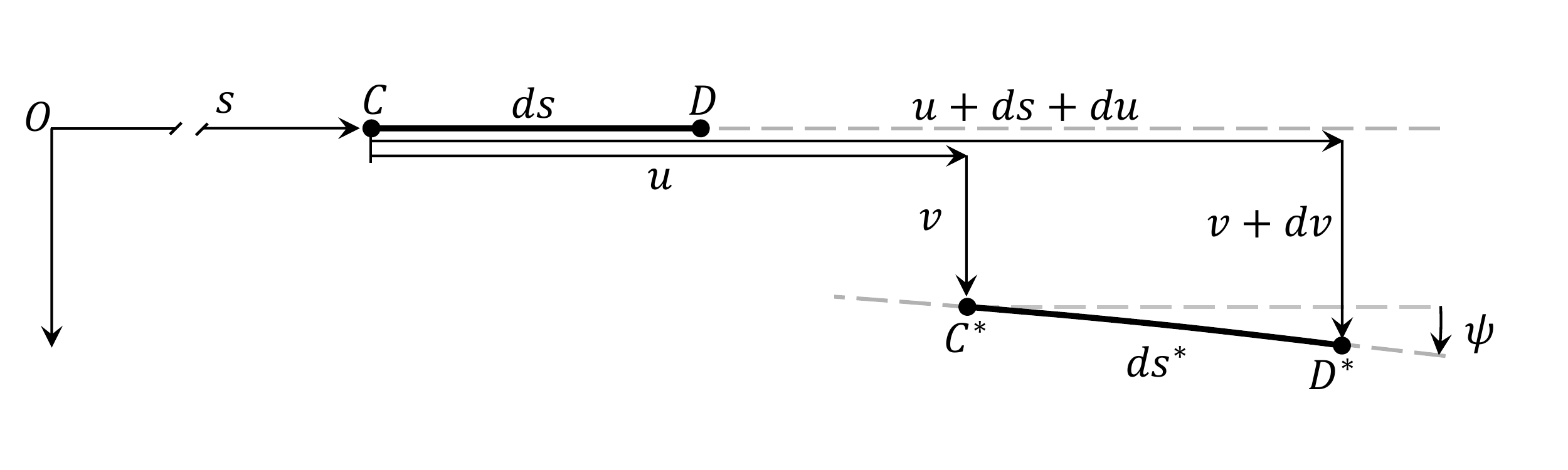}
	\caption{Deformation of an arbitrary element of the beam. The initial 
		configuration $CD$ extends, traverses, and rotates to another configuration 
		$C^*D^*$.}
	\label{Fig: inextensibility}
\end{figure}

\noindent We also let an arbitrary element $CD$ of the beam neutral axis, which 
is of length $ds$ and located at a distance $s$ from the origin $O$, to undergo 
a deformation to an arbitrary configuration $C^* D^*$ \textit{(see Fig. 
	\ref{Fig: inextensibility})}. The displacement components of points $C$ and $D$ 
are denoted by the pairs $(u,v)$ and $(u+du,v+dv)$, respectively. The axial 
strain 
$e(s,t)$ at 
the arbitrary point $C$ is then given by
\begin{align}
\label{Eq: strain}
e=\frac{ds^*-ds}{ds}=\frac{\sqrt{(ds+du)^2+dv^2} - ds}{ds} = \sqrt{(1+  \frac {\partial{u}}{\partial{s}})^2+ ( \frac {\partial{v}}{\partial{s}})^2 } - 1.
\end{align}
Applying the inextensionality constraint, i.e. $e = 0$, (\ref{Eq: strain}) 
becomes
\begin{align}
\label{Eq: inextensionality}
1 +  \frac {\partial{u}}{\partial{s}} = \left(1 - ( \frac {\partial{v}}{\partial{s}})^2\right)^{1/2}.
%u^{\prime} = (1 - {\frac{\partial{v}}{\partial{s}}}^2)^{1/2} - 1 \simeq 1 - \frac{1}{2} {\frac{\partial{v}}{\partial{s}}}^2 - 1 \simeq - \frac{1}{2} {\frac{\partial{v}}{\partial{s}}}^2.
%
\end{align}
Moreover, based on the assumption of no transverse shear deformation and using \eqref{Eq: inextensionality}, we have 
\begin{align}
\label{Eq: no-shear}
\psi = \text{tan}^{-1}\frac{ \frac {\partial{v}}{\partial{s}}}{1+ \frac {\partial{u}}{\partial{s}}} = \text{tan}^{-1}\frac{ \frac {\partial{v}}{\partial{s}}}{\left(1 - ( \frac {\partial{v}}{\partial{s}})^2\right)^{1/2}} .
\end{align}
Using the expansion $\tan^{-1}(x) = x -\frac{1}{3}x^3 + \cdots$, the curvature can be approximated up to cubic term as
\begin{align}
\label{Eq: rotation}
\psi 
& =  \frac {\partial{v}}{\partial{s}}(1 - ( \frac {\partial{v}}{\partial{s}})^2)^{-1/2} - \frac{1}{3} ( \frac {\partial{v}}{\partial{s}})^3(1 - ( \frac {\partial{v}}{\partial{s}})^2)^{-3/2} + \cdots
\\ \nonumber
& \simeq  \frac {\partial{v}}{\partial{s}}(1 + \frac{1}{2} ( \frac {\partial{v}}{\partial{s}})^2) - \frac{1}{3} ( \frac {\partial{v}}{\partial{s}})^3
\simeq  \frac {\partial{v}}{\partial{s}}+ \frac{1}{6} ( \frac {\partial{v}}{\partial{s}})^3
\end{align}
Therefore, the angular velocity and curvature of the beam, i.e. $ \frac {\partial{\psi}}{\partial{t}}$ and $ \frac {\partial{\psi}}{\partial{s}}$, respectively, can be approximated as:
\begin{align}
\label{Eq: angular velocity}
& \frac {\partial{\psi}}{\partial{t}} \simeq \frac {\partial^2{v}}{\partial{t}\partial{s}} + \frac{1}{2} \, \frac {\partial^2{v}}{\partial{t}\partial{s}} \, (\frac {\partial{v}}{\partial{s}} )^2 \simeq  \frac {\partial^2{v}}{\partial{t}\partial{s}}( 1 + \frac{1}{2} (\frac {\partial{v}}{\partial{s}})^2) ,
\\ \label{Eq: curvature}
&\frac {\partial{\psi}}{\partial{s}}   \simeq \frac {\partial^2{v}}{\partial{s}^2} + \frac{1}{2}\frac {\partial^2{v}}{\partial{s}^2} (\frac {\partial{v}}{\partial{s}})^2 \simeq \frac {\partial^2{v}}{\partial{s}^2}(1 + \frac{1}{2} (\frac {\partial{v}}{\partial{s}})^2) .
\end{align}
By the Euler-Bernoulli beam assumptions a slender, no-transverse-shear with no strains in the plane of cross sectional plane, the strain-curvature relation takes the form
\begin{align}
\label{Eq: strain-curvature relation}
\varepsilon(s,t) = 
%-\eta \, \rho(s,t) = 
-\eta \,\frac {\partial{\psi}(s,t)}{\partial{s}}
\end{align}
%

%Using \eqref{Eq: inextensionality} and the expansion $\text{tan}^{-1}(x) = x -\frac{1}{3} x^3 + \cdots$, we have
%%
%\begin{align}
%\label{Eq: no-shear}
%%
%\psi = \text{tan}^{-1}\frac{\frac{\partial{v}}{\partial{s}}}{(1 - {\frac{\partial{v}}{\partial{s}}}^2)^{1/2}} \simeq \frac{\partial{v}}{\partial{s}} (1 - {\frac{\partial{v}}{\partial{s}}}^2)^{1/2} - \frac{1}{3} {\frac{\partial{v}}{\partial{s}}}^{3} (1 - {\frac{\partial{v}}{\partial{s}}}^2)^{-3/2}
%%
%\end{align}
%%

%++++++++++++++++++++++++++++++++++++++++++
%++++++++++++++++++++++++++++++++++++++++++

%\vspace{1 in}
%
%
%, $v(x,t) : \Omega \rightarrow \mathbb{R}$, of a 1-dimensional Euler-Bernoulli cantilever beam, shown in Fig. \ref{Fig: Cantilever Beam}, where the computational domain $\Omega = [0,L] \times [0,T]$. %A Cartesian coordinate is attached to the beam base, where $z$ is along its span and the downward lateral coordinate $y$ is in positive sense. 
%We also let $A$, $\rho$, and $f$ to be cross section area, mass per unit length, and external force per unit length, respectively.The governing equations can be obtained as
%%
%\begin{align}
%\label{Eq: Euler Bernouli Beam}
%&\rho \frac{\partial^2 v }{\partial t^2} - \frac{\partial^2 M }{\partial x^2} = f,
%\\
%\label{Eq: Moment}
%&M = \int_{A} y \sigma dA,
%\\
%\label{Eq: Starin}
%&\varepsilon = -y \frac{\partial^2 v}{\partial x^2}
%%
%\end{align}
%%

%\newpage

%
%%%%%%%%%%%%%%%%%%%%%%%%%%%
\subsection{Linear Viscoelasticity: Boltzmann Superposition Principle}
\label{Sec: Viscoelasticity}
%%%%%%%%%%%%%%%%%%%%%%%%%%%
%
In this section we introduce the fractional-order Kelvin-Voigt model employed 
in this work. We start with a \textit{bottom-up} derivation of our rheological 
building block, \textit{i.e.}, the Scott-Blair model through the Boltzmann 
superposition principle. Then, in a \textit{top-bottom} fashion, we demonstrate 
how the fractional Kelvin-Voigt model is obtained from a general 
distributed-order form. 

Many experimental observations in the literature show viscoelastic behavior of material in different environmental/
boundary conditions, meaning that they do not behave purely elastic and there 
exists some internal dissipation mechanism. In such cases, the resulting stress 
has a memory dependence on the velocity of all earlier deformations, which can 
be described by the Boltzmann superposition principle. When the specimen is 
under loading, the material instantaneously reacts elastically and then, 
immediately starts to relax; this is where dissipation takes place. Thus, as a 
step increase in elongation (from the stretch $\lambda = 1$ to some $\lambda$) 
is imposed, the developed stress in the material will be a function of time and 
the stretch: 
\begin{align}
\label{Eq: relaxation function}
K(\lambda , t ) = G(t) \, \sigma^{(e)}(\lambda),
\end{align}
where $G(t)$ is the \emph{reduced relaxation function} and $\sigma^{(e)}$ is the \emph{elastic response} (in absence of any viscosity). $\sigma^{(e)}$ can also be interpreted as tensile stress response in a sufficiently high rate loading experiment. The Boltzmann superposition principle states that the stresses from different small deformations are additive, meaning that the total tensile stress of the specimen at time $t$ is obtained from the superposition of infinitesimal changes in stretch at some prior time $\tau_j$, given as $G(t-\tau_j) \frac{\partial \sigma^{(e)}[\lambda(\tau_j)]}{\partial \lambda} \delta \lambda(\tau_j)$. Therefore,
\begin{align}
\label{Eq: tensile stress - 1}
\sigma ( t ) 
= \sum_{\tau_j<t} G(t-\tau_j) \frac{\partial \sigma^{(e)}[\lambda(\tau_j)]}{\partial \lambda} \frac{\delta \lambda(\tau_j)}{\delta \tau_j} \delta \tau_j,
\end{align}
where in the limiting case $\delta \tau_j \rightarrow 0$ gives the integral form of the equation as
\begin{align}
\label{Eq: tensile stress - 2}
\sigma ( t ) 
=\int_{-\infty}^{t} G(t - \tau) \, \frac{\partial \sigma^{(e)}[\lambda(\tau)]}{\partial \lambda} \, \frac{\partial \lambda}{\partial \tau} \, d\tau
= \int_{-\infty}^{t} G(t - \tau) \, (\frac {\partial{\sigma}}{\partial{t}})^{(e)} \, d\tau.
\end{align}

\begin{rem}
	We note that \ref{Eq: tensile stress - 2} introduces a mathematical convolution type integro-differential operator as the constitutive equation that describes the stress-strain relation. The kernel of operator $G(t)$ is the relaxation function, which inherently is the mechanical property of the material under consideration and in general is obtained from experimental observations.
\end{rem}

\begin{figure}[t]
	\centering
	%	\begin{subfigure}{0.45\textwidth}
	%	\centering
	%	\includegraphics[width=1\linewidth]{figures/MaxwellModel.pdf}
	%	\end{subfigure}
	%	\begin{subfigure}{0.45\textwidth}
	%		\centering
	%		\includegraphics[width=1\linewidth]{figures/KLModel.pdf}
	%	\end{subfigure}
	\includegraphics[width=0.57\linewidth]{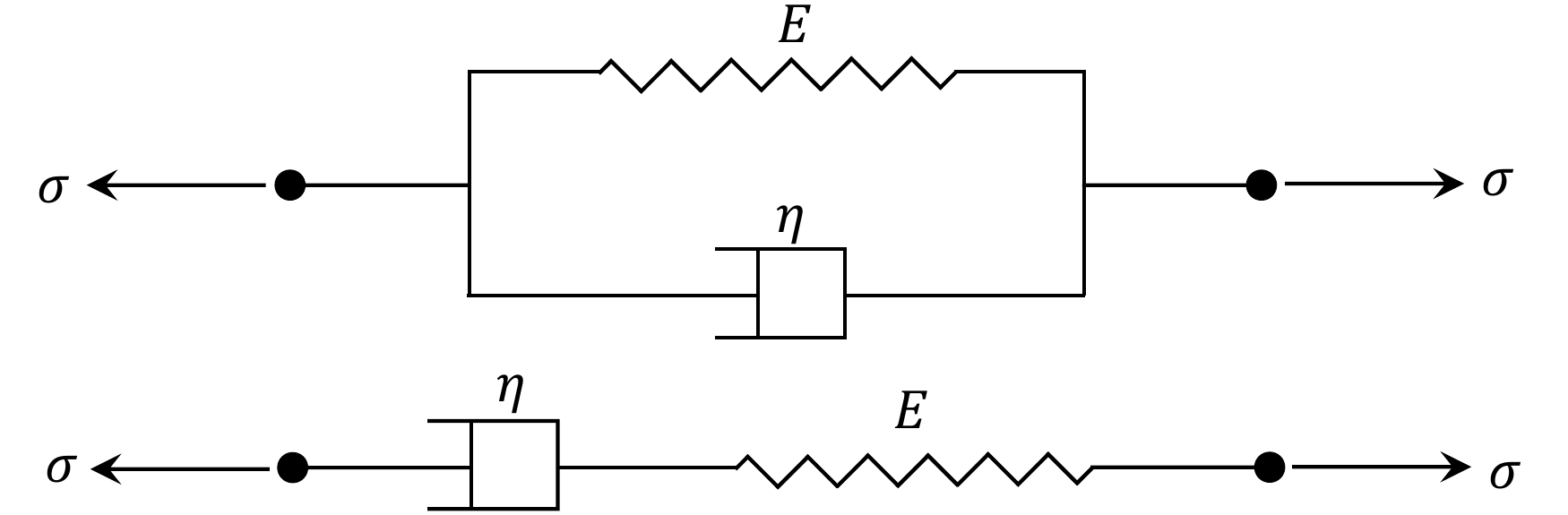}
	\caption{Classical visco-elastic models as a combination of spring (purely elastic) and dash-pot (purely viscous) elements. Kelvin-Voigt (top) and Maxwell (bottom) rheological models.}
	\label{Fig: Maxwell}
\end{figure}

\vspace{0.2 in}
%%%=============================================================================%%
\noindent \textbf{Exponential Relaxation, Classical Models:}
The relaxation function $G(t)$ is traditionally expressed as the summation of 
exponential functions with different exponents and constants, which yields the 
\textit{so-called} generalized Maxwell form as: 
\begin{align}
\label{Eq: exponential relaxation}
G(t) = \frac{\sum C_i e^{- t/\tau_i}}{\sum C_i}.
\end{align}
For the simple case of a single exponential term (a single Maxwell branch), we 
have $G(t) = e^{- t/\tau}$. Therefore, in the case of zero initial strain 
$(\varepsilon(0) = 0)$, we have:
\begin{align}
\sigma ( t ) 
= \int_{0}^{t}  e^{- (t - \tilde{t})/\tau} \, E \frac {\partial{\varepsilon}}{\partial{t}} \, d\tilde{t},
\end{align}
which solves the integer-order differential equation $ 
\frac{\partial{\varepsilon}}{\partial{t}} = \frac{1}{E}  
\frac{\partial{\sigma}}{\partial{t}} + \frac{1}{\eta} \sigma$, where the 
relaxation time constant $\tau = \eta / E$ is obtained from experimental 
observations. The Maxwell model is in fact a combination of purely elastic and 
purely viscous elements in series, see Fig. \ref{Fig: Maxwell}. Other different 
combinations of these purely elastic/viscous elements in both series and/or 
parallel give rise to various rheological models with distinctive properties, 
each of which can be used to model different types of material 
\cite{pipkin2012lectures, christensen2012theory}. However, one of the key 
issues in such modeling is that they require complex mechanical arrangements of 
a large number of Hookean spring and Newtonian dashpots in order to adequately 
model the complex hereditary behavior of power-law materials. They generally 
cannot fully capture such behavior as the standard building blocks do not 
reflect long-memory dependence in the material response. More importantly, 
these models introduce a relatively large number of model parameters, which 
adversely affect the condition of ill-posed inverse problem of parameter 
estimation and model fitting \cite{aster2018parameter}.

%By different combinations of these two building blocks in both series and parallel, one can obtain various material models with distinctive properties. Yet, in order to capture complex behavior of materials, one need to design a complicated combinations of elastic and viscous element. While this inherently do not reflect any history dependencies of material response, the excessive number of model parameters require exploiting expensive scheme to fit data.  

%
%++++++++++++++++++++++++++++++++++++++++++
%++++++++++++++++++++++++++++++++++++++++++
%

\vspace{0.2 in}
%%%=============================================================================%%
\noindent  \textbf{Power-Law Relaxation, Fractional Models:}
The mechanical stress appeared at the deformation of viscoelastic materials 
decreases as power-law functions in time \cite{bagley1989power}, suggesting 
that relaxation of stress obeys a power law behavior rather than exponential, 
and thus, the relaxation time can not be described with single time scale 
anymore \cite{mainardi2011creep}. By letting the relaxation function (kernel) 
in \eqref{Eq: tensile stress - 2} have a modulated power-law form $G(t) = 
g(\alpha) ( \, t-\tau \, )^{-\alpha}$ and elastic part response $\sigma^{(e)} = 
E \varepsilon $, the tensile stress takes the form of 
\begin{align}
\label{Eq: fractional stress}
\sigma ( t ) 
= \int_{-\infty}^{t}  \frac{g(\alpha)}{( \, t-\tau \, )^{\alpha}} \, E  \frac{\partial{\varepsilon}}{\partial{t}} \, d\tau
= E \, g(\alpha) \,  \int_{-\infty}^{t}  \frac{\frac{\partial{\varepsilon}}{\partial{t}}}{( \, t-\tau \, )^{\alpha}} \, d\tau.
\end{align}
If we choose the modulation $g(\alpha) = \frac{1}{\Gamma(1-\alpha)}$, then the 
integro-differential operator \eqref{Eq: fractional stress} gives the 
Liouville-Weyl fractional derivative \cite{mainardi2010fractional}. Although 
the lower integration limit of (\ref{Eq: fractional stress}) is taken as 
$-\infty$, under hypothesis of causal histories, which states that the 
viscoelastic body is quiescent for all time prior to some starting point $t=0$, 
\eqref{Eq: fractional stress} can be re-written as 
\begin{align}
\label{Eq: fractional stress - 2}
\sigma ( t ) 
& = \varepsilon(0^{+}) \, \frac{E \, g(\alpha)}{t ^{\alpha}} + E \, g(\alpha) \,  \int_{0}^{t}  \frac{\frac{\partial{\varepsilon}}{\partial{t}}}{( \, t-\tau \, )^{\alpha}} \, d\tau,
\\ \nonumber
& = \varepsilon(0^{+}) \, \frac{E \, g(\alpha)}{t ^{\alpha}} + E \, \prescript{C}{0}{\mathcal D}_{t}^{\alpha} \, \varepsilon,
\\ \nonumber
& = E\, \prescript{RL}{0}{\mathcal D}_{t}^{\alpha} \, \varepsilon,
\end{align}
where $\prescript{C}{0}{\mathcal D}_{t}^{\alpha}$ and 
$\prescript{RL}{0}{\mathcal D}_{t}^{\alpha}$ denote, respectively, the Caputo 
and Riemann-Liouville fractional derivatives \cite{mainardi2010fractional}. 
Both definitions are equivalent here due to homogeneous initial conditions for 
the strain. 

\begin{rem}
	\label{Rem: Scott Blair}
	The constitutive equation \eqref{Eq: fractional stress - 2} can be thought of as an interpolation between a pure elastic (spring) and a pure viscous (dash-pot) elements, i.e., the Scott Blair element \cite{rogosin2014george,mainardi2008time,mainardi2011creep,suzuki2016fractional}. It should be noted that in the limiting cases of $\alpha \rightarrow 0$ and $\alpha \rightarrow 1$, the relation \eqref{Eq: fractional stress - 2} recovers the corresponding equations for spring and dash-pot, respectively.
\end{rem}

\vspace{0.2 in}
%%%=============================================================================%%
\noindent  \textbf{Multi-Scale Power-Laws, Distributed-Order Models:}
In the most general sense, materials intrinsically possess a spectrum of 
power-law relaxations, and therefore we need a distributed-order representation 
for the stress-strain relationship. Consequently, the relaxation function 
$G(t)$ in \eqref{Eq: tensile stress - 2} does not only a single power-law 
as in \eqref{Eq: fractional stress}, but rather contains a distribution over a 
range of 
values. The use of distributed-order models to generalize the stress-strain 
relation of inelastic media and the Fick’s law was proposed in 
\cite{caputo1995m,caputo2001m}. The connection of such operators with 
diffusion-like equations was established in \cite{chechkin2002, 
	chechkin2008generalized}. The application of distributed-oder models are also 
discussed in rheological models 
\cite{lorenzo2002variable,Atanackovic2009distributional}, economic processes 
with distributed memory fading \cite{tarasova2018concept}, continuous time 
random walk \cite{jiao2012distributed,sandev2015distributed}, time domain 
analysis of control, filtering and signal processing 
\cite{li2011distributed,li2014lyapunov}, vibration \cite{duan2018steady}, 
frequency domain analysis \cite{bagley2000}, and uncertainty quantification 
\cite{kharazmi2017petrov, kharazmi2018fractional}. Considering nonlinear 
visco-elasticity with material heterogeneities, the distributed order 
constitutive equations over $t > 0$ with orders $\alpha \in [\alpha_{min}, 
\alpha_{max}]$ and $\beta \in [\beta_{min}, \beta_{max}]$ can be expressed in 
the general form as
\begin{align}
\label{Eq: Dis const law}
\int_{\beta_{min}}^{\beta_{max}} \, \Phi(\beta; x, t, \sigma) \, 
\prescript{*}{0}{\mathcal D}_{t}^{\beta} \sigma(t) \, d\beta 
=
\int_{\alpha_{min}}^{\alpha_{max}} \, \Psi(\alpha; x, t, \varepsilon) \, 
\prescript{*}{0}{\mathcal D}_{t}^{\alpha} \varepsilon(t) \, d\alpha ,
\end{align}
in which the prescript $*$ stands for any type of fractional derivative. The 
functions $\Phi(\beta; x, t, \sigma)$ and $\Psi(\alpha; x, t, \varepsilon)$ can 
be thought of as distribution functions, where $\alpha \mapsto \Psi(\alpha; x, 
t, \varepsilon)$ 
and $\beta \mapsto \Phi(\beta; x, t, \sigma)$  are continuous mappings in 
$[\alpha_{min} , 
\alpha_{max}]$ and $[\beta_{min} , \beta_{max}]$. Furthermore, the dependence 
of the distributions on the (thermodynamically) conjugate pair $(\sigma, 
\varepsilon)$ introduces the notion of nonlinear viscoelasticity, and the 
dependence on a material coordinate $x$ induces material heterogeneties in 
space.

\begin{rem}
	The pairs ($\alpha_{min}$, $\alpha_{max}$) and ($\beta_{min}$, $\beta_{max}$) 
	are only the theoretical lower and upper terminals in the definition of 
	distributed order models. In general, the distribution function $\Phi(\beta; x, 
	t, \sigma)$ and $\Psi(\alpha; x, t, \varepsilon)$ can arbitrarily confine the 
	domain of integration in each realization of practical rheological problems and 
	material design. If we let the distribution be summation of some delta 
	functions, then, the distributed order model becomes the following multi-term 
	model: 
	\begin{align*}
	\label{Eq: frac const law}
	\left( 1+ \sum_{k=1}^{p_{\sigma}} \, a_k \, \prescript{}{0}{\mathcal D}_{t}^{\beta_k} \right) \, \sigma(t)
	=
	\left( c + \sum_{k=1}^{p_{\varepsilon}} \, b_k \, \prescript{}{0}{\mathcal D}_{t}^{\alpha_k} \right) \, \varepsilon(t).
	\end{align*}
\end{rem}

In order to obtain the fractional Kelvin-Voigt model, we let $\Phi(\beta) = 
\delta(\beta)$ and $\Psi(\alpha) = E_{\infty} \delta(\alpha) + E_{\alpha} 
\delta(\alpha-\alpha_0)$ in \eqref{Eq: Dis const law}, and therefore,
\begin{align}
\label{Eq: frac KV}
\sigma(t)
=
E_{\infty} \, \varepsilon(t) + E_{\alpha} \, \prescript{RL}{0}{\mathcal D}_{t}^{\alpha} \, \varepsilon(t),
\end{align}
where $\alpha \in (0,1)$. Since we only have one single derivative order 
$\alpha$, we drop the subscript zero for the sake of simplification.

%
%%%%%%%%%%%%%%%%%%%%%%%%%%%%%%%%%%%%%
\subsection{Extended Hamilton's Principle}
%%%%%%%%%%%%%%%%%%%%%%%%%%%%%%%%%%%%%
%
We derive the equations of motion by employing the extended Hamilton's principle $$\int_{t_1}^{t_2} \left (\delta T - \delta W\right) \, dt = 0,$$ where $\delta T$ and $\delta W$ are the variations of kinetic energy and total work \cite{meirovitch2010fundamentals}. The only source of external input to our system of interest is the base excitation, which superposes base velocity to the beam velocity, and thus contributes to the kinetic energy. Hence, the total work only includes the internal work done by the induced stresses, and its variation can be expressed in the general form as \cite{bonet1997nonlinear}
\begin{align}
\label{Eq: work var}
\delta W = \int_{\mathbb{V}} \sigma \, \delta \varepsilon \, dv,
\end{align}
where the integral is taken over the whole system volume $\mathbb{V}$. The 
stress $\sigma$ includes both the conservative part, $\sigma_{c}$, due to 
elastic and the non-conservative part, $\sigma_{nc}$, due to viscous 
deformation, where the former constitutes the potential energy of the system. 

\begin{rem}
	It is remarked in \ref{Rem: Scott Blair} that the fractional derivative exhibit 
	both elasticity and viscosity. There has been some attempts in the literature 
	to separate the conservative (elastic) and non-conservative (viscous) parts of 
	fractional constitutive equations to define the free energy of the system 
	\cite{lion1997thermodynamics}. We note that as this separation in the time 
	domain is not trivial for sophisticated fractional constitutive equations. 
	Furthermore, as we do not deal with the free energy of our system, we would 
	rather leave the total work not separated and thus do not compute the potential 
	energy and work done by non-conservative forces separately. 
\end{rem}

The full derivation of governing equation using the extended Hamilton's principle is given in Appendix \ref{Sec: App. extended Hamilton}. We recall that $M$ and $J$ are the mass and rotatory inertia of the lumped mass at the tip of beam, $\rho$ is the mass per unit length of the beam, $I = \int_{A} \eta^2 \, dA$, and let $m = \frac{\rho}{ E_{\infty} \, I}$ and $E_r = \frac{E_{\alpha}}{E_{\infty}}$. We approximate the nonlinear terms up to third order and use the following dimensionless variables
\begin{align*}
s^{*} = \frac{s}{L}, \,\, 
v^{*} = \frac{v}{L}, \,\,
t^{*} = t \left(\frac{1}{m L^4}\right)^{1/2}, \,\, 
E_r^{*} = E_r \left(\frac{1}{m L^4}\right)^{\alpha/2}, \,\, 
J^{*} = \frac{J}{\rho L^3}, \,\,
M^{*} = \frac{M}{\rho L}, \,\,
V_b^* =  \frac{V_b}{L},
\end{align*}
and derive the strong form of the equation of motion. Therefore, by choosing proper function space $V$, the problem reads as: find $v \in V$ such that
\begin{align}
\label{Eq: eqn of motion - dimless - body}
\frac{\partial^2{v}}{\partial{t}^2} 
& + 
\frac{\partial^2}{\partial{s}^2} \left(  
\frac{\partial^2{v}}{\partial{s}^2}  + \frac{\partial^2{v}}{\partial{s}^2}  (\frac{\partial{v}}{\partial{s}} )^2 
+ E_r \prescript{RL}{0}{\mathcal D}_{t}^{\alpha} \, \Big[\frac{\partial^2{v}}{\partial{s}^2}(1 + \frac{1}{2} (\frac{\partial{v}}{\partial{s}} )^2)\Big] 
+ \frac{1}{2}  E_r (\frac{\partial{v}}{\partial{s}} )^2 \, \prescript{RL}{0}{\mathcal D}_{t}^{\alpha} \, \frac{\partial^2{v}}{\partial{s}^2} 
\right)
\\ \nonumber
& - 
\frac{\partial}{\partial{s}}\left( 
\frac{\partial{v}}{\partial{s}} \, (\frac{\partial^2{v}}{\partial{s}^2})^2  
+ E_r \frac{\partial{v}}{\partial{s}} \, \frac{\partial^2{v}}{\partial{s}^2} \,  \prescript{RL}{0}{\mathcal D}_{t}^{\alpha} \,  \frac{\partial^2{v}}{\partial{s}^2}
\right)
= - \ddot{V_b},
\end{align}
which is subject to the boundary conditions:
\begin{align}
\label{Eq: bc - dimless - body}
& v \, \Big|_{s=0} = \frac{\partial{v}}{\partial{s}} \, \Big|_{s=0} = 0 ,
\\ \nonumber
& 
J \left(  \frac{\partial^3{v}}{\partial{t}^2\partial{s}} ( 1 + ( \frac{\partial{v}}{\partial{s}})^2) +  \frac{\partial{v}}{\partial{s}}  (\frac{\partial^2{v}}{\partial{s}\partial{t}})^2  \right)
+
\\ \nonumber
& 
\left( 
\frac{\partial^2{v}}{\partial{s}^2}+  \frac{\partial^2{v}}{\partial{s}^2} ( \frac{\partial{v}}{\partial{s}})^2
+ E_r \, \prescript{RL}{0}{\mathcal D}_{t}^{\alpha} \,  \Big[\frac{\partial^2{v}}{\partial{s}^2}(1 + \frac{1}{2}  (\frac{\partial{v}}{\partial{s}})^2) \Big]
+ \frac{1}{2}  E_r \,  (\frac{\partial{v}}{\partial{s}})^2 \, \prescript{RL}{0}{\mathcal D}_{t}^{\alpha} \,  \frac{\partial^2{v}}{\partial{s}^2}
\right) \,\, \Bigg|_{s=1} = 0 ,
\\ \nonumber
& 
M ( \frac{\partial^2{v}}{\partial{t}^2}\!+\! \ddot{V_b})
\!-\! \frac{\partial{v}}{\partial{s}}\left(
\frac{\partial^2{v}}{\partial{s}^2}\!+\! \frac{\partial^2{v}}{\partial{s}^2} (\frac{\partial{v}}{\partial{s}})^2
\!+\! E_r \, \prescript{RL}{0}{\mathcal D}_{t}^{\alpha} \, \Big[\frac{\partial^2{v}}{\partial{s}^2}(1 \!+\! \frac{1}{2} (\frac{\partial{v}}{\partial{s}})^2)\Big]
\!+\! \frac{1}{2}  E_r \, (\frac{\partial{v}}{\partial{s}})^2 \, \prescript{RL}{0}{\mathcal D}_{t}^{\alpha} \, \frac{\partial^2{v}}{\partial{s}^2}
\right)
\\ \nonumber
& 
+ \left( 
\frac{\partial{v}}{\partial{s}} \, (\frac{\partial^2{v}}{\partial{s}^2})^2  + E_r \, \frac{\partial{v}}{\partial{s}} \, \frac{\partial^2{v}}{\partial{s}^2}  \, \prescript{RL}{0}{\mathcal D}_{t}^{\alpha} \, \frac{\partial^2{v}}{\partial{s}^2} 
\right) \,\, \Bigg|_{s=1}  = 0.
\end{align}
% 

%and we define the test space as 

%$$\mathfrak {B}^{\tau,\nu} \, (\Omega):=\, \prescript{r}{0}{H}^{\tau}\,\Big(I;L^2\,(\Lambda_d)\Big) \cap L^2(I ; \raisebox{\depth}{\scalebox{1}{$\chi $}}_d) $$

%%$$\|v\|_{\mathcal{B}^{\tau,\nu}(\Omega)}=\|v\|^2_{\prescript{r}{0}{H}^{\tau}(I; \, L^2(\Lambda_d))}+\|v\|^2_{L^{2}(I; \, \raisebox{\depth}{\scalebox{0.6}{$\chi $}}_d)}$$
%
%%%%%%%%%%%%%%%%%%%%%%%%%%%
\subsection{Weak Formulation}
\label{Sec: Weak Formulation}
%%%%%%%%%%%%%%%%%%%%%%%%%%%
%
The common practice in analysis of numerical methods for PDEs are mostly concerned with linear equations. The analysis for linear PDEs are well-developed and well-defined, however the nonlinear PDEs still lack such analysis. The linear theories are usually applicable to nonlinear problems if the solution is sufficiently smooth \cite{tadmor2012review}. We do not intend to investigate/develop analysis for our proposed nonlinear model. Instead, by assuming smooth solution, we employ the developed linear theories in our analysis. Let $v: \mathbb{R}^{1+1}\rightarrow\mathbb{R}$ for $\alpha \in (0,1)$ and $\Omega = [0,T]\times [0,L]$. Here, we construct the solution space, $\mathcal{B}^{\alpha} \, (\Omega)$, endowed with proper norms\cite{samiee2019unified2}, in which the corresponding weak form of \eqref{Eq: eqn of motion - dimless - body} can be formulated. If we recall the equation \eqref{Eq: eqn of motion - dimless - body} as E, then:
\begin{align}
\mathcal{B}^{\alpha} \, (\Omega):=\Big\{v\in\prescript{l}{0}{H}^{\alpha}(\Omega) \, \Big|  \int_{\Omega} E \,d\Omega< \infty \Big\}
\end{align}
where
$$\prescript{l}{0}{H}^{\alpha}(\Omega)=\prescript{l}{0}{H}^{\alpha}\,\Big(I;L^2\,(\Omega)\Big) \cap L^2(I ;\prescript{}{0}{H}^{2}(\Omega)) $$
and
$$\prescript{}{0}{H}^{2}(\Omega)=\Big\{v\in{H}^{2}(\Omega) \, \Big| \, \,v \, \Big|_{s=0} = \frac{\partial{v}}{\partial{s}} \, \Big|_{s=0} = 0  \Big\}$$

We obtain the weak form of the problem by multiplying the strong form \eqref{Eq: eqn of motion - dimless - body} with proper test functions {$\tilde{v}(s) \in \mathcal {B}^{\alpha} (\Omega)$ and integrating over the dimensionless spatial computational domain $\Omega_s = [0,L]$.} The test function satisfies the geometric boundary conditions, i.e. $\tilde{v}(0) = \frac{\partial{\tilde{v}}}{\partial{s}}(0) = 0 $. Therefore, by changing the order of integral and temporal derivatives, and through integration by parts, the weak form of problem can be written as
\begin{align}
\label{Eq: weak form - 1}
&\! \int_{0}^{1}\! \frac{\partial^2{v}}{\partial{t}^2}  \tilde{v} ds  
\!+\! \!\int_{0}^{1}\! 
\frac{\partial^2}{\partial{s}^2}\left(  
\frac{\partial^2{v}}{\partial{s}^2} \!+\! \frac{\partial^2{v}}{\partial{s}^2}(\frac{\partial{v}}{\partial{s}})^2 
\!+\! E_r \prescript{RL}{0}{\mathcal D}_{t}^{\alpha}  \Big[\frac{\partial^2{v}}{\partial{s}^2}(1 \!+\! \frac{1}{2} (\frac{\partial{v}}{\partial{s}})^2)\Big] \!+\! \frac{1}{2}  E_r 
(\frac{\partial{v}}{\partial{s}})^2  \prescript{RL}{0}{\mathcal D}_{t}^{\alpha}  \frac{\partial^2{v}}{\partial{s}^2}\right)
\tilde{v}  ds
\\ \nonumber
&
- \int_{0}^{1} 
\frac{\partial}{\partial{s}}\left( 
\frac{\partial{v}}{\partial{s}} \,  (\frac{\partial^2{v}}{\partial{s}^2})^2  
+ E_r
\frac{\partial{v}}{\partial{s}} \, \frac{\partial^2{v}}{\partial{s}^2} \,  \prescript{RL}{0}{\mathcal D}_{t}^{\alpha} \, \frac{\partial^2{v}}{\partial{s}^2}
\right)
\, \tilde{v} \, ds
= - \int_{0}^{1} \ddot{V_b} \, \tilde{v} \, ds ,
\end{align}
where we transfer the spatial derivative load to the test function through integration by parts as
\begin{align}
\label{Eq: weak form - 2}
& \frac{\partial^2}{\partial {t}^2}  \int_{0}^{1} v \tilde{v}  ds  
\!+\! 
\int_{0}^{1} 
\left(  
\frac{\partial^2{v}}{\partial{s}^2}+ \frac{\partial^2{v}}{\partial{s}^2} (\frac{\partial{v}}{\partial{s}})^2 
\!+\! E_r
\prescript{RL}{0}{\mathcal D}_{t}^{\alpha}  \Big[\frac{\partial^2{v}}{\partial{s}^2}(1 \!+\! \frac{1}{2} (\frac{\partial{v}}{\partial{s}})^2)\Big] 
\!+\! \frac{1}{2}  E_r 
(\frac{\partial{v}}{\partial{s}})^2  \prescript{RL}{0}{\mathcal D}_{t}^{\alpha}  \frac{\partial^2{v}}{\partial{s}^2}
\right) \frac{\partial^2{\tilde{v}}}{\partial{s}^2}  ds
\\ \nonumber
&
+
\int_{0}^{1} 
\left( 
\frac{\partial{v}}{\partial{s}} \,  (\frac{\partial^2{v}}{\partial{s}^2})^2  
+ E_r
\frac{\partial{v}}{\partial{s}} \, \frac{\partial^2{v}}{\partial{s}^2}\,  \prescript{RL}{0}{\mathcal D}_{t}^{\alpha} \, \frac{\partial^2{v}}{\partial{s}^2}
\right)
\, \frac{\partial{\tilde{v}}}{\partial{s}} \, ds
+ M (\frac{\partial^2{v}}{\partial{t}^2} + \ddot{V_b}) \,\tilde{v} \, \Big|_{s=1}
\\ \nonumber
& 
+ J \left(  \frac{\partial^3{v}}{\partial{t}^2\partial{s}}( 1 + (\frac{\partial{v}}{\partial{s}})^2) + \frac{\partial{v}}{\partial{s}} (\frac{\partial^2{v}}{\partial{t}\partial{s}})^2  \right) \frac{\partial{\tilde{v}}}{\partial{s}} \, \Big|_{s=1}
= - \ddot{V_b} \, \int_{0}^{1} \tilde{v} \, ds .
\end{align}
By rearranging the terms, we get
\begin{align}
\label{Eq: weak form - 3}
& \frac{\partial^2}{\partial {t}^2} \, \left( \int_{0}^{1} v \, \tilde{v} \, ds +  M \, v \, \tilde{v} \, \Big|_{s=1} + J \, \frac{\partial{v}}{\partial{s}} \, \frac{\partial{\tilde{v}}}{\partial{s}} \, \Big|_{s=1}  \right) 
\\ \nonumber
& 
+
J \left(  \frac{\partial^3{v}}{\partial{t}^2\partial{s}} (\frac{\partial{v}}{\partial{s}})^2 + \frac{\partial{v}}{\partial{s}} \ (\frac{\partial^2{v}}{\partial{t}\partial{s}})^2 \right) \frac{\partial{\tilde{v}}}{\partial{s}} \, \Big|_{s=1}
+ 
\int_{0}^{1} \frac{\partial^2{v}}{\partial{s}^2}\, \frac{\partial^2{\tilde{v}}}{\partial{s}^2} \, ds
+ E_r 
\int_{0}^{1} \prescript{RL}{0}{\mathcal D}_{t}^{\alpha} \, \Big[\frac{\partial^2{v}}{\partial{s}^2}\Big]\,\, \frac{\partial^2{\tilde{v}}}{\partial{s}^2} \, ds
\\ \nonumber
& 
+ 
\int_{0}^{1} \frac{\partial^2{v}}{\partial{s}^2}(\frac{\partial{v}}{\partial{s}})^2 \,\, \frac{\partial^2{\tilde{v}}}{\partial{s}^2} \, ds
+
\int_{0}^{1} \frac{\partial{v}}{\partial{s}} \,  (\frac{\partial^2{v}}{\partial{s}^2})^2  \,\, \frac{\partial{\tilde{v}}}{\partial{s}} \, ds
+ \frac{E_r}{2}
\int_{0}^{1} \prescript{RL}{0}{\mathcal D}_{t}^{\alpha} \, \Big[\frac{\partial^2{v}}{\partial{s}^2}(\frac{\partial{v}}{\partial{s}})^2\Big] \,\, \frac{\partial^2{\tilde{v}}}{\partial{s}^2} \, ds
\\ \nonumber
& 
\!+\! \frac{E_r}{2}  
\int_{0}^{1} (\frac{\partial{v}}{\partial{s}})^2  \prescript{RL}{0}{\mathcal D}_{t}^{\alpha} \Big[\frac{\partial^2{v}}{\partial{s}^2}\Big] \frac{\partial^2{\tilde{v}}}{\partial{s}^2}  ds
\!+\! E_r
\int_{0}^{1} \frac{\partial{v}}{\partial{s}} \frac{\partial^2{v}}{\partial{s}^2} \prescript{RL}{0}{\mathcal D}_{t}^{\alpha}  \Big[\frac{\partial^2{v}}{\partial{s}^2}\Big] \frac{\partial{\tilde{v}}}{\partial{s}}  ds
= - \ddot{V_b} \left( \int_{0}^{1} \tilde{v}  ds \!+\! M  \tilde{v}  \Big|_{s=1} \right).
\end{align}
%
%%%%%%%%%%%%%%%%%%%%%%%%%%%
\subsection{Assumed Mode: A Spectral Galerkin Approximation In Space}
\label{Sec: Spectral Galerkin Method}
%%%%%%%%%%%%%%%%%%%%%%%%%%%
%
We employ the following modal discretization to obtain a reduced-order model of the beam. Therefore,
\begin{align}
\label{Eq: assumed mode}
v(s,t) \simeq v_N(s,t) =  \sum_{n=1}^{N} q_n(t) \, \phi_n(s),
\end{align}
where the spatial functions $\phi_n(s), \,\, n=1,2,\cdots,N$ are assumed \textit{a priori} and the temporal functions $q_n(t), \,\, n=1,2,\cdots,N$ are the unknown modal coordinates. The assumed modes $\phi_n(s)$ in discretization \eqref{Eq: assumed mode} are obtained in Appendix \ref{Sec: App. Eigenvalue Problem of Linear Model} by solving the corresponding eigenvalue problem of linear counterpart of the obtained nonlinear model. Subsequently, we construct the proper finite dimensional spaces of basis/test functions as:
\begin{align}
\label{Eq: Solution/Test Space}
V_N = \tilde{V}_N =  span \, \Big\{ \,\, \phi_n(x) : n = 1,2, \cdots, N \,\, \Big\}.
\end{align}
Since $V_N = \tilde{V}_N \subset V = \tilde{V} $, problem \eqref{Eq: weak form - 3} read as: find $v_N \in V_N$ such that
\begin{align}
\label{Eq: weak form - discrete}
& \frac{\partial^2}{\partial {t}^2}  \left( \int_{0}^{1} v_N  \tilde{v}_N  ds \!+ \! M  v_N  \tilde{v}_N  \Big|_{s=1} \!+\! J  \frac{\partial{{v}_N}}{\partial{s}}  \frac{\partial{\tilde{v}_N}}{\partial{s}} \, \Big|_{s=1}  \right) 
\!+\! 
J \left(  \frac{\partial^3{{v}_N}}{\partial{t}^2\partial{s}} (\frac{\partial{{v}_N}}{\partial{s}})^2 \!+\! \frac{\partial{{v}_N}}{\partial{s}}(\frac{\partial^2{{v}_N}}{\partial{t}\partial{s}})^2  \right) \frac{\partial{\tilde{v}_N}}{\partial{s}}  \Big|_{s=1}
\\ \nonumber
&
+ 
\int_{0}^{1} \frac{\partial^2{{v}_N}}{\partial{s}^2} \, \frac{\partial^2{\tilde{v}_N}}{\partial{s}^2} \, ds
+ E_r 
\int_{0}^{1} \prescript{RL}{0}{\mathcal D}_{t}^{\alpha} \, \Big[\frac{\partial^2{{v}_N}}{\partial{s}^2}\Big] \,\,  \frac{\partial^2{\tilde{v}_N}}{\partial{s}^2} \, ds
+ 
\int_{0}^{1} \frac{\partial^2{{v}_N}}{\partial{s}^2} (\frac{\partial{{v}_N}}{\partial{s}})^2 \,\,  \frac{\partial^2{\tilde{v}_N}}{\partial{s}^2} \, ds
\\ \nonumber
&
+
\int_{0}^{1} \frac{\partial{{v}_N}}{\partial{s}} \,  (\frac{\partial^2{{v}_N}}{\partial{s}})^2  \,\, \frac{\partial{\tilde{v}_N}}{\partial{s}}  \, ds
+ \frac{E_r}{2}
\int_{0}^{1} \prescript{RL}{0}{\mathcal D}_{t}^{\alpha} \, \Big[\frac{\partial^2{{v}_N}}{\partial{s}^2} (\frac{\partial{v_N}}{\partial{s}} )^2\Big] \,\, \frac{\partial^2{\tilde{v}_N}}{\partial{s}^2}\, ds
\\ \nonumber
& 
+ \frac{E_r}{2} 
\int_{0}^{1} (\frac{\partial{{v}_N}}{\partial{s}})^2 \, \prescript{RL}{0}{\mathcal D}_{t}^{\alpha} \, \Big[\frac{\partial^2{{v}_N}}{\partial{s}^2}\Big] \,\, \frac{\partial^2{\tilde{v}_N}}{\partial{s}^2}  \, ds
+ E_r
\int_{0}^{1} \frac{\partial{{v}_N}}{\partial{s}} \, \frac{\partial^2{v_{N}}}{\partial{s}^2} \,  \prescript{RL}{0}{\mathcal D}_{t}^{\alpha} \, \Big[\frac{\partial^2{{v}_N}}{\partial{s}^2}\Big] \,\, \frac{\partial{\tilde{v}_N}}{\partial{s}}  \, ds
\\ \nonumber
& 
= -\ddot{V_b} \left( \int_{0}^{1} \tilde{v}_N \, ds + M \, \tilde{v}_N \, \Big|_{s=1} \right),
\end{align}
for all $ \tilde{v}_N \in \tilde{V}_N$.

%
%%%%%%%%%%%%%%%%%%%%%%%%%%%
\subsection{Single Mode Approximation} \label{subsection}
%%%%%%%%%%%%%%%%%%%%%%%%%%%
%
In general, the modal discretization \eqref{Eq: assumed mode} in \eqref{Eq: weak form - discrete} leads to coupled nonlinear system of fractional ordinary differential equations. We note that while the fractional operators already impose excessive numerical challenges, the nonlinearity will further adds to the complications, leading to failure of existing numerical schemes to solve the coupled system. However, without loss of generality, we can assume that only one mode (primary mode) of motion is involved in the dynamics of system of interest, and thus isolate the rest of modes. In this case, we further reduce the approximation to lower fidelity model by considering single mode approximation.

We assume that the only active mode of vibration is the primary one, which encapsulates most of the fundamental dynamics of our complex system. Therefore, we replace \eqref{Eq: assumed mode} with unimodal discretization $v_N = q(t) \, \phi(s)$, where we let $N=1$ and drop subscript $1$ for simplicity. Upon substituting in \eqref{Eq: weak form - discrete}, we obtain the unimodal governing equation of motion as
\begin{align}
\label{Eq: weak form - discrete 2}
\mathcal{M}  \ddot{q} \!+\! \mathcal{J} (\ddot q  q^2 \!+\! q  {\dot q}^2)
\!+\! \mathcal{K}_l  q \!+\! E_r  \mathcal{C}_l  \prescript{RL}{0}{\mathcal D}_{t}^{\alpha} q  
\!+\! 2 \mathcal{K}_{nl} \, q^3 \!+\! \frac{E_r  \mathcal{C}_{nl}}{2}  \left( \prescript{RL}{0}{\mathcal D}_{t}^{\alpha} q^3 \!+\! 3  q^2  \prescript{RL}{0}{\mathcal D}_{t}^{\alpha} q  \right)
\!=\! -\mathcal{M}_b  \ddot{V_b},
\end{align}
in which
\begin{align}
\label{Eq: coeff unimodal}
& \mathcal{M} = \int_{0}^{1} \phi^2 \, ds +  M \, \phi^2(1) + J \, {\phi^{\prime}}^2(1),  \quad
\mathcal{J} =  J \, {\phi^{\prime}}^4(1),
\\ \nonumber
& \mathcal{K}_l = \mathcal{C}_l = \int_{0}^{1}  {\phi^{\prime\prime}}^2 \,\, ds, \quad
\mathcal{K}_{nl} = \mathcal{C}_{nl} = \int_{0}^{1}  {\phi^{\prime}}^2 \, {\phi^{\prime\prime}}^2 \,\, ds,
\\ \nonumber
& \mathcal{M}_b = \int_{0}^{1} \phi \, ds + M \, \phi(1) .
\end{align}

\begin{rem}
	We note that one can isolate any mode of vibration $\phi_n(s), \,\, n = 1,2, \cdots, N$ (and not necessarily the primary mode) by assuming that $\phi_n(s)$ is the only active one, and thus, end up with similar equation of motion as \eqref{Eq: weak form - discrete 2}, where the coefficients in \eqref{Eq: coeff unimodal} are obtained based on the active mode $\phi_n(s)$. Therefore, we can also make sense of \eqref{Eq: weak form - discrete 2} as a decoupled equation of motion associated with mode $\phi_n(s)$, in which the interaction with other inactive modes is absent.
	
\end{rem}

%\newpage
%
%%%%%%%%%%%%%%%%%%%%%%%%%%%
\section{Linearized Equation: Direct Numerical Time Integration}
\label{Sec: Direct Numerical Time Integration}
\begin{figure}[t]
	\centering
	\begin{subfigure}{0.5\textwidth}
		\centering
		\includegraphics[width=1\linewidth]{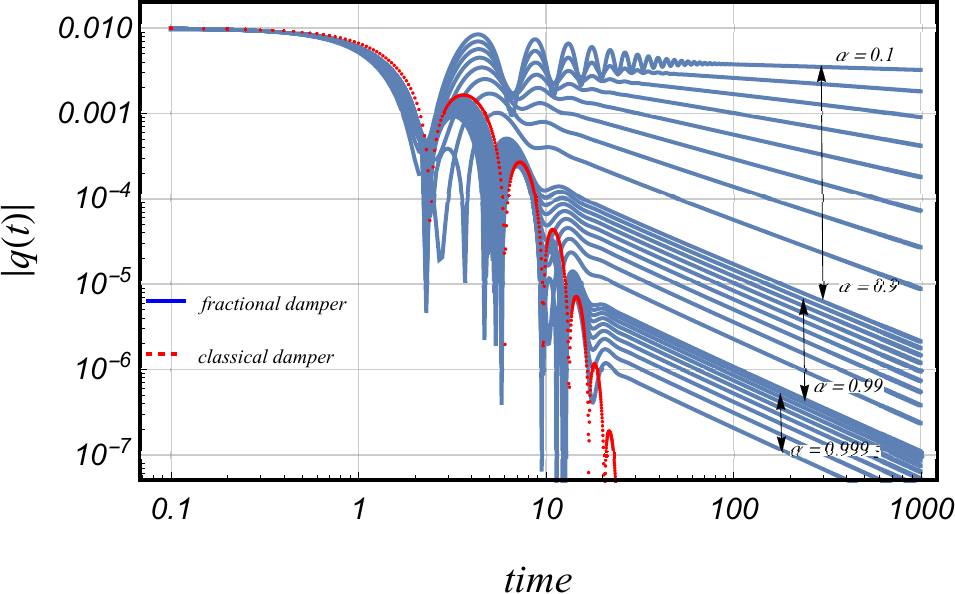}
	\end{subfigure}
	\begin{subfigure}{0.35\textwidth}
		\centering
		\includegraphics[width=1\linewidth]{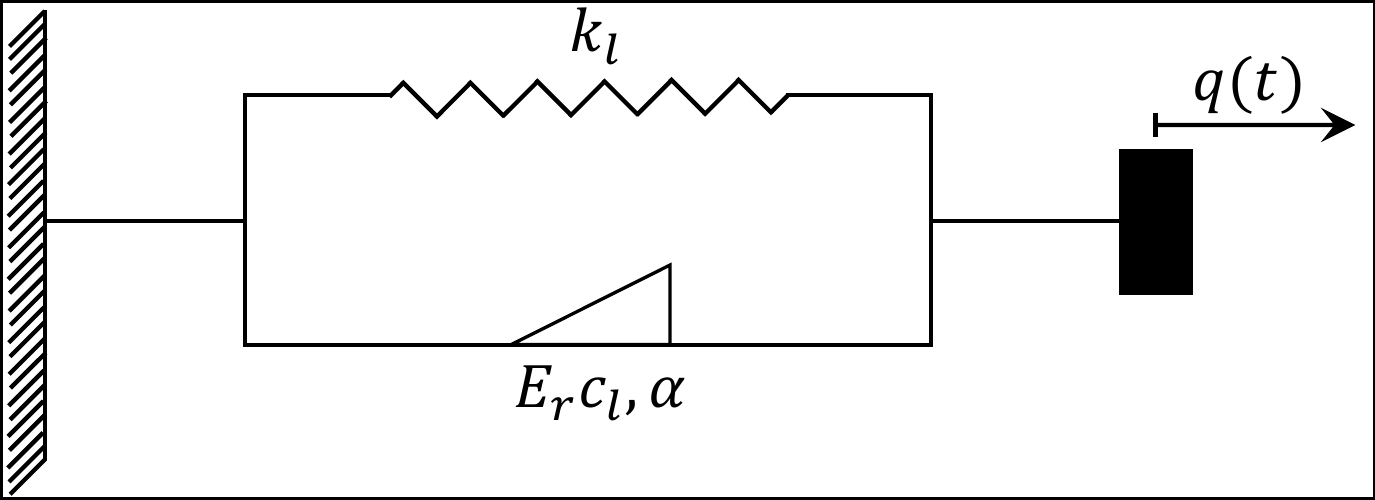}
	\end{subfigure}
	\caption{Power-Law Decay: Time response of linear fractionally damped oscillator using Newmark and $L1$ scheme. The fractional damper has two constants $E_r \, c_l$ and $\alpha$ as the coefficient and derivative order of fractional operator. }
	\label{Fig: fracOscil}
\end{figure}

%%%%%%%%%%%%%%%%%%%%%%%%%%%
%
Since the source of nonlinearity in our problem is coming from geometry, we linearize the equation of single mode approximation, which governs the time evolution of the active mode of vibration in appendix \ref{Sec: App. Linearization}. Therefore, in the absence of base excitation, \eqref{Eq: weak form - discrete 2} takes the following form 
\begin{align}
\label{Eq: weak form linear - discrete 2}
\ddot{q} + E_r \, c_l \, \prescript{RL}{0}{\mathcal D}_{t}^{\alpha} q  + k_l \, q  = 0
\end{align}
in which the coefficients $c_l = \frac{\mathcal{C}_l}{\mathcal{M}}$ and $k_l = \frac{\mathcal{K}_l}{\mathcal{M}}$ are given in \eqref{Eq: coeff unimodal}. The linearized equation \eqref{Eq: weak form linear - discrete 2} can be thought of as a fractionally damped oscillator, shown schematically in Fig. \ref{Fig: fracOscil} (right). This setting describes the vibration of a lumped fractional Kelvin-Voigt model. By letting $E_r = 1$, the dimensionless parameters $c_l = k_l = 1.24$ with a unit mass at the tip, i.e. $ M = 1$. We find the time response of the linearized model \eqref{Eq: weak form linear - discrete 2} using a direct finite difference time integration scheme, which employs $L1$ scheme \cite{lin2007finite,li2012finite} and Newmark method to approximate the fractional derivative and the inertial term, respectively. The Newmark method is of second order accuracy and thus the overall accuracy of the developed scheme is governed by the error level of $L1$ scheme, which is of order $2-\alpha$. 

Figure \ref{Fig: fracOscil} (left) shows the time response of free vibration of 
a fractionally damped oscillator. The absolute value of $q(t)$ versus time for 
different values of $\alpha$ is plotted in Log-Log scale. We observe that in 
the long time, the amplitude of oscillation decays with a power-law, whose rate 
is governed by the order of fractional derivative $\alpha$ and increases by 
increasing $\alpha$ (see blue lines in the figure). By replacing the fractional 
damper with a classical integer-order one, we see that the amplitude decays 
exponentially and not anymore by a power-law (see the dotted red line in the 
figure). These results are in perfect agreement with the power-law and 
exponential relaxation kernel, described in Sec. \ref{Sec: Viscoelasticity}. We 
note that since the fractional element is inherently a viscoelastic element 
that interpolates between the two spring and dash-pot elements (see Sec. 
\ref{Sec: Viscoelasticity} for more discussion and references), it contributes 
both in the stiffness and damping ratio of the system. As $\alpha$ increases, 
the fractional element converges to purely viscous element, and thus the system 
becomes softer (less stiff), resulting in frequency reduction. This frequency 
shift can be seen from the drift of consecutive amplitude peaks to the right as 
$\alpha$ is increased. The fractional linear oscillators are also considered in 
\cite{svenkeson2016spectral} as a case of systems with memory, where their 
interaction with a fluctuating environment causes the time evolution of the 
system to be intermittent. The authors in \cite{svenkeson2016spectral} apply 
the Koopman operator theory to the corresponding integer order system and then 
make a L$\grave{\text{e}}$vy transformation in time to recover long-term memory 
effects; they observe a power-law behavior in the amplitude decay of the 
system's response. Such an anomalous decay rate has also been investigated in 
\cite{shoshani2017anomalous} for an extended theory of decay of classical 
vibrational models brought into nonlinear resonances. The authors report a 
``non-exponential" decay in variables describing the dynamics of the system in 
the presence of dissipation and also a sharp change in the decay rate close to 
resonance. 

\begin{rem}
	The change in fractional derivative order, $\alpha$, is a notion of stiffening/softening of a viscoelastic material modeled via fractional constitutive equations. As shown in Fig. \ref{Fig: fracOscil} (left), the value of $\alpha$ directly affects the decay rate of free vibration. This strong relation can be used to develop a prediction framework, which takes time series of free vibrations as an input, and returns an estimation of the level of material stiffness as a reflection of the health of the system of interest.
\end{rem}

%
%++++++++++++++++++++++++++++++++++++++++++
%++++++++++++++++++++++++++++++++++++++++++
%

%%
%%++++++++++++++++++++++++++++++++++++++++++
%\begin{figure}[t]
%	\centering
%	\includegraphics[width=0.57\linewidth]{figures/LinFO.eps}
%	\caption{Linear fractionally damped oscillator. The fractional damper has two constants $E_r \, c_l$ and $\alpha$ as the coefficient and derivative order of fractional operator.  }
%	\label{Fig: Linear Frac Osc}
%\end{figure}
%%++++++++++++++++++++++++++++++++++++++++++
%%
%%
%%++++++++++++++++++++++++++++++++++++++++++
%\begin{figure}[t]
%	\centering
%	\includegraphics[width=0.57\linewidth]{figures/FracOsci.pdf}
%	\caption{Linear fractionally damped oscillator. The fractional damper has two constants $E_r \, c_l$ and $\alpha$ as the coefficient and derivative order of fractional operator.  }
%	\label{Fig: fracOscil}
%\end{figure}
%%++++++++++++++++++++++++++++++++++++++++++
%%
%

%
%%%%%%%%%%%%%%%%%%%%%%%%%%%
\section{Perturbation Analysis of Nonlinear Equation}
\label{Sec: Perturbatin Analysis}
%%%%%%%%%%%%%%%%%%%%%%%%%%%
%
Nonlinear terms in the fractional differential equation \eqref{Eq: weak form - discrete 2} give rise to expensive time integration schemes.

We use perturbation analysis to investigate the behavior of a nonlinear system, where we reduce the nonlinear fractional differential equation to an algebraic equation to solve for the steady state amplitude and phase of vibration.
%
%%%%%%%%%%%%%%%%%%%%%%%%%%%
\subsection{Method of Multiple Scales}
%%%%%%%%%%%%%%%%%%%%%%%%%%%
%
To investigate the dynamics of the system described by \eqref{Eq: weak form - discrete 2}, we use the method of multiple scales \cite{nayfeh2008nonlinear,rossikhin2010application}. The new independent time scales and the integer-order derivative with respect to them are defined as
\begin{equation}
\label{Eq: time scales}
T_m = \epsilon^{m}\, t, \qquad D_m = \frac{\partial}{\partial T_m} , \qquad m = 0, 1, 2, \cdots.
\end{equation}
It is also convenient to utilize another representation of the fractional derivative as in equation (5.82) in \cite{samko1993fractional}, which according to the Rieman-Liouville fractional derivative, is equivalent to the fractional power of the operator of conventional time-derivative, i.e. $\prescript{RL}{0}{\mathcal D}_{t}^{\alpha} = (\frac{d}{dt})^{\alpha}$. Therefore,  
\begin{align}
\label{Eq: time derivatives}
\frac{d}{dt} &= D_0 + \epsilon D_1 + \cdots ,
\\ \nonumber
\frac{d^2}{dt^2} &= D^2_0 + 2 \epsilon D_0 D_1 + \cdots ,
\\ \nonumber
\prescript{RL}{0}{\mathcal D}_{t}^{\alpha} &= (\frac{d}{dt})^{\alpha} = D^{\alpha}_0 + \epsilon \alpha D^{\alpha-1}_0 D_1 + \cdots ,
\end{align}
The solution $q(t)$ can then be represented in terms of series 
\begin{align}
\label{Eq: MMS expansion}
q(T_0,T_1,\cdots) = q_0(T_0,T_1,\cdots) + \epsilon q_1(T_0,T_1,\cdots) + \epsilon^2 q_2(T_0,T_1,\cdots) + \cdots
\end{align}
We assume that the coefficients in the equation of motion has the following scaling
\begin{align}
\label{Eq: coeff scaling}
\frac{\mathcal{J}}{\mathcal{M}} = \epsilon \, m_{nl} , \quad
\frac{\mathcal{K}_l}{\mathcal{M}} = k_l = \omega^2_0  , \quad
\frac{\mathcal{C}_l}{\mathcal{M}} = \epsilon \, c_l , \quad
\frac{\mathcal{K}_{nl}}{\mathcal{M}} = \epsilon \, k_{nl} , \quad
\frac{\mathcal{C}_{nl}}{\mathcal{M}} = \epsilon \, c_{nl} , 
\end{align}
and the base excitation $- \frac{\mathcal{M}_b}{\mathcal{M}} \, \ddot{V_b} $ is a harmonic function of form $\epsilon \, F \cos(\Omega \, t)$. Thus, \eqref{Eq: weak form - discrete 2} can be expanded as
\begin{align}
\label{Eq: governing eqn unimodal - 4}
&\quad (D^2_0 + 2 \epsilon D_0 D_1 + \cdots)(q_0 + \epsilon q_1 + \cdots) 
\\ \nonumber
&+ \epsilon \, m_{nl} (D^2_0 + 2 \epsilon D_0 D_1 + \cdots)(q_0 + \epsilon q_1 + \cdots)  \,\times\, (q_0 + \epsilon q_1 + \cdots)^2 
\\ \nonumber
&+ \epsilon \, m_{nl} (q_0 + \epsilon q_1 + \cdots)  \,\times\, \left( (D_0 + \epsilon D_1 + \cdots)(q_0 + \epsilon q_1 + \cdots) \right)^2 
\\ \nonumber
&+ \omega^2_0 \, (q_0 + \epsilon q_1 + \cdots) 
\\ \nonumber
&+ \epsilon \, E_r \, c_l \, (D^{\alpha}_0 + \epsilon \alpha D^{\alpha - 1}_0 D_1 + \cdots)(q_0 + \epsilon q_1 + \cdots) 
\\ \nonumber
&+ 2 \epsilon \, k_{nl} \, (q_0 + \epsilon q_1 + \cdots)^3 
\\ \nonumber
&+ \frac{1}{2} \epsilon \, E_r \, c_{nl} \, (D^{\alpha}_0 + \epsilon \alpha D^{\alpha - 1}_0 D_1 + \cdots)(q_0 + \epsilon q_1 + \cdots)^3 
\\ \nonumber
&+ \frac{3}{2} \epsilon \, E_r \, c_{nl} \, (q_0 + \epsilon q_1 + \cdots)^2 \left[ (D^{\alpha}_0 + \epsilon \alpha D^{\alpha - 1}_0 D_1 + \cdots)(q_0 + \epsilon q_1 + \cdots) \right]
\\ \nonumber
& = \epsilon \, F \cos(\Omega \, T_0).
\end{align}
By collecting similar coefficients of zero-th and first orders of $\epsilon$, we obtain the following equations
\begin{alignat}{2}
\mathcal{O}(\epsilon^0) &: \quad D_0^2 q_{0} + \omega_0^2 q_{0} &\,\,=\,\,& 0 ,                   \label{Eq: order zero} \\
\mathcal{O}(\epsilon^{1}) &: \quad D_0^2 q_{1} + \omega_0^2 q_{1} &\,\,=\,\,& -2 D_0 D_1 q_{0} - m_{nl} \, \left( q_0^2 D_0^2 q_0 + q_0 (D_0 q_0)^2 \right) \nonumber \\
&&&- E_r \, c_l \, D^{\alpha}_0 q_{0} - 2 \, k_{nl} \, q^3_{0} - \frac{1}{2} \, E_r \, c_{nl} \, D^{\alpha}_0 q^3_{0} \nonumber \\
&&&- \frac{3}{2} \, E_r \, c_{nl} \, q^2_0 D^{\alpha}_0 q_{0} + F \cos(\Omega \, T_0).         \label{Eq: first order}
\end{alignat}
The solution to \eqref{Eq: order zero} is of the form
\begin{align}
\label{Eq: order zero solution}
q_0(T_0,T_1) = A(T_1) \, e^{i \, \omega_0 \, T_0} + c.c 
%- \frac{f}{2i} \left( e^{i \, \Omega \, T_0} + c.c \right),  \quad f = \frac{F}{\omega^2_0 - 1}
%
\end{align}
where ``c.c" denotes the complex conjugate. By substituting \eqref{Eq: order zero solution} into the right-hand-side of \eqref{Eq: first order}, we observe that different resonance cases are possible. In each case, we obtain the corresponding solvability conditions by removing the secular terms, i.e. the terms that grow in time unbounded. Then, we write $A$ in the polar form $A = \frac{1}{2} a \, e^{i \, \varphi}$, where the real valued functions $a$ and $\varphi$ are the amplitude and phase lag of time response, respectively. Thus, the solution $q(t)$ becomes
\begin{align}
\label{Eq: order zero solution - 2}
q(t) = a(\epsilon \, t) \, \text{cos}(\omega_0 \, t + \varphi(\epsilon \, t)) + \mathcal{O}(\epsilon),
\end{align}
where the governing equations of $a$ and $\varphi$ are obtained by separating the real and imaginary parts.

%
%%%%%%%%%%%%%%%%%%%%%%%%%%%
\subsubsection{Case 1: No Lumped Mass At The Tip}
\label{Sec: case 1}
\begin{figure}[h]
	\centering
	\begin{subfigure}{0.45\textwidth}
		\centering
		\includegraphics[width=1\linewidth]{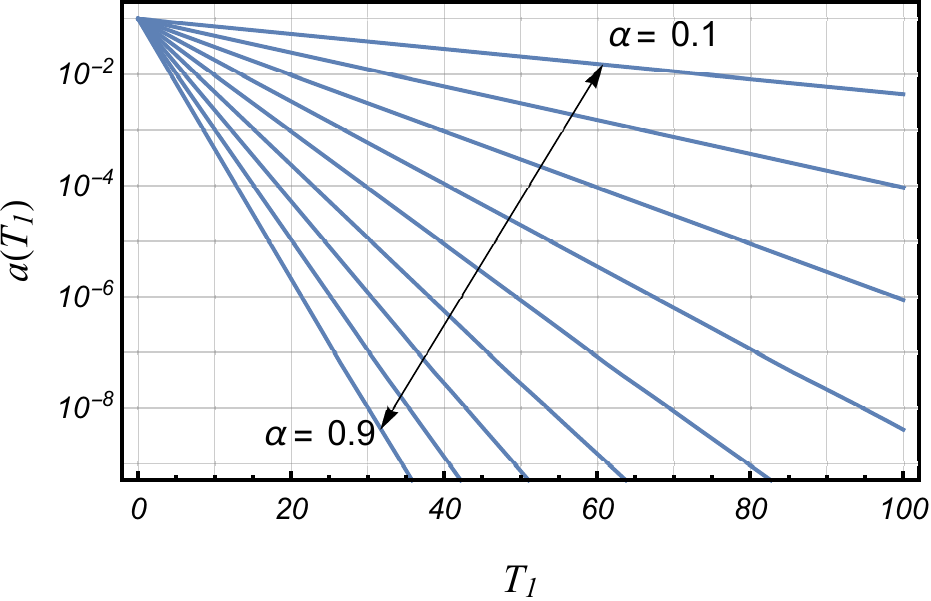}
	\end{subfigure}
	\begin{subfigure}{0.45\textwidth}
		\centering
		\includegraphics[width=1\linewidth]{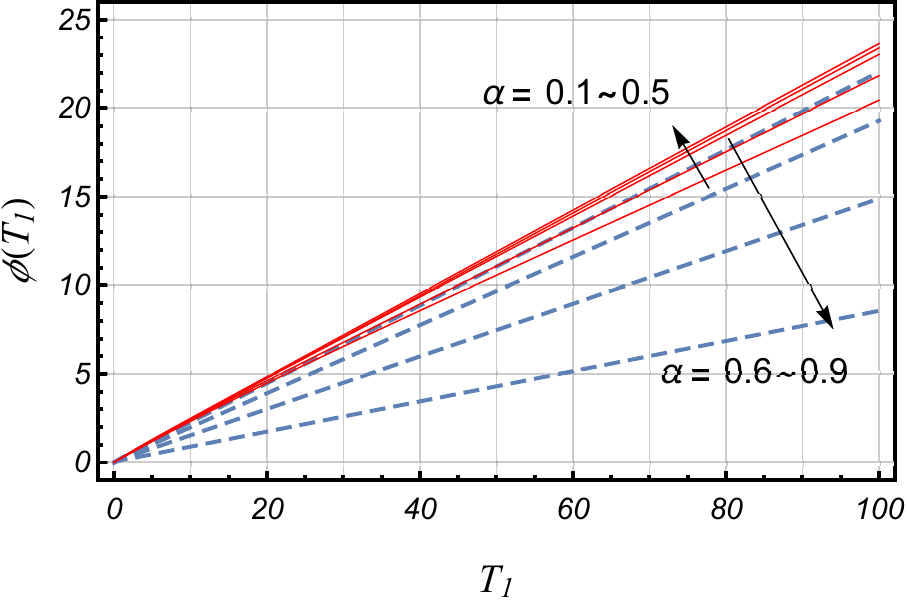}
	\end{subfigure}
	%		\begin{subfigure}{0.32\textwidth}
	%		\centering
	%		\includegraphics[width=1\linewidth]{figures/perturbation/FreeVibphase2.eps}
	%		\end{subfigure}
	\caption{Free vibration of the viscoelastic cantilever beam with no lumped mass at the tip. The rate of decay of amplitude strongly depends on the fractional derivative order $\alpha$ and the coefficient $E_r$. The left figure (log-linear scale) shows the rapid increase in amplitude decaying as $\alpha$ is increased and $E_r = 0.1$. The right figure (linear scale) shows the phase lag $\varphi(\epsilon t)$, where its increase rate decreases as $\alpha$ is increased.}
	\label{Fig: perturbation free vib}
\end{figure}

\begin{figure}[h]
	\centering
	\includegraphics[width=0.5\linewidth]{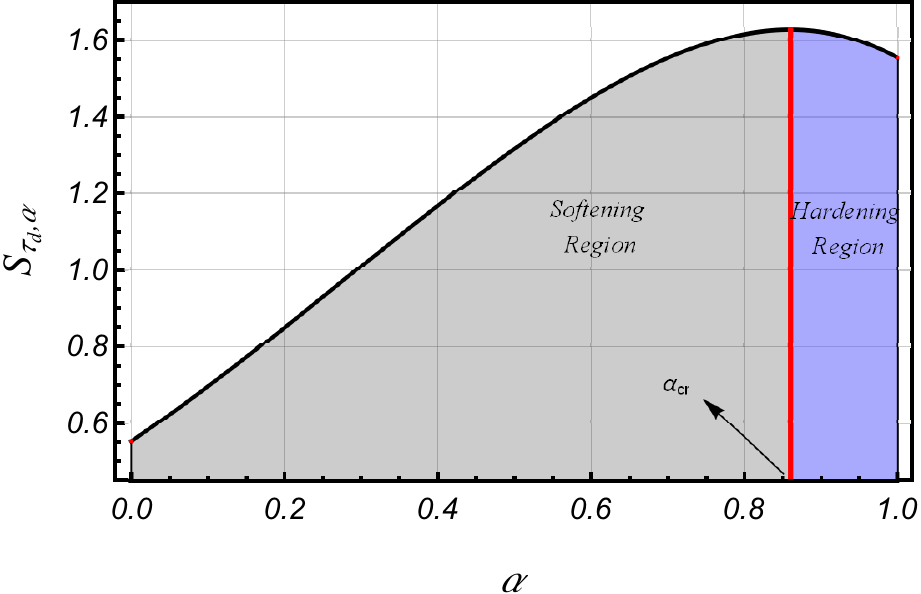}
	\caption{Free vibration of the viscoelastic cantilever beam with no lumped mass at the tip. This graph shows sensitivity of the decay rate $\tau_d$ with respect to change of $\alpha$. Increasing $\alpha$ when $\alpha < \alpha_{cr}$ leads to higher dissipation and decay rate. The reverse effect is observed when $\alpha > \alpha_{cr}$. By softening and hardening we reflect to the regions where increasing $\alpha$ (introducing extra viscosity) leads to higher and lower decay rate, respectively.   }
	\label{Fig: perturbation free vib Sen}
\end{figure}

%%%%%%%%%%%%%%%%%%%%%%%%%%%
%
In this case, $M = J = 0$, and thus, given the functions $\varphi_1(x)$ in Appendix \ref{Sec: App. Eigenvalue Problem of Linear Model}, the coefficients are computed as $\mathcal{M} = 1$, $\mathcal{K}_l = \mathcal{C}_l =12.3624$, $\mathcal{M}_b = 0.782992$, and $\mathcal{K}_{nl} = \mathcal{C}_{nl} = 20.2203$. We consider the following cases:

%========================================
\vspace{0.2 in}
\noindent$\bullet$ Free Vibration, $F = 0$: Super Sensitivity to $\alpha$\\
In this case, the beam is not externally excited and thus, $F = 0$. By removing the secular terms that are the coefficients of $e^{i \, \omega_0 \, T_0}$ in the solvability condition, we find the governing equations of solution amplitude and phase as 
%(the superscript prime denotes the derivative with respect to time scale $T_1$):
%
\begin{align}
\label{Eq: free amp}
&\frac{d a}{d T_1}= 
- E_r \, \omega_0^{\alpha-1} \sin (\alpha \frac{\pi }{2}) 
\left( \frac{1}{2} \, c_l \, a +  \frac{3}{8} \, c_{nl} \, a^3 \right),
\\ \label{Eq: free phase}
& 
\frac{d \varphi}{d T_1} =
\frac{1}{2} c_l \, E_r \, \omega _0^{\alpha-1} \, \cos \left(\frac{\pi  \alpha }{2}\right)  
+ \frac{3}{4} c_{nl} \, E_r \, \omega _0^{\alpha-1} \, \cos \left(\frac{\pi  \alpha }{2}\right) \, a^2 
+\frac{3}{4} \, \omega _0^{-1} \, k_{nl} \, a^2 .
\end{align}
We can see from the first equation \eqref{Eq: free amp} that the amplitude of free vibration decays out, where the decay rate $\tau_d = c_l \, E_r \, \omega_0^{\alpha-1} \sin (\alpha \frac{\pi }{2})$ directly depends on values of the fractional derivative $\alpha$ and the coefficients $E_r$ (see Fig. \ref{Fig: perturbation free vib}). 
%
%++++++++++++++++++++++++++++++++++++++++++
%++++++++++++++++++++++++++++++++++++++++++
%
We introduce the sensitivity index $S_{\tau_d, \alpha}$ as the partial derivative of decay rate with respect to $\alpha$, i.e. 
\begin{align}
\label{Eq: Free Vib Amp Sen}
S_{\tau_d, \alpha} = \frac{d \tau_d}{d \alpha} = \frac{\pi}{2} c_l \, E_r \, \omega_0^{\alpha -1} \, \cos(\alpha \frac{\pi }{2}) + c_l \, E_r \, \omega _0^{\alpha -1} \sin(\alpha \frac{\pi }{2})  \log(\omega _0).
\end{align}
The sensitivity index is computed and plotted in Fig. \ref{Fig: perturbation free vib Sen} for the same set of parameters as in Fig. \ref{Fig: perturbation free vib}. There exists a critical value
\begin{align}
\label{Eq: Free Vib Amp Sen critical}
%
%\frac{\pi}{2} \, \cos(\alpha \frac{\pi }{2}) 
%+ \sin(\alpha \frac{\pi }{2})  \log(\omega _0) = 0
%\quad 
\alpha_{cr}  = 
-\frac{2}{\pi}\tan^{-1} \left(\frac{\pi}{2 \, \log(\omega _0)} \right)  ,
\end{align}
where $(dS_{\tau_d, \alpha} / d\alpha) = 0 $. We observe in Fig.\ref{Fig: 
	perturbation free vib Sen} that by increasing $\alpha$ when $\alpha < 
\alpha_{cr}$, i.e. introducing more viscosity to the system, the dissipation 
rate, and thus decay rate, increases; this can be thought of as a softening 
(stiffness-decreasing) region. Further increasing $\alpha$ when $\alpha > 
\alpha_{cr}$, will reversely results in decrease of decay rate; this can be 
thought of as a hardening (more stiffening) region. We also note that 
$\alpha_{cr}$ solely depends on value of $\omega_0$, given in \eqref{Eq: coeff 
	scaling}, and even though the value of $E_r$ affects decay rate, it does not 
change the value of $\alpha_{cr}$.  
%
%++++++++++++++++++++++++++++++++++++++++++
%++++++++++++++++++++++++++++++++++++++++++
%

Although the observed hardening response after a critical value of $\alpha$ in Fig. \ref{Fig: perturbation free vib Sen} might seem counter-intuitive at first, we remark that here the notions of softening and hardening have a mixed nature regarding energy dissipation and time-scale dependent material stress response, which have anomalous nature for fractional viscoelasticity. Therefore, we demonstrate two numerical tests by purely utilizing the constitutive response of the fractional Kelvin-Voigt model (\ref{Eq: frac KV}) to justify the observed behavior in Fig. \ref{Fig: perturbation free vib Sen} by employing the \textit{tangent loss} and the stress-strain response under monotone loads/relaxation.

\begin{figure}[h]
	\centering
	\begin{subfigure}{0.48\textwidth}
		\centering
		\includegraphics[width=1\linewidth]{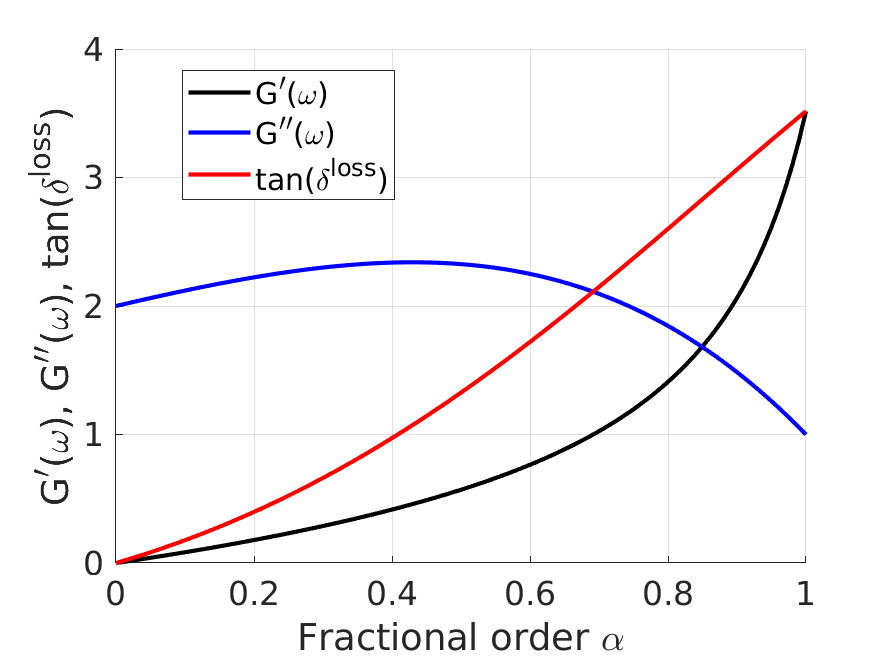}
	\end{subfigure}
	\begin{subfigure}{0.48\textwidth}
		\centering
		\includegraphics[width=1\linewidth]{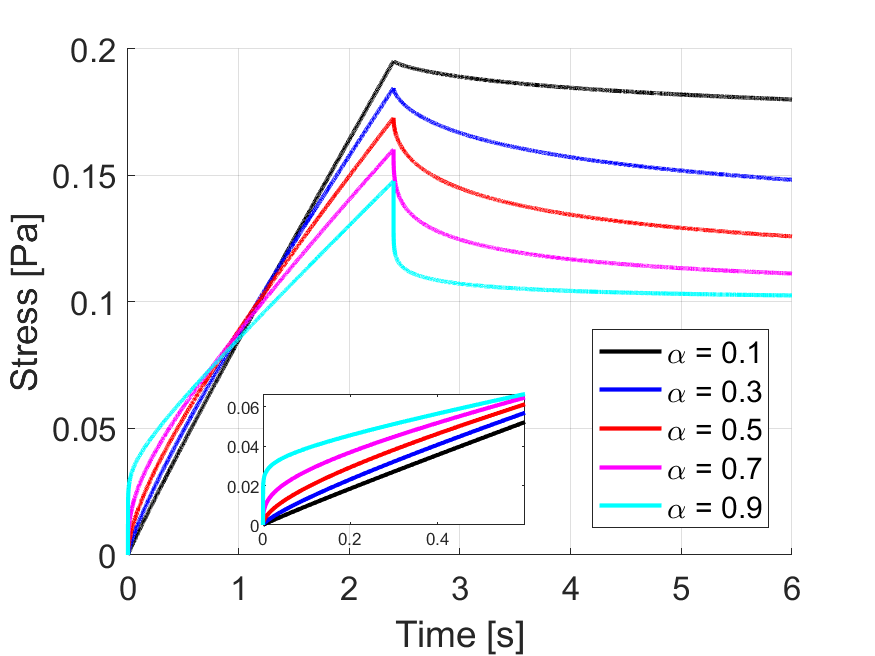}
	\end{subfigure}
	%		\begin{subfigure}{0.32\textwidth}
	%		\centering
	%		\includegraphics[width=1\linewidth]{figures/perturbation/FreeVibphase2.eps}
	%		\end{subfigure}
	\caption{\textit{(Left)} Storage and loss moduli, and tangent loss for the fractional Kelvin-Voigt model at $\omega_0$ with varying fractional-orders and $E_r = 1$. \textit{(Right)} Stress-time response under a monotone load with constant strain rate undergoing ballistic hardening response for short-time and higher $\alpha$, followed by a stress relaxation.}
	\label{Fig:Constitutive_response}
\end{figure}

\noindent\textbf{Dissipation \textit{via} tangent loss}: By taking the Fourier transform of (\ref{Eq: frac KV}), we obtain the so-called complex modulus $G^*$ \cite{mainardi2010fractional}, which is given by:
\begin{equation}
G^*(\omega) = E_\infty + E_\alpha \omega^\alpha \left( \cos\left(\alpha \frac{\pi}{2}\right) + i \sin\left(\alpha \frac{\pi}{2}\right) \right),
\end{equation}
from which the real and imaginary parts yield, respectively, the storage and loss moduli, as follows:
\begin{equation}
G^\prime(\omega) = E_\infty + E_\alpha \omega^\alpha \cos\left(\alpha \frac{\pi}{2}\right), \qquad
G^{\prime\prime}(\omega) = E_\alpha \omega^\alpha \sin\left(\alpha \frac{\pi}{2}\right),
\end{equation}
which represent, respectively, the energy stored and dissipated for each loading cycle. Finally, we define the tangent loss, which represents ratio between dissipated and stored energy cycle, and therefore related to the mechanical damping of the system, as:
\begin{equation}\label{Eq:Tangent_loss}
\tan \delta^{loss} = \frac{G^{\prime\prime}(\omega)}{G^\prime(\omega)} = \frac{E_r \omega^\alpha \sin\left(\alpha\frac{\pi}{2}\right)}{1 + E_r\omega^\alpha\cos\left(\alpha\frac{\pi}{2}\right)}
\end{equation}
We set $\omega = \omega_0$ and $E_r = 1$ and demonstrate the results for (\ref{Eq:Tangent_loss}) with varying fractional orders $\alpha$. We present the obtained results in Fig.\ref{Fig:Constitutive_response} \textit{(left)}, where we observe that increasing fractional orders lead to increased dissipation per loading cycle with the increase of the tangent loss, and the hardening part ($\alpha > \alpha_{c}$) is not associated with higher storage in the material. Instead, the increasing dissipation with $\alpha$ suggests an increasing damping of the mechanical structure.

\noindent\textbf{Stress-time response for monotone loads/relaxation}: In this test, we demonstrate how increasing fractional orders for the fractional model leads to increased hardening for sufficiently high strain rates. Therefore, we directly evaluate (\ref{Eq: frac KV}) with $E_\infty = 1$, $E_\alpha = 1$, and the following strain function: $\varepsilon(t) = (1/24) t$, for $0 \le t < 2.5$ (monotone stress/strain), and $\varepsilon(t) = 1/10$ for $2.5 \le t \le 6$ (stress relaxation). The obtained results are illustrated in Fig.\ref{Fig:Constitutive_response} \textit{(right)}, where we observe that even for relatively low strain rates, there is a ballistic region nearby the initial time where higher fractional orders present higher values of stresses, characterizing a hardening response. However, due to the dissipative nature of fractional rheological elements, after a given time, the initially higher-stress material softens due to its faster relaxation nature.
%

%
%++++++++++++++++++++++++++++++++++++++++++
%++++++++++++++++++++++++++++++++++++++++++

%%========================================
%\vspace{0.2 in}
%\noindent$\bullet$ \textbf{Non-Resonant Case: $\Omega \neq \omega_0$}\\
%In this case, the exciting frequency $\Omega$ is away from the natural frequency $\omega_0$. Thus, the secular terms are the ones containing $e^{i \, \omega_0 \, T_0}$. The governing equations of solution amplitude and phase are:
%%%
%%\begin{align}
%%%
%%\label{Eq: nonresonant amp}
%%a' & = - \frac{c_n \, \omega_0^{\alpha-1}}{2} \,  \text{sin}\left(\frac{\pi \, \alpha}{2}\right) \, \left(48 f^2 + \frac{c_l}{c_n} \right) \, a 
%%- \frac{3 c_n \, \omega_0^{\alpha-1}}{8} \, \text{sin}\left(\frac{\pi \, \alpha}{2}\right) \, \left(2 + 2^{\alpha} \right)  \, a^3
%%%
%%\\ \label{Eq: nonresonant phase}
%%\omega_0 \, \varphi' & = 12 f^2 c_n \left(3 + \frac{k_n}{c_n} + 2 \Omega^2 \, \text{cos}\left(\frac{\pi \, \alpha}{2}\right) \right)
%%+ \frac{c_n \, \omega_0^{\alpha}}{2} \, \text{cos}\left(\frac{\pi \, \alpha}{2}\right) \, \left(48 f^2 + \frac{c_l}{c_n}\right)
%%\\ \nonumber
%%&  \quad + \frac{3 c_n \, \omega_0^{\alpha}}{8} \, \text{cos}\left(\frac{\pi \, \alpha}{2}\right) \, \left(2 + \frac{k_n}{c_n} + \omega_0^2 \left(4 + 2^{\alpha}\right) \, \text{cos}\left(\frac{\pi \, \alpha}{2}\right) \right) \, a^2.
%%%
%%\end{align}
%%%
%
%In this case the amplitude $a$ and thus, the first term in \eqref{Eq: order zero solution} decays rapidly and the steady state amplitude is $f = \frac{F}{\omega^2_0 - 1}$.

\begin{figure}[h]
	\centering
	\begin{subfigure}{0.45\textwidth}
		\centering
		\includegraphics[width=1\linewidth]{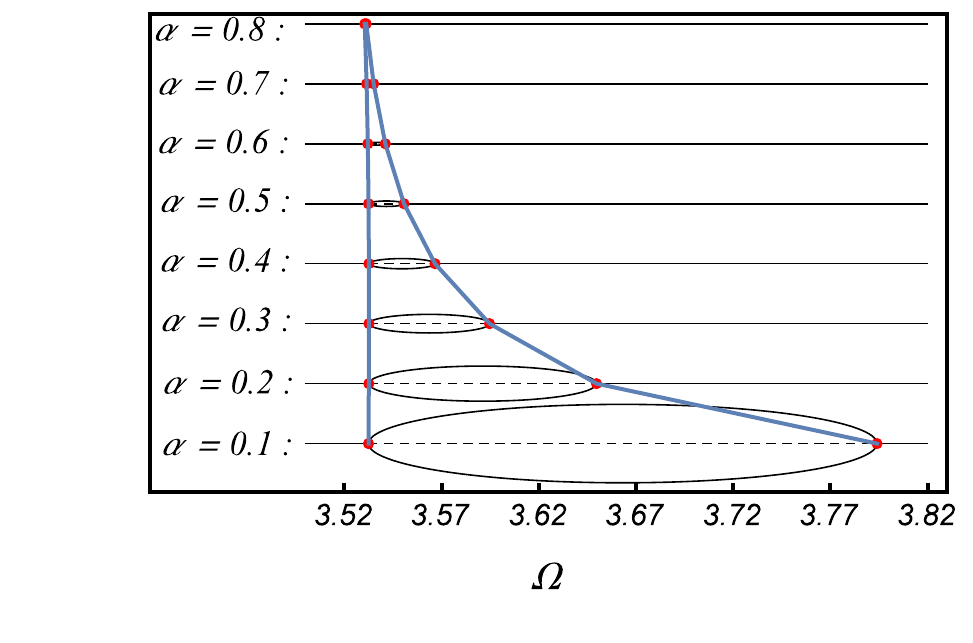}
	\end{subfigure}
	\begin{subfigure}{0.45\textwidth}
		\centering
		\includegraphics[width=1\linewidth]{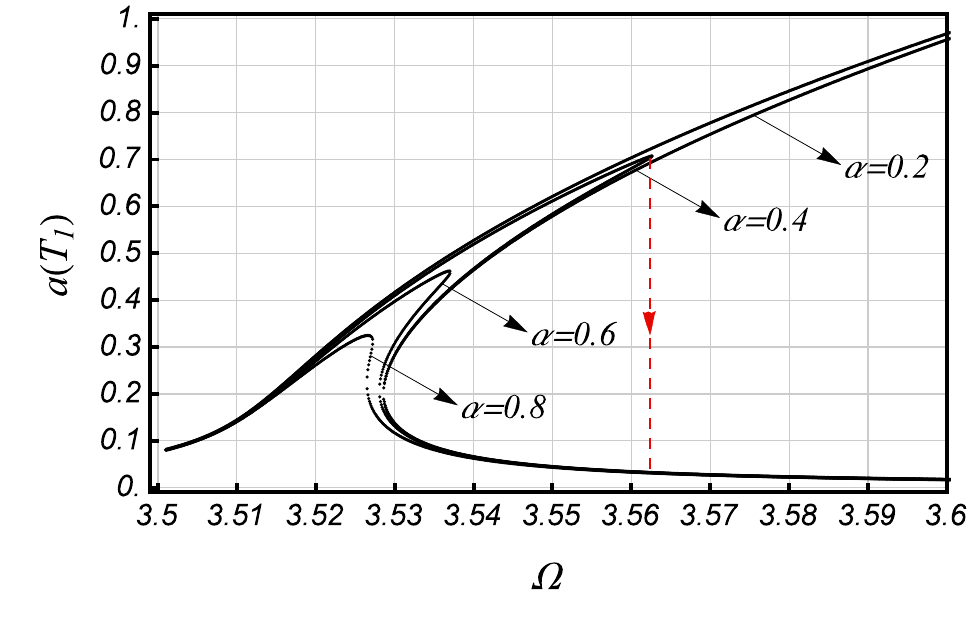}
	\end{subfigure}
	\caption{Primary resonance of the viscoelastic cantilever beam with no lumped mass at the tip. Steady state amplitude (right) and its bifurcation diagram (left) by changing the detuning parameter $\Delta$ for different values of $\alpha$ and $E_r=0.3, f=1$.   }
	\label{Fig: biff diag}
\end{figure}

%========================================
\vspace{0.2 in}
\noindent$\bullet$ Primary Resonance Case, $\Omega \approx \omega_0$\\
In the case of primary resonance, the excitation frequency is close to the natural frequency of the system. We let $\Omega = \omega_0 + \epsilon \, \Delta$, where $\Delta$ is called the detuning parameter and thus, write the force function as $\frac{1}{2} F \, e^{i \, \Delta \, T_1} \, e^{i \, \omega_0 \, T_0} + c.c$ . In this case, the force function also contributes to the secular terms. Therefore, we find the governing equations of solution amplitude and phase as
\begin{align}
\label{Eq: primary amp}
&\frac{d a}{d T_1}
= 
- E_r \, \omega_0^{\alpha-1} \sin (\alpha \frac{\pi }{2}) 
\left( \frac{1}{2} \, c_l  \, a + \frac{3}{8} \, c_{nl}  \, a^3 \right) 
+ \frac{1}{2} f \, \omega^{-1}_0 \, \sin (\Delta  T_1 - \varphi),
\\ \label{Eq: primary phase}
& a \, \frac{d \varphi}{d T_1} =
\frac{1}{2} c_l \, E_r \, \omega _0^{\alpha-1} \, \cos(\frac{\pi  \alpha }{2})  \, a
+ \frac{3}{4} c_{nl} \, E_r \, \omega _0^{\alpha-1} \, \cos(\frac{\pi  \alpha }{2}) \, a^3 
+\frac{3}{4} \, \omega _0^{-1} \, k_{nl} \, a^3 
\\ \nonumber
& \qquad - \frac{1}{2} f \, \omega^{-1}_0 \, \cos (\Delta  T_1 - \varphi) ,
\end{align}
in which the four parameters $\{\alpha,E_r,f,\Delta\}$ mainly change the frequency response of the system. The equations \eqref{Eq: primary amp} and \eqref{Eq: primary phase} can be transformed into an autonomous system, where the $T_1$ does not appear explicitly, by letting $$ \gamma = \Delta \, T_1 - \varphi .$$ The steady state solution occur when $\frac{d a}{d T_1} = \frac{d \varphi}{d T_1} =0$, that gives
\begin{align}
\label{Eq: primary amp ss}
&E_r \, \omega_0^{\alpha-1} \sin (\frac{\pi  \alpha }{2}) 
\left( \frac{c_l}{2} a + \frac{3 c_{nl}}{8}  a^3 \right) 
= \frac{f}{2 \, \omega_0} \sin (\gamma),
\\ \label{Eq: primary phase ss}
& \left(\Delta - \frac{c_l}{2} E_r \, \omega _0^{\alpha-1} \, \cos(\frac{\pi  \alpha }{2}) \right)  a
- \frac{3}{4} \left( c_{nl} \, E_r \, \omega _0^{\alpha-1} \, \cos(\frac{\pi  \alpha }{2}) + \omega _0^{-1} \, k_{nl} \right) a^3 
%\\ \nonumber
%& \qquad 
= \frac{f}{2 \, \omega_0} \cos (\gamma) ,
\end{align}
and thus, by squaring and adding these two equations, we get
\begin{align}
\label{Eq: primary ss}
&\left[
\frac{c_l}{2} E_r \, \omega_0^{\alpha-1} \sin (\frac{\pi  \alpha }{2})  \, a
+\frac{3 c_{nl}}{8}  E_r \, \omega_0^{\alpha-1} \sin (\frac{\pi  \alpha }{2}) \, a^3 
\right]^2
\\ \nonumber
+ & \left[
\left(\Delta - \frac{c_l}{2} E_r \, \omega _0^{\alpha-1} \, \cos(\frac{\pi  \alpha }{2}) \right)  a
- \frac{3}{4} \left( c_{nl} \, E_r \, \omega _0^{\alpha-1} \, \cos(\frac{\pi  \alpha }{2}) + \omega _0^{-1} \, k_{nl} \right) a^3 
\right]^2
%\\ \nonumber
%& \qquad 
= \frac{f^2}{4 \, \omega^2_0}.
\end{align}
This can be written in a simpler way as
\begin{align}
\label{Eq: primary ss - 2}
\left[ A_1 \, a + A_2 \, a^3 \right]^2
+ \left[B_1 \, a + B_2 \, a^3 \right]^2
= C,
\end{align}
where
\begin{alignat*}{4}
&A_1 = \frac{c_l}{2} E_r \, \omega_0^{\alpha-1} \sin (\frac{\pi  \alpha }{2}), 
&& \quad
A_2 = \frac{3 c_{nl}}{8}  E_r \, \omega_0^{\alpha-1} \sin (\frac{\pi  \alpha }{2}), 
\quad
C =\frac{f^2}{4 \, \omega^2_0},
\\
&B_1 = \Delta - \frac{c_l}{2} E_r \, \omega _0^{\alpha-1} \, \cos(\frac{\pi  \alpha }{2}), 
&& \quad
B_2 = - \frac{3}{4} \left( c_{nl} \, E_r \, \omega _0^{\alpha-1} \, \cos(\frac{\pi  \alpha }{2}) + \omega _0^{-1} \, k_{nl} \right).
\end{alignat*}
Hence, the steady state response amplitude is the admissible root of 
\begin{align}
\label{Eq: SS amp}
(A_2^2 + B_2^2) a^6 + (2 A_1 A_2 + 2 B_1 B_2)a^4 + (A_1^2 + B_1^2) a^2 - C =0  ,
\end{align}
which is a cubic equation in $a^2$. The discriminant of a cubic equation of the form $a x^3 + b x^2 + c x + d = 0$ is given as $\vartheta = 18 abcd - 4 b^3 d + b^2 c^2 - 4 a c^3 - 27 a^2 d^2$. The cubic equation \eqref{Eq: SS amp} has one real root when $\vartheta<0$ and three distinct real roots when $\vartheta>0$. The main four parameters $\{\alpha,E_r,f,\Delta\}$ dictate the value of coefficients $\{A_1,A_2,B_1,B_2,C\}$, the value of discriminant $\vartheta$, and thus the number of admissible steady state amplitudes. We see that for fixed values of $\{\alpha,E_r,f\}$, by sweeping the detuning parameter $\Delta$ from lower to higher excitation frequency, the stable steady state amplitude bifurcates into two stable branches and one unstable branch, where they converge back to a stable amplitude by further increasing $\Delta$. Fig. \ref{Fig: biff diag} (left) shows the bifurcation diagram by sweeping the detuning parameter $\Delta$ and for different values of $\alpha$ when $E_r=0.3$ and $f=1$. The solid and dashed black lines are the stable and unstable amplitudes, respectively. The blue lines connect the bifurcation points (red dots) for each value of $\alpha$. We see that the bifurcation points are strongly related to the value of $\alpha$, meaning that by introducing extra viscosity to the system, i.e. increasing the value of $\alpha$, the amplitudes bifurcate and then converge back faster. The right panel of Fig. \ref{Fig: biff diag} shows the frequency response of the system, i.e. the magnitude of steady state amplitudes versus excitation frequency. As the excitation frequency is swept to the right, the steady state amplitude increases, reaches a peak value, and then jumps down (see e.g. red dashed line for $\alpha=0.4$). The peak amplitude and the jump magnitude decreases as $\alpha$ is increased. 
%
%++++++++++++++++++++++++++++++++++++++++++
%++++++++++++++++++++++++++++++++++++++++++
%

\begin{figure}[h]
	\centering
	\begin{subfigure}{0.24\textwidth}
		\centering
		\includegraphics[width=1\linewidth]{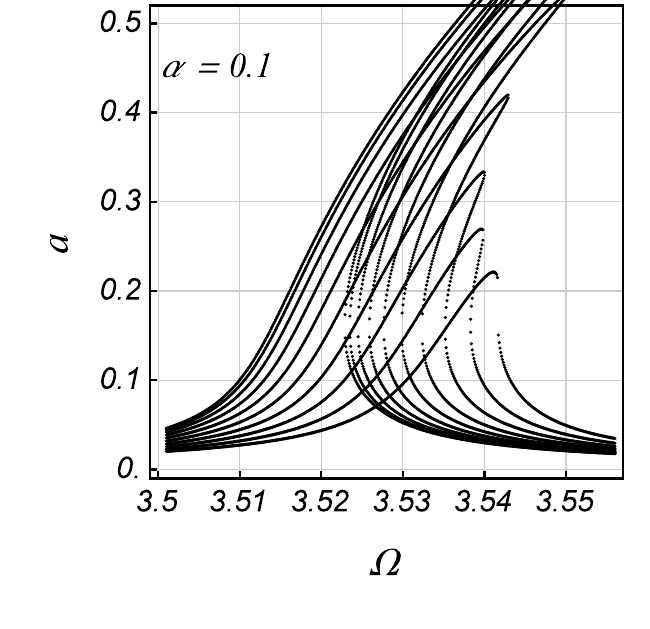}
	\end{subfigure}
	\begin{subfigure}{0.24\textwidth}
		\centering
		\includegraphics[width=1\linewidth]{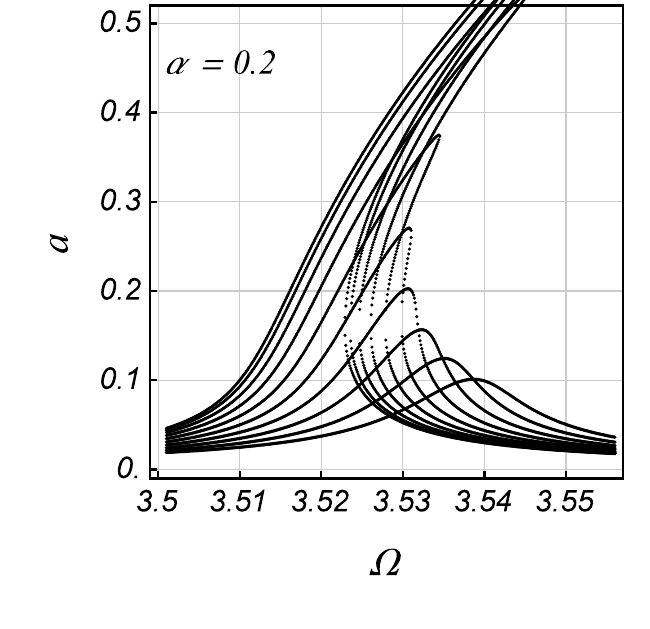}
	\end{subfigure}
	\begin{subfigure}{0.24\textwidth}
		\centering
		\includegraphics[width=1\linewidth]{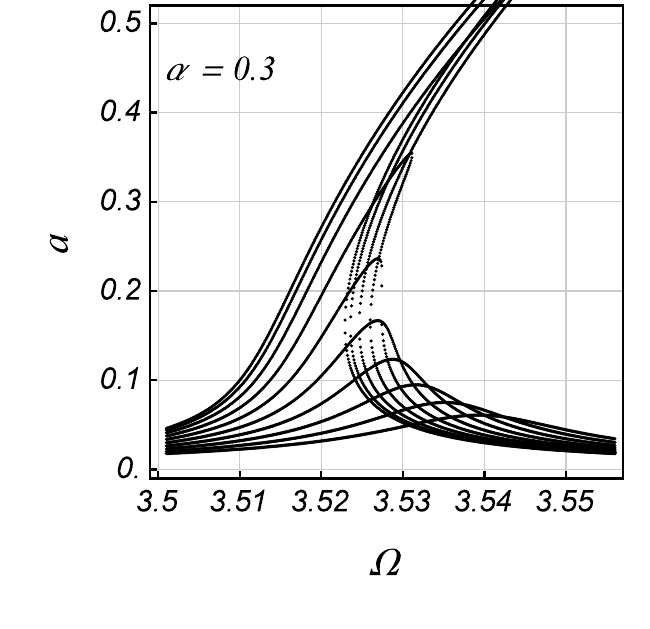}
	\end{subfigure}
	\begin{subfigure}{0.24\textwidth}
		\centering
		\includegraphics[width=1\linewidth]{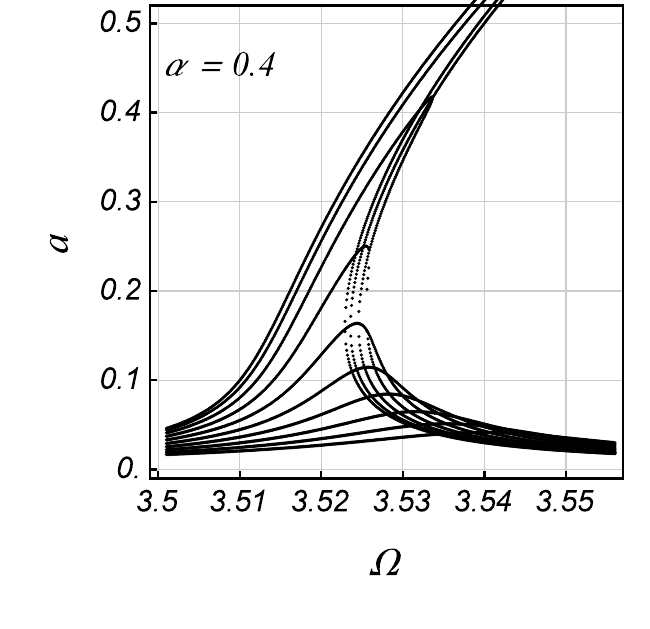}
	\end{subfigure}
	\begin{subfigure}{0.24\textwidth}
		\centering
		\includegraphics[width=1\linewidth]{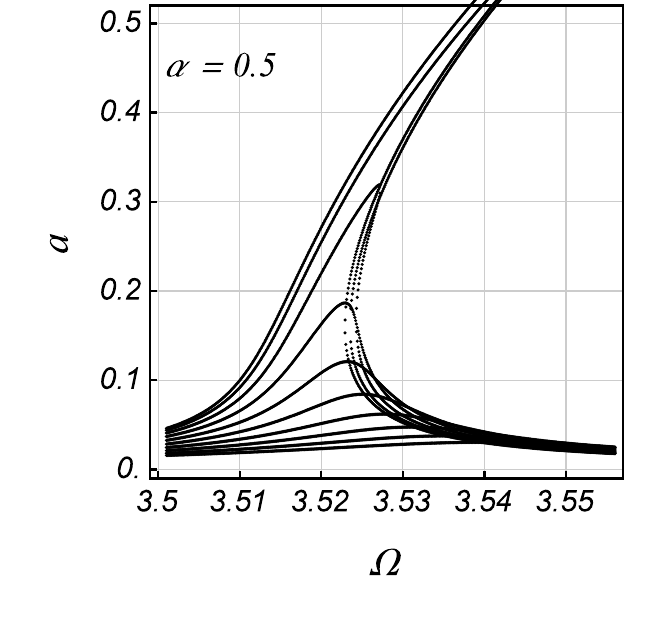}
	\end{subfigure}
	\begin{subfigure}{0.24\textwidth}
		\centering
		\includegraphics[width=1\linewidth]{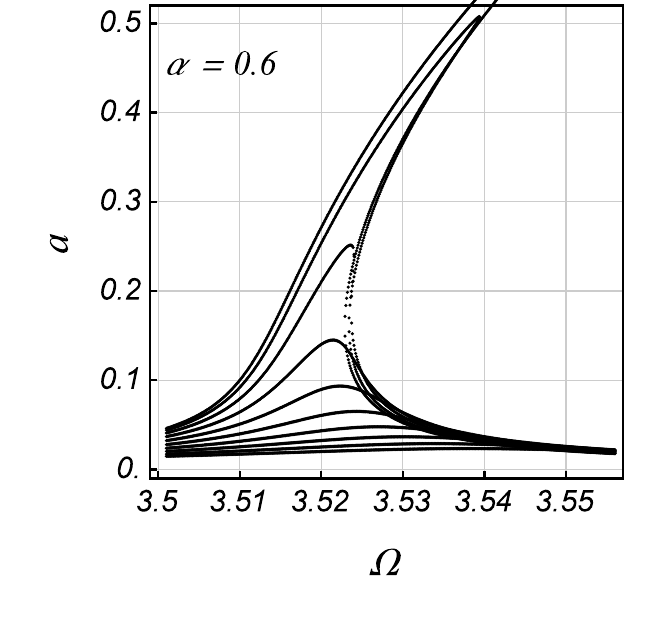}
	\end{subfigure}
	\begin{subfigure}{0.24\textwidth}
		\centering
		\includegraphics[width=1\linewidth]{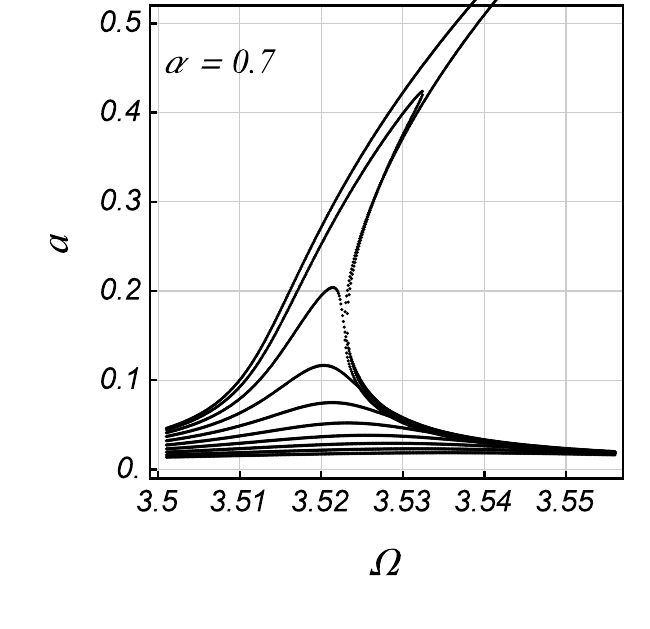}
	\end{subfigure}
	\begin{subfigure}{0.24\textwidth}
		\centering
		\includegraphics[width=1\linewidth]{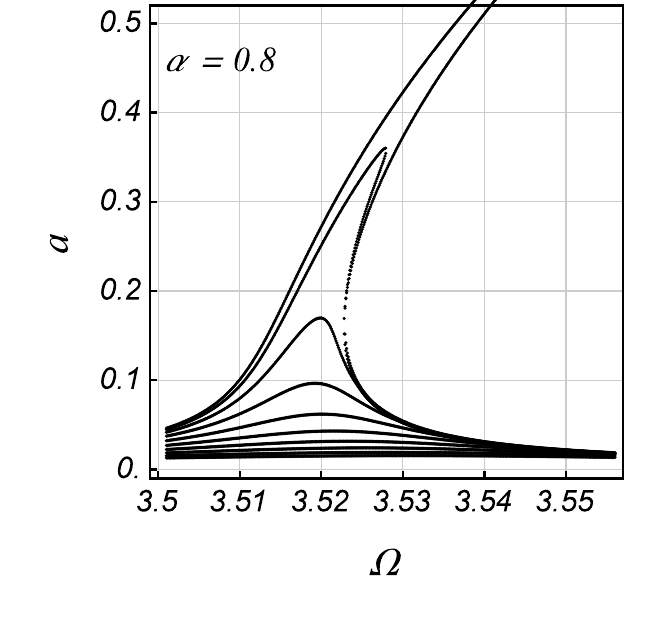}
	\end{subfigure}
	\caption{Frequency-Response curve for the case of primary resonance in the viscoelastic cantilever beam with no lumped mass at the tip. Each sub-figure corresponds to a fixed value of $\alpha$ and $f$ when $E_r = \{0.1, 0.2, \cdots, 1\}$. As effect of fractional element becomes more pronounced, i.e. $\alpha$ and $E_r$ increase, the curve moves down and drift to left. }
	\label{Fig: amp diag}
\end{figure}

The coefficient $E_r = \frac{E_{\infty}}{E_{\alpha}}$ is the proportional contribution of fractional and pure elastic element. At a certain value while increasing this parameter, we see that the bifurcation disappears and the frequency response of system slightly changes. Fig. \ref{Fig: amp diag} shows the frequency response of the system for different values of $\{\alpha,E_r\}$ when $f = 0.5$. In each sub-figure, we let $\alpha$ be fixed and then plot the frequency response for $E_r = \{0.1,0.2,\cdots,1\}$; the amplitude peak moves down as $E_r$ is increased. For higher values of $E_r$, we see that as $\alpha$ is increased, the amplitude peaks drift back to the left, showing a softening behavior in the system response.

%The rows and columns are associated with $E_r = \{0.1,0.2,0.3,0.4,0.5\}$ and $\alpha = \{0.05, 0.45,0.95\}$, respectively. In each figure, the frequency response is plotted for $f= \{0.1,0.2,\cdots,1\}$, while the detuning parameter $\Delta$ is swept from $-1.5$ to $2.5$ with step $0.02$, making the excitation frequency to fall in the range $\Omega \in \left[ 3.501 , 3.541  \right]$. We see a jump in the response amplitude for small values of $\alpha$ and $E_r$. 

%
%++++++++++++++++++++++++++++++++++++++++++
%++++++++++++++++++++++++++++++++++++++++++
%

%\newpage
%
%%%%%%%%%%%%%%%%%%%%%%%%%%%
\subsubsection{Case 2: Lumped Mass At The Tip}
\label{Sec: case 2}
%%%%%%%%%%%%%%%%%%%%%%%%%%%
%
In this case, $M = J = 1$, and thus, given the functions $\phi_1(x)$ in Appendix \ref{Sec: App. Eigenvalue Problem of Linear Model}, the coefficients are computed as $\mathcal{M} = 1 + 70.769 J + 7.2734 M$, $\mathcal{J} = 5008.25$, $\mathcal{K}_l = \mathcal{C}_l = 98.1058$, $\mathcal{M}_b = -0.648623 - 2.69692 M$, and $\mathcal{K}_{nl} = \mathcal{C}_{nl} = 2979.66$. Similar to Case 1, we consider the following cases:

%========================================
\vspace{0.2 in}
\noindent$\bullet$ Free Vibration, $F = 0$\\
Following the same steps as in Case 1, we see that the equation governing amplitude preserve its structure, but the governing equation of phase contains an extra term accommodating the $m_{nl}$.
\begin{align}
\label{Eq: free amp TipMass}
\frac{d a}{d T_1} 
& = 
- E_r \, \omega_0^{\alpha-1} \sin (\alpha \frac{\pi }{2}) 
\left( \frac{1}{2} \, c_l \, a +  \frac{3}{8} \, c_{nl} \, a^3 \right),
\\ \label{Eq: free phase TipMass}
\frac{d \varphi}{d T_1} 
& =
\frac{1}{2} c_l \, E_r \, \omega _0^{\alpha-1} \, \cos \left(\frac{\pi  \alpha }{2}\right)  
+ \frac{3}{4} c_{nl} \, E_r \, \omega _0^{\alpha-1} \, \cos \left(\frac{\pi  \alpha }{2}\right) \, a^2 
\\ \nonumber
& 
+\frac{3}{4} \, \omega _0^{-1} \, k_{nl} \, a^2 
-\frac{1}{4} \, m_{nl} \,  \omega_0 \,  a^2.
\end{align}
This extra term does not significantly alter the behavior of phase and the whole system.

%========================================
\vspace{0.2 in}
\noindent$\bullet$ Primary Resonance Case, $\Omega \approx \omega_0$\\
Similar to the free vibration, we see that the equation governing amplitude preserves its structure while the governing equation of phase contains an extra term accommodating the $m_{nl}$
\begin{align}
\label{Eq: primary amp TipMass}
&\frac{d a}{d T_1} =
- E_r \, \omega_0^{\alpha-1} \sin (\alpha \frac{\pi }{2}) 
\left( \frac{1}{2} \, c_l  \, a + \frac{3}{8} \, c_{nl}  \, a^3 \right) 
+ \frac{1}{2} f \, \omega^{-1}_0 \, \sin (\Delta  T_1 - \varphi),
\\ \label{Eq: primary phase TipMass}
& a \, \frac{d \varphi}{d T_1} =
\frac{1}{2} c_l \, E_r \, \omega _0^{\alpha-1} \, \cos(\frac{\pi  \alpha }{2})  \, a
+ \frac{3}{4} c_{nl} \, E_r \, \omega _0^{\alpha-1} \, \cos(\frac{\pi  \alpha }{2}) \, a^3 
+\frac{3}{4} \, \omega _0^{-1} \, k_{nl} \, a^3 
\\ \nonumber
& \qquad - \frac{1}{2} f \, \omega^{-1}_0 \, \cos (\Delta  T_1 - \varphi) -\frac{1}{4} \, m_{nl} \,  \omega_0 \,  a^3 .
\end{align}
Transforming the equations into an autonomous system by letting $ \gamma = \Delta \, T_1 - \varphi$, we obtain the governing equation of steady state solution as
\begin{align}
\label{Eq: primary ss TipMass}
&\left[
\frac{c_l}{2} E_r \, \omega_0^{\alpha-1} \sin (\frac{\pi  \alpha }{2})  \, a
+\frac{3 c_{nl}}{8}  E_r \, \omega_0^{\alpha-1} \sin (\frac{\pi  \alpha }{2}) \, a^3 
\right]^2
\\ \nonumber
& \left[
\left(\Delta - \frac{c_l}{2} E_r \, \omega _0^{\alpha-1} \, \cos(\frac{\pi  \alpha }{2}) \right)  a
- \frac{3}{4} \left( c_{nl} \, E_r \, \omega _0^{\alpha-1} \, \cos(\frac{\pi  \alpha }{2}) + \omega _0^{-1} \, k_{nl} + \frac{1}{3} \, m_{nl} \,  \omega_0  \right) a^3 
\right]^2
%\\ \nonumber
%& \qquad 
= \frac{f^2}{4 \, \omega^2_0},
\end{align}
which, similar to Case 1, can be written as $$ (A_2^2 + B_2^2) a^6 + (2 A_1 A_2 + 2 B_1 B_2)a^4 + (A_1^2 + B_1^2) a^2 - C =0 ,$$ where all the $A_1$, $A_2$, $B_1$, and $C$ are the same as in Case 1, but $$B_2 = - \frac{3}{4} \left( c_{nl} \, E_r \, \omega _0^{\alpha-1} \, \cos(\frac{\pi  \alpha }{2}) + \omega _0^{-1} \, k_{nl} + \frac{1}{3} \, m_{nl} \,  \omega_0 \right).$$ The corresponding cubic equation can be solved to obtain the bifurcation diagram and also the frequency response of the system. However, in addition to Case 1, we have an extra parameter $m_{nl}$ which affects the response of the system.

\section{Summary and Discussion}
\label{Sec: Summery}
In this work we investigated the nonlinear vibration characteristics of a 
fractional viscoelastic cantilever beam, utilizing a fractional Kelvin-Voigt 
constitutive model. The spectral Galerkin method was employed for spatial 
discretization of the governing equation of motion, reducing it to a set of 
nonlinear fractional ordinary differential equations. The corresponding system 
was linearized and its time-fractional integration was done through direct 
finite-difference scheme, together with a Newmark method. Furthermore, a method 
of multiple scale was employed and the time response of viscoelastic cantilever 
beam subjected to a base excitation was obtained. We performed a set of 
numerical experiments, where we analyzed the response of the system under 
distinct fractional order values, representing different stages of material 
evolution/damage, where we observed:
\begin{itemize}
	\item Super sensitivity of response amplitude to the fractional model parameters at free vibration.
	\item Sensitivity of decay rate with respect to change of $\alpha$, where 
	increasing values of fractional order yielded higher decay rates 
	(softening) before a critical value $\alpha_{cr}$. Beyond such critical 
	values, we observed lower decay rates (stress hardening).
	\item A bifurcation behavior under steady-state amplitude at primary 
	resonance case.
\end{itemize}
The choice of a fractional Kelvin-Voigt model in this work allowed us to 
describe a material in the intersection between anomalous and standard 
constitutive behavior, where the contribution of the Scott-Blair element yields 
the power-law material response, while the Hookean spring reflects the 
instantaneous response of many engineering materials. Furthermore, the changes 
in frequency response due to material evolution provides us a link to the 
development of damage. In terms of modifications of the current model, 
different material distribution functions could be chosen, leading to 
application-based material design for a wide range of 
structural materials and anomalous systems, including microelectromechanical 
systems (MEMS). 
Furthermore, in terms of numerical discretizations, one possibility is to 
utilize more than one active mode of vibration, as well as faster 
time-fractional integration methods, in order to better capture the fundamental 
dynamics of the presented cantilever beam.

\section*{Acknowledgments}
This work was supported by the ARO Young Investigator Program Award (W911NF-19-1-0444), 
and the National Science Foundation Award (DMS-1923201), also partially by MURI/ARO (W911NF-15-1-0562) and the AFOSR Young Investigator Program Award (FA9550-17-1-0150).

\newpage
\appendix

%
%%%%%%%%%%%%%%%%%%%%%%%%%%%%%%%%%%%%%%%%%
\section{Derivation of Governing Equation Using Extended Hamilton's Principle}
\label{Sec: App. extended Hamilton} 
%%%%%%%%%%%%%%%%%%%%%%%%%%%%%%%%%%%%%%%%%
%

%
%%%%%%%%%%%%%%%%%%%%%%%%%%%
\subsection{Equation of Motion}
%%%%%%%%%%%%%%%%%%%%%%%%%%%
%

We recast the integral \eqref{Eq: work var} as $\delta W = \int_{0}^{L} \int_{A} \sigma \, \delta \varepsilon \, dA \,\, ds$ for the considered cantilever beam, in which the variation of strain is $\delta \varepsilon = -\eta \, \delta \frac{\partial{\psi}}{\partial{s}}$, using \eqref{Eq: strain-curvature relation}. 

Therefore, by assuming the constitutive equation \eqref{Eq: frac KV}, the variation of total work is expressed as
\begin{align}
\label{Eq: total work var}
\delta w 
& = 
\int_{0}^{L} \int_{A} \left( -\eta \, E_{\infty} \, \frac{\partial{\psi}}{\partial{s}} -\eta \, E_{\alpha} \, \prescript{RL}{0}{\mathcal D}_{t}^{\alpha} \, \frac{\partial{\psi}}{\partial{s}} \right) \, (-\eta \, \delta \frac{\partial{\psi}}{\partial{s}}) \, dA  \,\, ds
\\ \nonumber
& = 
\int_{0}^{L} \left( E_{\infty} \, \left( \int_{A}  \, \eta^2 dA \right) \, \frac{\partial{\psi}}{\partial{s}} + E_{\alpha} \, \left( \int_{A}  \, \eta^2 dA \right) \, \prescript{RL}{0}{\mathcal D}_{t}^{\alpha} \, \frac{\partial{\psi}}{\partial{s}} \right) \, \delta \frac{\partial{\psi}}{\partial{s}} \,\, ds
\\ \nonumber
& = 
\int_{0}^{L} \left( E_{\infty} \, I \, \frac{\partial{\psi}}{\partial{s}} + E_{\alpha} \, I \, \prescript{RL}{0}{\mathcal D}_{t}^{\alpha} \, \frac{\partial{\psi}}{\partial{s}} \right) \, \delta \frac{\partial{\psi}}{\partial{s}} \,\, ds
\end{align}
where $I = \int_{A} \eta^2 \, dA$. By approximation \eqref{Eq: curvature}, we write the variation of curvature as
\begin{align}
\label{Eq: strain variation}
\delta \frac{\partial{\psi}}{\partial{s}} 
= (1 + \frac{1}{2} (\frac{\partial{v}}{\partial{s}})^2) \, \delta \frac{\partial^2{v}}{\partial{s}^2}+ \frac{\partial^2{v}}{\partial{s}^2} \, \frac{\partial{v}}{\partial{s}} \, \delta \frac{\partial{v}}{\partial{s}}.
\end{align}
Therefore, the variation of total energy becomes
\begin{align}
\label{Eq: total work var - 2}
\delta w 
& \!=\!
\int_{0}^{L}  
\left( E_{\infty}  I  \frac{\partial^2{v}}{\partial{s}^2} (1 + \frac{1}{2} (\frac{\partial{v}}{\partial{s}})^2)  + E_{\alpha}  I  \prescript{RL}{0}{\mathcal D}_{t}^{\alpha}  \Big[\frac{\partial^2{v}}{\partial{s}^2}(1 + \frac{1}{2} (\frac{\partial{v}}{\partial{s}})^2)\Big]\right) 
(1 + \frac{1}{2} (\frac{\partial{v}}{\partial{s}})^2)  \delta \frac{\partial^2{v}}{\partial{s}^2} \,\, ds 
\\ \nonumber
& +
\int_{0}^{L}  
\left( E_{\infty}  I  \frac{\partial^2{v}}{\partial{s}^2} (1 + \frac{1}{2} (\frac{\partial{v}}{\partial{s}})^2)  + E_{\alpha}  I  \prescript{RL}{0}{\mathcal D}_{t}^{\alpha}  \Big[\frac{\partial^2{v}}{\partial{s}^2}(1 + \frac{1}{2} (\frac{\partial{v}}{\partial{s}})^2)\Big]\right) 
\frac{\partial^2{v}}{\partial{s}^2}  \frac{\partial{v}}{\partial{s}}  \delta \frac{\partial{v}}{\partial{s}}  \,\, ds 
\end{align}
By expanding the terms and integrating by parts, we have
\begin{align}
\label{Eq: total work var - 3}
\delta w 
& =
\int_{0}^{L}  
\frac{\partial^2}{\partial{s}^2}\left(
\left( E_{\infty}  I  \frac{\partial^2{v}}{\partial{s}^2} (1 + \frac{1}{2} (\frac{\partial{v}}{\partial{s}})^2)  + E_{\alpha}  I  \prescript{RL}{0}{\mathcal D}_{t}^{\alpha}  \Big[\frac{\partial^2{v}}{\partial{s}^2}(1 + \frac{1}{2} (\frac{\partial{v}}{\partial{s}})^2)\Big]\right) 
(1 + \frac{1}{2} (\frac{\partial{v}}{\partial{s}})^2) 
\right) 
\delta v \,\, ds 
\\ \nonumber
& -
\int_{0}^{L}  
\frac{\partial}{\partial{s}}\left(
\left( E_{\infty}  I  \frac{\partial^2{v}}{\partial{s}^2} (1 + \frac{1}{2} (\frac{\partial{v}}{\partial{s}})^2)  + E_{\alpha} I  \prescript{RL}{0}{\mathcal D}_{t}^{\alpha} \Big[\frac{\partial^2{v}}{\partial{s}^2}(1 + \frac{1}{2} (\frac{\partial{v}}{\partial{s}})^2)\Big]\right) 
\frac{\partial^2{v}}{\partial{s}^2}  \frac{\partial{v}}{\partial{s}} 
\right)
\delta v  \,\, ds 
\\ \nonumber
& +
\left( E_{\infty}  I  \frac{\partial^2{v}}{\partial{s}^2} (1 + \frac{1}{2} (\frac{\partial{v}}{\partial{s}})^2)  + E_{\alpha}  I  \prescript{RL}{0}{\mathcal D}_{t}^{\alpha}  \Big[\frac{\partial^2{v}}{\partial{s}^2}(1 + \frac{1}{2} (\frac{\partial{v}}{\partial{s}})^2)\Big]\right) 
(1 + \frac{1}{2} (\frac{\partial{v}}{\partial{s}})^2)  \delta \frac{\partial{v}}{\partial{s}} \Bigg|_{0}^{L}
\\ \nonumber
& -
\frac{\partial}{\partial{s}}\left(
\left( E_{\infty}  I  \frac{\partial^2{v}}{\partial{s}^2} (1 + \frac{1}{2} (\frac{\partial{v}}{\partial{s}})^2)  + E_{\alpha}  I  \prescript{RL}{0}{\mathcal D}_{t}^{\alpha}  \Big[\frac{\partial^2{v}}{\partial{s}^2}(1 + \frac{1}{2} (\frac{\partial{v}}{\partial{s}})^2)\Big]\right) 
(1 + \frac{1}{2} (\frac{\partial{v}}{\partial{s}})^2) \right) \delta v \Bigg|_{0}^{L}
\\ \nonumber
& +
\left( E_{\infty}  I  \frac{\partial^2{v}}{\partial{s}^2} (1 + \frac{1}{2} (\frac{\partial{v}}{\partial{s}})^2)  + E_{\alpha}  I  \prescript{RL}{0}{\mathcal D}_{t}^{\alpha} \Big[\frac{\partial^2{v}}{\partial{s}^2}(1 + \frac{1}{2} (\frac{\partial{v}}{\partial{s}})^2)\Big]\right) 
\frac{\partial^2{v}}{\partial{s}^2}  \frac{\partial{v}}{\partial{s}}  \delta v \Bigg|_{0}^{L}
\end{align}
The prescribed geometry boundary conditions at the base of the beam, $s=0$, allow the variation of deflection and its first derivative to be zero at $s=0$, i.e. $\delta v(0,t) = \delta \frac{\partial{v}}{\partial{s}}(0,t) = 0$. Therefore,
\begin{align}
\label{Eq: total work var - 4}
\delta w 
& \!=\!
\!\int_{0}^{L}\!  
\frac{\partial^2}{\partial{s}^2}\left(
\left( E_{\infty}  I  \frac{\partial^2{v}}{\partial{s}^2} (1 + \frac{1}{2} (\frac{\partial{v}}{\partial{s}})^2)  + E_{\alpha}  I  \prescript{RL}{0}{\mathcal D}_{t}^{\alpha}  \Big[\frac{\partial^2{v}}{\partial{s}^2}(1 + \frac{1}{2} (\frac{\partial{v}}{\partial{s}})^2)\Big]\right) 
(1 + \frac{1}{2} (\frac{\partial{v}}{\partial{s}})^2) 
\right) 
\delta v \,\, ds 
\\ \nonumber
& -
\int_{0}^{L}  
\frac{\partial}{\partial{s}}\left(
\left( E_{\infty}  I  \frac{\partial^2{v}}{\partial{s}^2} (1 + \frac{1}{2} (\frac{\partial{v}}{\partial{s}})^2)  + E_{\alpha}  I  \prescript{RL}{0}{\mathcal D}_{t}^{\alpha}  \Big[\frac{\partial^2{v}}{\partial{s}^2}(1 + \frac{1}{2} (\frac{\partial{v}}{\partial{s}})^2)\Big]\right) 
\frac{\partial^2{v}}{\partial{s}^2}  \frac{\partial{v}}{\partial{s}}  
\right)
\delta v  \,\, ds 
\\ \nonumber
& +
\left( E_{\infty}  I  \frac{\partial^2{v}}{\partial{s}^2} (1 + \frac{1}{2} (\frac{\partial{v}}{\partial{s}})^2)  + E_{\alpha}  I  \prescript{RL}{0}{\mathcal D}_{t}^{\alpha}  \Big[\frac{\partial^2{v}}{\partial{s}^2}(1 + \frac{1}{2} (\frac{\partial{v}}{\partial{s}})^2)\Big]\right) 
(1 + \frac{1}{2} (\frac{\partial{v}}{\partial{s}})^2) \Bigg|_{s=L}  \delta \frac{\partial{v}}{\partial{s}}(L,t) 
\\ \nonumber
& \!-\!
\frac{\partial}{\partial{s}}\left(
\left( E_{\infty}  I  \frac{\partial^2{v}}{\partial{s}^2} (1 \!+\! \frac{1}{2} (\frac{\partial{v}}{\partial{s}})^2)  \!+\! E_{\alpha}  I  \prescript{RL}{0}{\mathcal D}_{t}^{\alpha}  \Big[\frac{\partial^2{v}}{\partial{s}^2}(1 \!+\! \frac{1}{2} (\frac{\partial{v}}{\partial{s}})^2)\Big]\right) 
(1 \!+\! \frac{1}{2} (\frac{\partial{v}}{\partial{s}})^2) \right) \Bigg|_{s=L}  \delta v(L,t) 
\\ \nonumber
& +
\left( E_{\infty}  I  \frac{\partial^2{v}}{\partial{s}^2} (1 + \frac{1}{2} (\frac{\partial{v}}{\partial{s}})^2)  + E_{\alpha}  I  \prescript{RL}{0}{\mathcal D}_{t}^{\alpha}  \Big[\frac{\partial^2{v}}{\partial{s}^2}(1 + \frac{1}{2} (\frac{\partial{v}}{\partial{s}})^2\Big]\right) 
\frac{\partial^2{v}}{\partial{s}^2}  \frac{\partial{v}}{\partial{s}} \Bigg|_{s=L}  \delta v(L,t) 
\end{align}
Let $\varrho$ be mass per unit volume of the beam, $M$ and $J$ be the mass and rotatory inertia of the lumped mass at the tip of beam. By considering the displacement and velocity of the beam given in \eqref{Eq: displacement} and \eqref{Eq: velocity}, respectively, the kinetic energy is obtained as
\begin{align}
\label{Eq: kinetic}
T 
& = \frac{1}{2} \, \int_{0}^{L} \int_{A} \varrho  \, (\frac{\partial{\textbf{r}}}{\partial{t}})^2 \, dA \, ds 
+ \frac{1}{2}M \left( (\frac{\partial{u}}{\partial{t}})^2  + (\frac{\partial{v}}{\partial{t}} + \dot V_b)^2 \right) \Big|_{s=L} + \frac{1}{2} J (\frac{\partial{\psi}}{\partial{t}})^2 \Big|_{s=L} , 
\\ \nonumber 
& = \frac{1}{2} \, \int_{0}^{L} \int_{A} \varrho  \, \Big\{ (\frac{\partial{u}}{\partial{t}} - \eta \, \frac{\partial{\psi}}{\partial{t}} \, \cos(\psi) )^2 + (\frac{\partial{v}}{\partial{t}}+ \dot V_b - \eta \, \frac{\partial{\psi}}{\partial{t}} \, \sin(\psi))^2 \Big\}  \,  dA \, ds 
\\ \nonumber
&
+ \frac{1}{2}M \left((\frac{\partial{u}}{\partial{t}})^2 + (\frac{\partial{v}}{\partial{t}} + \dot V_b)^2 \right) \Big|_{s=L} + \frac{1}{2} J (\frac{\partial{\psi}}{\partial{t}})^2 \Big|_{s=L} , 
\\ \nonumber
& = 
\frac{1}{2} \, \int_{0}^{L} \int_{A} \varrho  \,
\Big\{ (\frac{\partial{u}}{\partial{t}} )^2 - 2 \eta \, \frac{\partial{u}}{\partial{t}}  \, \frac{\partial{\psi}}{\partial{t}}  \, \cos(\psi) + \eta^2 \, (\frac{\partial{\psi}}{\partial{t}})^2 \, \cos^2(\psi) 
+ (\frac{\partial{v}}{\partial{t}} )^2 + {\dot V_b}^2 + 2 \frac{\partial{v}}{\partial{t}}  \, \dot V_b
\\ \nonumber
& 
- 2 \eta \, (\frac{\partial{v}}{\partial{t}} + \dot V_b) \frac{\partial{\psi}}{\partial{t}}  \, \sin(\psi) + \eta^2 (\frac{\partial{\psi}}{\partial{t}})^2 \, \sin^2(\psi) 
\Big\}  \,  dA \, ds + \frac{1}{2}M \left((\frac{\partial{u}}{\partial{t}})^2 + (\frac{\partial{v}}{\partial{t}} + \dot V_b)^2 \right) \Big|_{s=L}
\\ \nonumber
&
+ \frac{1}{2} J (\frac{\partial{\psi}}{\partial{t}})^2 \Big|_{s=L} , 
\\ \nonumber
& = 
\frac{1}{2} \, \int_{0}^{L} \int_{A} \varrho  \,
\Big\{ (\frac{\partial{u}}{\partial{t}})^2 + (\frac{\partial{v}}{\partial{t}})^2 + {\dot V_b}^2 + 2 \frac{\partial{v}}{\partial{t}} \, \dot V_b 
- 2 \eta \, \frac{\partial{u}}{\partial{t}}\, \frac{\partial{\psi}}{\partial{t}} \, \cos(\psi) + \eta^2 \, (\frac{\partial{\psi}}{\partial{t}})^2 
\\ \nonumber
&
- 2 \eta \, (\frac{\partial{v}}{\partial{t}} + \dot V_b) \frac{\partial{\psi}}{\partial{t}} \, \sin(\psi)  
\Big\}  \,  dA \, ds + \frac{1}{2}M \left((\frac{\partial{u}}{\partial{t}} )^2 + (\frac{\partial{v}}{\partial{t}}  + \dot V_b)^2 \right) \Big|_{s=L} + \frac{1}{2} J (\frac{\partial{\psi}}{\partial{t}})^2 \Big|_{s=L} .
\end{align}
Let
\begin{align*}
\rho = \int_{A} \varrho dA , \quad
\mathcal{J}_1 = \int_{A} \varrho \, \eta \,  dA , \quad
\mathcal{J}_2 = \int_{A} \varrho \, \eta^2 dA.
\end{align*}
$\rho$ is the mass per unit length of the beam, $\mathcal{J}_1$ is the first moment of inertia and is zero because the reference point of coordinate system attached to the cross section coincides with the mass centroid, and $\mathcal{J}_2$ is the second moment of inertia, which is very small for slender beam and can be ignored \cite{hamdan1997large}. Assuming that the velocity along the length of the beam, $\dot u$, is relatively small compared to the lateral velocity $\dot v + \dot V_b$, the kinetic energy of the beam can be reduced to
\begin{align}
\label{Eq: kinetic - 2}
T = 
\frac{1}{2} \, \rho \int_{0}^{L} ( \frac{\partial{v}}{\partial{t}} + \dot V_b)^2  \,   ds 
+ \frac{1}{2}M (\frac{\partial{v}}{\partial{t}} + \dot V_b)^2  \Big|_{s=L} + \frac{1}{2} J (\frac{\partial{\psi}}{\partial{t}})^2 \Big|_{s=L}, 
\end{align}
where its variation can be taken as
\begin{align}
\label{Eq: kinetic variation}
\delta T 
& = \rho \int_{0}^{L} (\frac{\partial{v}}{\partial{t}} + \dot{V_b}) \, \delta \frac{\partial{v}}{\partial{t}} \, ds + M (\frac{\partial{v}}{\partial{t}} + \dot{V_b}) \, \delta \frac{\partial{v}}{\partial{t}} \Big|_{s=L} +  J \frac{\partial{\psi}}{\partial{t}} \, \delta \frac{\partial{\psi}}{\partial{t}} \Big|_{s=L} ,
\end{align}
in which $\frac{\partial{\psi}}{\partial{t}}$ is given in \eqref{Eq: angular velocity} and $\delta \frac{\partial{\psi}}{\partial{t}}$ can be obtained as $\delta \frac{\partial{\psi}}{\partial{t}} \simeq  ( 1 + \frac{1}{2} (\frac{\partial{v}}{\partial{s}})^2) \delta \frac{\partial^2{v}}{\partial{t}\partial{s}}+ \frac{\partial{v}}{\partial{s}} \frac{\partial^2{v}}{\partial{t}\partial{s}} \delta \frac{\partial{v}}{\partial{s}}$. Therefore, 
\begin{align}
\label{Eq: kinetic variation - 2}
\delta T 
& \!\simeq\! \rho \!\int_{0}^{L}\! (\frac{\partial{v}}{\partial{t}} \!+\! \dot{V_b})  \delta \frac{\partial{v}}{\partial{t}}  ds 
\!+\! M (\frac{\partial{v}}{\partial{t}} \!+\! \dot{V_b})  \delta \frac{\partial{v}}{\partial{t}} \Big|_{s=L} 
\!+\! J \left( \frac{\partial^2{v}}{\partial{t}\partial{s}} ( 1 \!+\! (\frac{\partial{v}}{\partial{s}})^2) \delta \frac{\partial^2{v}}{\partial{t}\partial{s}} \!+\! \frac{\partial{v}}{\partial{s}} (\frac{\partial^2{v}}{\partial{t}\partial{s}})^2 \delta \frac{\partial{v}}{\partial{s}} \right) \Big|_{s=L} .
\end{align}
The time integration of $\delta T$ takes the following form through integration by parts
\begin{align}
\label{Eq: kinetic var - 2}
&\int_{t_1}^{t_2} \delta T \, dt 
\\ \nonumber
= & \int_{t_1}^{t_2} 
\Bigg\{
\rho \int_{0}^{L} (\frac{\partial{v}}{\partial{t}} + \dot{V_b}) \, \delta \frac{\partial{v}}{\partial{t}} \, ds 
+ M (\frac{\partial{v}}{\partial{t}} + \dot{V_b}) \, \delta \frac{\partial{v}}{\partial{t}} \Big|_{s=L} 
\\ \nonumber
& 
+ J \left( \frac{\partial^2{v}}{\partial{t}\partial{s}} ( 1 + (\frac{\partial{v}}{\partial{s}})^2) \delta \frac{\partial^2{v}}{\partial{t}\partial{s}} + \frac{\partial{v}}{\partial{s}} (\frac{\partial^2{v}}{\partial{t}\partial{s}})^2 \delta \frac{\partial{v}}{\partial{s}} \right) \Big|_{s=L}
\Bigg\} \, dt
\\ \nonumber
= &
\int_{t_1}^{t_2} \rho \int_{0}^{L} (\frac{\partial{v}}{\partial{t}} + \dot{V_b}) \, \delta \frac{\partial{v}}{\partial{t}} \, ds \, dt
+ M \int_{t_1}^{t_2} (\frac{\partial{v}}{\partial{t}} + \dot{V_b}) \, \delta \frac{\partial{v}}{\partial{t}} \Big|_{s=L} \, dt
\\ \nonumber
& 
+ J \int_{t_1}^{t_2} \left( \frac{\partial^2{v}}{\partial{t}\partial{s}} ( 1 + (\frac{\partial{v}}{\partial{s}})^2) \delta \frac{\partial^2{v}}{\partial{t}\partial{s}} + \frac{\partial{v}}{\partial{s}} (\frac{\partial^2{v}}{\partial{t}\partial{s}})^2 \delta \frac{\partial{v}}{\partial{s}}  \right) \Big|_{s=L} \, dt
\\ \nonumber
= &
\rho \int_{0}^{L} \int_{t_1}^{t_2}  (\frac{\partial{v}}{\partial{t}} + \dot{V_b}) \, \delta \frac{\partial{v}}{\partial{t}} \, dt \, ds
+ M \int_{t_1}^{t_2} (\frac{\partial{v}}{\partial{t}} + \dot{V_b}) \, \delta \frac{\partial{v}}{\partial{t}}  \, dt \, \Big|_{s=L}
\\ \nonumber
& 
+ J \int_{t_1}^{t_2} \left( \frac{\partial^2{v}}{\partial{t}\partial{s}} ( 1 + (\frac{\partial{v}}{\partial{s}})^2) \delta \frac{\partial^2{v}}{\partial{t}\partial{s}} + \frac{\partial{v}}{\partial{s}} (\frac{\partial^2{v}}{\partial{t}\partial{s}})^2 \delta \frac{\partial{v}}{\partial{s}} \right)  \, dt \, \Big|_{s=L}
\\ \nonumber
= &
\rho \!\int_{0}^{L}\! \left[ (\frac{\partial{v}}{\partial{t}} \!+\! \dot{V_b})  \delta{v} \Big|_{t_1}^{t_2} \!-\! \!\int_{t_1}^{t_2}\!  (\frac{\partial^2{v}}{\partial{t}^2} \!+\! \ddot{V_b})  \delta {v}  dt \right]  ds
\!+\! M (\frac{\partial{v}}{\partial{t}} \!+\! \dot{V_b})  \delta v \Big|_{s=L} \Big|_{t_1}^{t_2} 
\!-\! M \!\int_{t_1}^{t_2} (\frac{\partial^2{v}}{\partial{t}^2} \!+\! \ddot{V_b})  \delta v   dt  \Big|_{s=L}
\\ \nonumber
& \quad
\!+\! J  \frac{\partial^2{v}}{\partial{t}\partial{s}} ( 1 \!+\! (\frac{\partial{v}}{\partial{s}})^2) \delta \frac{\partial{v}}{\partial{s}}  \Big|_{s=L}  \Big|_{t_1}^{t_2}
\!-\! J \int_{t_1}^{t_2} 
\left(  \frac{\partial^3{v}}{\partial{t}^2\partial{s}} ( 1 \!+\! (\frac{\partial{v}}{\partial{s}})^2) \!+\! \frac{\partial{v}}{\partial{s}} (\frac{\partial^2{v}}{\partial{t}\partial{s}})^2 \right)  
\delta \frac{\partial{v}}{\partial{s}}   dt  \Big|_{s=L}
\\ \nonumber
\!=\! &
- \!\int_{t_1}^{t_2}\! \Bigg\{
\rho \int_{0}^{L}  (\frac{\partial^2{v}}{\partial{t}^2} \!+\! \ddot{V_b})  \delta {v}  ds 
\!+\! M (\frac{\partial^2{v}}{\partial{t}^2} \!+\! \ddot{V_b})  \delta v   \Big|_{s=L}
\\ \nonumber
&+ J \left(  \frac{\partial^3{v}}{\partial{t}^2\partial{s}}( 1 \!+\! (\frac{\partial{v}}{\partial{s}})^2) \!+\! \frac{\partial{v}}{\partial{s}} (\frac{\partial^2{v}}{\partial{t}\partial{s}})^2  \right)  \delta \frac{\partial{v}}{\partial{s}}   \Big|_{s=L}
\Bigg\}  dt,
\end{align}
where we consider that $\delta v = \delta \frac{\partial{v}}{\partial{s}} = 0$ at $t=t_1$ and $t=t_2$. Therefore, the extended Hamilton's principle takes the form
\begin{align}
\label{Eq: extended Hamilton}
&
\int_{t_1}^{t_2} \Bigg\{
\int_{0}^{L} \Bigg[
- \rho (\frac{\partial^2{v}}{\partial{t}^2}  + \ddot{V_b}) 
\\ \nonumber
&
-\frac{\partial^2}{\partial{s}^2} \left(
\left( E_{\infty} \, I \, \frac{\partial^2{v}}{\partial{s}^2} (1 + \frac{1}{2} (\frac{\partial{v}}{\partial{s}})^2)  + E_{\alpha} \, I \, \prescript{RL}{0}{\mathcal D}_{t}^{\alpha} \, \Big[\frac{\partial^2{v}}{\partial{s}^2}(1 + \frac{1}{2} (\frac{\partial{v}}{\partial{s}})^2)\Big]\right) 
(1 + \frac{1}{2} (\frac{\partial{v}}{\partial{s}})^2) 
\right)
\\ \nonumber
&
+\frac{\partial}{\partial{s}}\left( 
\left( E_{\infty} \, I \, \frac{\partial^2{v}}{\partial{s}^2} (1 + \frac{1}{2} (\frac{\partial{v}}{\partial{s}})^2)  + E_{\alpha} \, I \, \prescript{RL}{0}{\mathcal D}_{t}^{\alpha} \, \Big[\frac{\partial^2{v}}{\partial{s}^2}(1 + \frac{1}{2} (\frac{\partial{v}}{\partial{s}})^2)\Big]\right) 
\frac{\partial^2{v}}{\partial{s}^2} \, \frac{\partial{v}}{\partial{s}}\, 
\right)
\Bigg] \, \delta {v} \, ds
\\ \nonumber
& 
- M (\frac{\partial^2{v}}{\partial{t}^2} + \ddot{V_b}) \Big|_{s=L} \, \delta v(L,t) 
- J \left(  \frac{\partial^3{v}}{\partial{t}^2\partial{s}} ( 1 + (\frac{\partial{v}}{\partial{s}})^2) + \frac{\partial{v}}{\partial{s}} (\frac{\partial^2{v}}{\partial{t}\partial{s}})^2 \right) \Big|_{s=L} \,  \delta \frac{\partial{v}}{\partial{s}}(L,t) 
\\ \nonumber
& -
\left( E_{\infty} \, I \, \frac{\partial^2{v}}{\partial{s}^2} (1 + \frac{1}{2} (\frac{\partial{v}}{\partial{s}})^2)  + E_{\alpha} \, I \, \prescript{RL}{0}{\mathcal D}_{t}^{\alpha} \, \Big[\frac{\partial^2{v}}{\partial{s}^2}(1 + \frac{1}{2} (\frac{\partial{v}}{\partial{s}})^2)\Big]\right) 
(1 + \frac{1}{2} (\frac{\partial{v}}{\partial{s}})^2) \Bigg|_{s=L} \, \delta \frac{\partial{v}}{\partial{s}}(L,t) 
\\ \nonumber
&  +
\frac{\partial}{\partial{s}}\left(
\left( E_{\infty} \, I \, \frac{\partial^2{v}}{\partial{s}^2} (1 + \frac{1}{2} (\frac{\partial{v}}{\partial{s}})^2)  + E_{\alpha} \, I \, \prescript{RL}{0}{\mathcal D}_{t}^{\alpha} \, \Big[\frac{\partial^2{v}}{\partial{s}^2}(1 + \frac{1}{2} (\frac{\partial{v}}{\partial{s}})^2)\Big]\right) 
(1 + \frac{1}{2} (\frac{\partial{v}}{\partial{s}})^2) \right) \Bigg|_{s=L} \, \delta v(L,t) 
\\ \nonumber
&  -
\left( E_{\infty} \, I \, \frac{\partial^2{v}}{\partial{s}^2} (1 + \frac{1}{2} (\frac{\partial{v}}{\partial{s}})^2)  + E_{\alpha} \, I \, \prescript{RL}{0}{\mathcal D}_{t}^{\alpha} \, \Big[\frac{\partial^2{v}}{\partial{s}^2}(1 + \frac{1}{2} (\frac{\partial{v}}{\partial{s}})^2)\Big]\right) 
\frac{\partial^2{v}}{\partial{s}^2} \, \frac{\partial{v}}{\partial{s}} \Bigg|_{s=L} \, \delta v(L,t) 
\,\, \Bigg\} \, dt = 0.
\end{align}
Invoking the arbitrariness of virtual displacement $\delta v$, we obtain the strong form of the equation of motion as:
\begin{align}
\label{Eq: eqn of motion}
& \rho \, \frac{\partial^2{v}}{\partial{t}^2}  
+ 
E_{\infty} \, I \, \frac{\partial^2}{\partial{s}^2}\left(  \frac{\partial^2{v}}{\partial{s}^2} (1 + \frac{1}{2} (\frac{\partial{v}}{\partial{s}})^2)^2 \right)
+ 
E_{\alpha} \, I \, \frac{\partial^2}{\partial{s}^2}\left(  (1 + \frac{1}{2} (\frac{\partial{v}}{\partial{s}})^2) \, \prescript{RL}{0}{\mathcal D}_{t}^{\alpha} \, \Big[\frac{\partial^2{v}}{\partial{s}^2}(1 + \frac{1}{2} (\frac{\partial{v}}{\partial{s}})^2) \Big]\right)
\\ \nonumber
& - 
E_{\infty} \, I \, \frac{\partial}{\partial{s}}\left(  \frac{\partial{v}}{\partial{s}} \,  (\frac{\partial^2{v}}{\partial{s}^2})^2 (1 + \frac{1}{2} (\frac{\partial{v}}{\partial{s}})^2)  \right)
- 
E_{\alpha} \, I \, \frac{\partial}{\partial{s}} \left(  \frac{\partial{v}}{\partial{s}} \, \frac{\partial^2{v}}{\partial{s}^2} \,  \prescript{RL}{0}{\mathcal D}_{t}^{\alpha} \, \Big[\frac{\partial^2{v}}{\partial{s}^2}(1 + \frac{1}{2} (\frac{\partial{v}}{\partial{s}})^2)\Big] \right) 
= -\rho \,\ddot{V_b} ,
\end{align}
which is subject to the following natural boundary conditions:
\begin{align}
\label{Eq: natural bc}
& 
J \left(  \frac{\partial^3{v}}{\partial{t}^2\partial{s}}  ( 1 + (\frac{\partial{v}}{\partial{s}})^2) + \frac{\partial{v}}{\partial{s}} (\frac{\partial^2{v}}{\partial{t}\partial{s}})^2   \right)
\\ \nonumber
& \qquad
+ E_{\infty} \, I \, \frac{\partial^2{v}}{\partial{s}^2} (1 + \frac{1}{2} (\frac{\partial{v}}{\partial{s}})^2)^2  + E_{\alpha} \, I \, (1 + \frac{1}{2} (\frac{\partial{v}}{\partial{s}})^2) \, \prescript{RL}{0}{\mathcal D}_{t}^{\alpha} \, \Big[\frac{\partial^2{v}}{\partial{s}^2}(1 + \frac{1}{2} (\frac{\partial{v}}{\partial{s}})^2)\Big]
\,\, \Bigg|_{s=L} = 0 ,
\\ \nonumber
& 
M (\frac{\partial^2{v}}{\partial{t}^2}  + \ddot{V_b})
- \frac{\partial}{\partial{s}}\left(
E_{\infty} \, I \, \frac{\partial^2{v}}{\partial{s}^2} (1 + \frac{1}{2} (\frac{\partial{v}}{\partial{s}})^2)^2  + E_{\alpha} \, I \, (1 + \frac{1}{2} (\frac{\partial{v}}{\partial{s}})^2) \, \prescript{RL}{0}{\mathcal D}_{t}^{\alpha} \, \Big[\frac{\partial^2{v}}{\partial{s}^2}(1 + \frac{1}{2} (\frac{\partial{v}}{\partial{s}})^2)\Big] \right)
\\ \nonumber
& \qquad 
+ \left( E_{\infty} \, I \, \frac{\partial{v}}{\partial{s}} \, (\frac{\partial^2{v}}{\partial{s}^2})^2 (1 + \frac{1}{2} (\frac{\partial{v}}{\partial{s}})^2)  + E_{\alpha} \, I \, \frac{\partial{v}}{\partial{s}} \, \frac{\partial^2{v}}{\partial{s}^2}  \, \prescript{RL}{0}{\mathcal D}_{t}^{\alpha} \, \Big[\frac{\partial^2{v}}{\partial{s}^2}(1 + \frac{1}{2} (\frac{\partial{v}}{\partial{s}})^2) \Big]\right) 
\,\, \Bigg|_{s=L}  = 0 .
\end{align}
Following a similar approach as in \eqref{Eq: curvature} in deriving the beam curvature, we obtain the approximations below, where we only consider up to third order terms and remove the higher order terms (HOTs). 
\begin{align*}
\frac{\partial^2{v}}{\partial{s}^2} (1 + \frac{1}{2} (\frac{\partial{v}}{\partial{s}})^2)^2 
& = 
\frac{\partial^2{v}}{\partial{s}^2} + \frac{\partial^2{v}}{\partial{s}^2} (\frac{\partial{v}}{\partial{s}})^2+ \text{HOTs}
\\ 
(1 + \frac{1}{2} (\frac{\partial{v}}{\partial{s}})^2) \, \prescript{RL}{0}{\mathcal D}_{t}^{\alpha} \, \Big[\frac{\partial^2{v}}{\partial{s}^2}(1 + \frac{1}{2} (\frac{\partial{v}}{\partial{s}})^2)\Big] 
& = 
\prescript{RL}{0}{\mathcal D}_{t}^{\alpha}  \Big[\frac{\partial^2{v}}{\partial{s}^2}(1 \!+\! \frac{1}{2} (\frac{\partial{v}}{\partial{s}})^2)\Big] \!+\!\frac{1}{2} (\frac{\partial{v}}{\partial{s}})^2  \prescript{RL}{0}{\mathcal D}_{t}^{\alpha}  \frac{\partial^2{v}}{\partial{s}^2} \!+\! \text{HOTs}
\\
\frac{\partial{v}}{\partial{s}} \,  (\frac{\partial^2{v}}{\partial{s}^2})^2 (1 + \frac{1}{2} (\frac{\partial{v}}{\partial{s}})^2) 
& = \frac{\partial{v}}{\partial{s}} \,  (\frac{\partial^2{v}}{\partial{s}^2})^2 + \text{HOTs}
\\
\frac{\partial{v}}{\partial{s}} \, \frac{\partial^2{v}}{\partial{s}^2} \,  \prescript{RL}{0}{\mathcal D}_{t}^{\alpha} \, \Big[\frac{\partial^2{v}}{\partial{s}^2}(1 + \frac{1}{2} (\frac{\partial{v}}{\partial{s}})^2)\Big] 
& = \frac{\partial{v}}{\partial{s}}\, \frac{\partial^2{v}}{\partial{s}^2} \,  \prescript{RL}{0}{\mathcal D}_{t}^{\alpha} \, \frac{\partial^2{v}}{\partial{s}^2} + \text{HOTs}
\end{align*}
Therefore, the strong form can be approximated up to the third order and the problem then reads as: find $v \in V$ such that
\begin{align}
\label{Eq: eqn of motion - 1}
& m \, \frac{\partial^2{v}}{\partial{t}^2}  
+ 
\frac{\partial^2}{\partial{s}^2}\left(  \frac{\partial^2{v}}{\partial{s}^2} + \frac{\partial^2{v}}{\partial{s}^2} (\frac{\partial{v}}{\partial{s}})^2 \right)
- 
\frac{\partial}{\partial{s}}\left(  \frac{\partial{v}}{\partial{s}} \,  (\frac{\partial^2{v}}{\partial{s}^2})^2  \right)
\\ \nonumber
&
+ 
E_r  \frac{\partial^2}{\partial{s}^2}\left(  \prescript{RL}{0}{\mathcal D}_{t}^{\alpha}  \Big[\frac{\partial^2{v}}{\partial{s}^2}(1 + \frac{1}{2} (\frac{\partial{v}}{\partial{s}})^2)\Big] + \frac{1}{2}  (\frac{\partial{v}}{\partial{s}})^2  \prescript{RL}{0}{\mathcal D}_{t}^{\alpha}  \frac{\partial^2{v}}{\partial{s}^2} \right)
- 
E_r  \frac{\partial}{\partial{s}}\left(  \frac{\partial{v}}{\partial{s}} \frac{\partial^2{v}}{\partial{s}^2}   \prescript{RL}{0}{\mathcal D}_{t}^{\alpha} \frac{\partial^2{v}}{\partial{s}^2} \right) 
= -m \ddot{V_b} ,
\end{align}
\begin{align}
\label{Eq: eqn of motion - 2}
& m \, \frac{\partial^2{v}}{\partial{t}^2} 
+ 
\frac{\partial^2}{\partial{s}^2}\left(  
\frac{\partial^2{v}}{\partial{s}^2} + \frac{\partial^2{v}}{\partial{s}^2} (\frac{\partial{v}}{\partial{s}})^2 
+ E_r
\prescript{RL}{0}{\mathcal D}_{t}^{\alpha} \, \Big[\frac{\partial^2{v}}{\partial{s}^2}(1 + \frac{1}{2} (\frac{\partial{v}}{\partial{s}})^2)\Big] 
+  \frac{1}{2}  E_r 
(\frac{\partial{v}}{\partial{s}})^2 \, \prescript{RL}{0}{\mathcal D}_{t}^{\alpha} \, \frac{\partial^2{v}}{\partial{s}^2}
\right)
\\ \nonumber
& 
- 
\frac{\partial}{\partial{s}} \left( 
\frac{\partial{v}}{\partial{s}} \,  (\frac{\partial^2{v}}{\partial{s}^2})^2  
+ E_r
\frac{\partial{v}}{\partial{s}} \, \frac{\partial^2{v}}{\partial{s}^2} \,  \prescript{RL}{0}{\mathcal D}_{t}^{\alpha} \, \frac{\partial^2{v}}{\partial{s}^2}
\right)
= -m \,\ddot{V_b} ,
\end{align}
subject to the following boundary conditions:
\begin{align}
\label{Eq: bc}
& v  \Big|_{s=0} = \frac{\partial{v}}{\partial{s}}  \Big|_{s=0} = 0 ,
\\ \nonumber
& 
\frac{J m}{\rho} \left(  \frac{\partial^3{v}}{\partial{t}^2\partial{s}} ( 1 + (\frac{\partial{v}}{\partial{s}})^2) + \frac{\partial{v}}{\partial{s}} (\frac{\partial^2{v}}{\partial{t}\partial{s}})^2  \right)
+
\\ \nonumber
& 
\left( 
\frac{\partial^2{v}}{\partial{s}^2} + \frac{\partial^2{v}}{\partial{s}^2} (\frac{\partial{v}}{\partial{s}})^2
+ E_r  \prescript{RL}{0}{\mathcal D}_{t}^{\alpha}  \Big[\frac{\partial^2{v}}{\partial{s}^2}(1 + \frac{1}{2} (\frac{\partial{v}}{\partial{s}})^2)\Big] 
+ \frac{1}{2}  E_r  (\frac{\partial{v}}{\partial{s}})^2  \prescript{RL}{0}{\mathcal D}_{t}^{\alpha}  \frac{\partial^2{v}}{\partial{s}^2}  
\right)  \Bigg|_{s=L} = 0 ,
\\ \nonumber
& 
\frac{M m}{\rho} (\frac{\partial^2{v}}{\partial{t}^2} \!+\! \ddot{V_b})
\!-\! \frac{\partial}{\partial{s}} \left(
\frac{\partial^2{v}}{\partial{s}^2} \!+\! \frac{\partial^2{v}}{\partial{s}^2} (\frac{\partial{v}}{\partial{s}})^2
\!+\! E_r  \prescript{RL}{0}{\mathcal D}_{t}^{\alpha}  \Big[\frac{\partial^2{v}}{\partial{s}^2}(1 \!+\! \frac{1}{2} (\frac{\partial{v}}{\partial{s}})^2) \Big]
\!+\! \frac{1}{2}  E_r  (\frac{\partial{v}}{\partial{s}})^2  \prescript{RL}{0}{\mathcal D}_{t}^{\alpha}  \frac{\partial^2{v}}{\partial{s}^2}  
\right)
\\ \nonumber
& 
+ \left( 
\frac{\partial{v}}{\partial{s}}  (\frac{\partial^2{v}}{\partial{s}^2})^2  + E_r  \frac{\partial{v}}{\partial{s}}   \frac{\partial^2{v}}{\partial{s}^2}   \prescript{RL}{0}{\mathcal D}_{t}^{\alpha}  \frac{\partial^2{v}}{\partial{s}^2}	
\right)  \Bigg|_{s=L}  = 0 ,
\end{align}
where $m = \frac{\rho}{ E_{\infty} \, I}$ and $E_r = \frac{E_{\alpha}}{E_{\infty}}$. 

%
%%%%%%%%%%%%%%%%%%%%%%%%%%%
\subsection{Nondimensionalization}
%%%%%%%%%%%%%%%%%%%%%%%%%%%
%
Let the dimensionless variables
\begin{align}
s^{*} \!=\! \frac{s}{L}, \quad 
v^{*} \!=\! \frac{v}{L}, \quad 
t^{*} \!=\! t \left(\frac{1}{m L^4}\right)^{1/2}, \quad 
E_r^{*}\! =\! E_r \left(\frac{1}{m L^4}\right)^{\alpha/2}, \quad 
J^{*} \!=\! \frac{J}{\rho L^3}, \quad
M^{*} \!=\! \frac{M}{\rho L}, \quad
V_b^* \!=\!  \frac{V_b}{L}.
\end{align}
We obtain the following dimensionless equation by substituting the above dimensionless variables. 
\begin{align}
\label{Eq: eqn of motion - dimless - 0}
& m \frac{L}{mL^4} \, \frac{\partial^2 v^*}{\partial {t^*}^2}  
\\ \nonumber
& + 
\frac{1}{L^2} \frac{\partial^2}{\partial {s^*}^2} 
\Bigg[  
\frac{L}{L^2} \frac{\partial^2 v^*}{\partial {s^*}^2} 
+ \frac{L}{L^2} \frac{\partial^2 v^*}{\partial {s^*}^2} 
( \frac{L}{L} \frac{\partial v^*}{\partial {s^*}} )^2 
+ \frac{E_r^* (mL^4)^{\alpha/2}}{2} \frac{1}{(mL^4)^{\alpha/2}} \frac{L}{L^2} (\frac{L}{L})^2 \prescript{RL}{0}{\mathcal D}_{t^*}^{\alpha} \frac{\partial^2 v^*}{\partial {s^*}^2}(\frac{\partial v^*}{\partial {s^*}})^2 
\\ \nonumber
&
+ E_r^* (mL^4)^{\alpha/2} \frac{1}{(mL^4)^{\alpha/2}} \frac{L}{L^2} \prescript{RL}{0}{\mathcal D}_{t^*}^{\alpha} \frac{\partial^2 v^*}{\partial {s^*}^2}
+ \frac{1}{2}  E_r^* (mL^4)^{\alpha/2} (\frac{L}{L} \frac{\partial v^*}{\partial s^*})^2  \frac{1}{(mL^4)^{\alpha/2}} \frac{L}{L^2} \prescript{RL}{0}{\mathcal D}_{t^*}^{\alpha} \frac{\partial^2 v^*}{\partial {s^*}^2}
\Bigg]
\\ \nonumber
& - 
\frac{1}{L} \frac{\partial}{\partial {s^*}} 
\Bigg[
\frac{L}{L} \frac{\partial v^*}{\partial {s^*}} (\frac{L}{L^2} \frac{\partial^2 v^*}{\partial {s^*}^2})^2 
+ E_r^* (mL^4)^{\alpha/2} \frac{L}{L} \frac{\partial v^*}{\partial s^*} \frac{L}{L^2} \frac{\partial^2 v^*}{\partial {s^*}^2}  \frac{1}{(mL^4)^{\alpha/2}} \frac{L}{L^2} \prescript{RL}{0}{\mathcal D}_{t^*}^{\alpha} \frac{\partial^2 v^*}{\partial {s^*}^2}
\Bigg]
\\ \nonumber
&
= -m \frac{L}{mL^4} \, \frac{\partial^2 V_b^*}{\partial {t^*}^2}  ,
\end{align}
which can be simplified to
\begin{align}
\label{Eq: eqn of motion - dimless - 1}
\frac{\partial^2 v^*}{\partial {t^*}^2}  
& \!+\! 
\frac{\partial^2}{\partial {s^*}^2} 
\Bigg[  
\frac{\partial^2 v^*}{\partial {s^*}^2} 
\!+\! \frac{\partial^2 v^*}{\partial {s^*}^2} (\frac{\partial v^*}{\partial {s^*}} )^2 
\!+\! \frac{E_r^*}{2} \prescript{RL}{0}{\mathcal D}_{t^*}^{\alpha} \frac{\partial^2 v^*}{\partial {s^*}^2}(\frac{\partial v^*}{\partial {s^*}})^2 
\!+\! E_r^* \prescript{RL}{0}{\mathcal D}_{t^*}^{\alpha} \frac{\partial^2 v^*}{\partial {s^*}^2}
\\ \nonumber
&
\!+\! \frac{1}{2}  E_r^* (\frac{\partial v^*}{\partial s^*})^2 \prescript{RL}{0}{\mathcal D}_{t^*}^{\alpha} \frac{\partial^2 v^*}{\partial {s^*}^2}
\Bigg]
\!-\!
\frac{\partial}{\partial {s^*}} 
\Bigg[
\frac{\partial v^*}{\partial {s^*}} (\frac{\partial^2 v^*}{\partial {s^*}^2})^2 
\!+\! E_r^* \frac{\partial v^*}{\partial s^*} \frac{\partial^2 v^*}{\partial {s^*}^2}  \prescript{RL}{0}{\mathcal D}_{t^*}^{\alpha} \frac{\partial^2 v^*}{\partial {s^*}^2}
\Bigg]
= -\frac{\partial^2 V_b^*}{\partial {t^*}^2}  ,
\end{align}
The dimensionless boundary conditions are also obtained by substituting dimensionless variables in \eqref{Eq: bc}. We can show similarly that they preserve their structure as: 
\begin{align*}
& v^{*} \, \Big|_{s^{*}=0} = \frac{\partial v^{*}}{\partial s^{*}} \, \Big|_{s^{*}=0} = 0 ,
\\ \nonumber
& 
\frac{J^* \rho L^3 m}{\rho} \frac{1}{mL^4}
\Bigg[ 
\frac{\partial^3 v^{*}}{\partial t^{*} \partial^2 s^{*}} \left( 1 + \left(\frac{\partial v^{*}}{\partial s^{*}}\right)^2 \right) 
+\frac{\partial v^{*}}{\partial s^{*}} \left(\frac{\partial^2 v^{*}}{\partial t^{*} \partial s^{*}}\right)^2  
\Bigg]
\!+\!
\frac{1}{L} \Bigg[  
\frac{\partial^2 v^*}{\partial {s^*}^2} 
\!+\! \frac{\partial^2 v^*}{\partial {s^*}^2} (\frac{\partial v^*}{\partial {s^*}} )^2 
\\ \nonumber
&
\!+\! \frac{E_r^*}{2} \prescript{RL}{0}{\mathcal D}_{t^*}^{\alpha} \frac{\partial^2 v^*}{\partial {s^*}^2}(\frac{\partial v^*}{\partial {s^*}})^2 
\!+\! E_r^* \prescript{RL}{0}{\mathcal D}_{t^*}^{\alpha} \frac{\partial^2 v^*}{\partial {s^*}^2}
\!\!+ \frac{1}{2}  E_r^* (\frac{\partial v^*}{\partial s^*})^2 \prescript{RL}{0}{\mathcal D}_{t^*}^{\alpha} \frac{\partial^2 v^*}{\partial {s^*}^2}
\Bigg]
\Bigg|_{s^{*}=1} \!=\! 0 ,
\\ \nonumber
& 
\frac{M^* \rho L m}{\rho} \frac{L}{mL^4} \left( \frac{\partial^2 v^{*}}{\partial^2 t^{*}}  + \frac{\partial^2 V_b^{*}}{\partial^2 t^{*}}  \right)
\!-\!\frac{1}{L^2} \frac{\partial v^{*}}{\partial s^{*}} \Bigg[  
\frac{\partial^2 v^*}{\partial {s^*}^2} 
\!+\! \frac{\partial^2 v^*}{\partial {s^*}^2} (\frac{\partial v^*}{\partial {s^*}} )^2 
\!+\! \frac{E_r^*}{2} \prescript{RL}{0}{\mathcal D}_{t^*}^{\alpha} \frac{\partial^2 v^*}{\partial {s^*}^2}(\frac{\partial v^*}{\partial {s^*}})^2 
\\ \nonumber
& 
\!+\! E_r^* \prescript{RL}{0}{\mathcal D}_{t^*}^{\alpha} \frac{\partial^2 v^*}{\partial {s^*}^2}
\!+\! \frac{1}{2}  E_r^* (\frac{\partial v^*}{\partial s^*})^2 \prescript{RL}{0}{\mathcal D}_{t^*}^{\alpha} \frac{\partial^2 v^*}{\partial {s^*}^2}
\Bigg]
\!+\!\frac{1}{L^2} 
\Bigg[
\frac{\partial v^*}{\partial {s^*}} (\frac{\partial^2 v^*}{\partial {s^*}^2})^2 
\!+\! E_r^* \frac{\partial v^*}{\partial s^*} \frac{\partial^2 v^*}{\partial {s^*}^2}  \prescript{RL}{0}{\mathcal D}_{t^*}^{\alpha} \frac{\partial^2 v^*}{\partial {s^*}^2}
\Bigg] \Bigg|_{s^{*}=1}  \!=\! 0 ,
\end{align*}
Therefore, the dimensionless equation of motion becomes (after dropping $^*$ for the sake of simplicity)
\begin{align}
\label{Eq: eqn of motion - dimless - 2}
\frac{\partial^2{v}}{\partial{t}^2} 
& + 
\frac{\partial^2}{\partial{s}^2} \left(  
\frac{\partial^2{v}}{\partial{s}^2} + \frac{\partial^2{v}}{\partial{s}^2} (\frac{\partial{v}}{\partial{s}})^2
+ E_r \prescript{RL}{0}{\mathcal D}_{t}^{\alpha} \, \Big[\frac{\partial^2{v}}{\partial{s}^2}(1 + \frac{1}{2} (\frac{\partial{v}}{\partial{s}})^2)\Big] 
+ \frac{1}{2}  E_r (\frac{\partial{v}}{\partial{s}})^2 \, \prescript{RL}{0}{\mathcal D}_{t}^{\alpha} \, \frac{\partial^2{v}}{\partial{s}^2}
\right)
\\ \nonumber
& - 
\frac{\partial}{\partial{s}} \left( 
\frac{\partial{v}}{\partial{s}}  \,  (\frac{\partial^2{v}}{\partial{s}^2})^2  
+ E_r \frac{\partial{v}}{\partial{s}}  \, \frac{\partial^2{v}}{\partial{s}^2} \,  \prescript{RL}{0}{\mathcal D}_{t}^{\alpha} \, \frac{\partial^2{v}}{\partial{s}^2}
\right)
= -\ddot{V_b} ,
\end{align}
which is subject to the following dimensionless boundary conditions
\begin{align}
\label{Eq: bc - dimless}
& v \, \Big|_{s=0} = \frac{\partial{v}}{\partial{s}} \, \Big|_{s=0} = 0 ,
\\ \nonumber
& 
J \left(  \frac{\partial^3{v}}{\partial{t}^2\partial{s}}  ( 1 + (\frac{\partial{v}}{\partial{s}})^2) + \frac{\partial{v}}{\partial{s}}  (\frac{\partial^2{v}}{\partial{s}\partial{t}})^2   \right)
+
\\ \nonumber
& 
\left( 
\frac{\partial^2{v}}{\partial{s}^2} + \frac{\partial^2{v}}{\partial{s}^2} (\frac{\partial{v}}{\partial{s}})^2
+ E_r \, \prescript{RL}{0}{\mathcal D}_{t}^{\alpha} \, \Big[\frac{\partial^2{v}}{\partial{s}^2}(1 + \frac{1}{2} (\frac{\partial{v}}{\partial{s}})^2)\Big] 
+ \frac{1}{2}  E_r \, (\frac{\partial{v}}{\partial{s}})^2 \, \prescript{RL}{0}{\mathcal D}_{t}^{\alpha} \, \frac{\partial^2{v}}{\partial{s}^2}  
\right) \,\, \Bigg|_{s=1} = 0 ,
\\ \nonumber
& 
M (\frac{\partial^2{v}}{\partial{t}^2} + \ddot{V_b})
- \frac{\partial}{\partial{s}} \left(
\frac{\partial^2{v}}{\partial{s}^2} + \frac{\partial^2{v}}{\partial{s}^2} (\frac{\partial{v}}{\partial{s}})^2
+ E_r \, \prescript{RL}{0}{\mathcal D}_{t}^{\alpha} \, \Big[\frac{\partial^2{v}}{\partial{s}^2}(1 + \frac{1}{2} (\frac{\partial{v}}{\partial{s}})^2) \Big]
+ \frac{1}{2}  E_r \, (\frac{\partial{v}}{\partial{s}})^2 \, \prescript{RL}{0}{\mathcal D}_{t}^{\alpha} \, \frac{\partial^2{v}}{\partial{s}^2}  
\right)
\\ \nonumber
& 
+ \left( 
\frac{\partial{v}}{\partial{s}} \, (\frac{\partial^2{v}}{\partial{s}^2})^2  + E_r \, \frac{\partial{v}}{\partial{s}} \, \frac{\partial^2{v}}{\partial{s}^2}  \, \prescript{RL}{0}{\mathcal D}_{t}^{\alpha} \, \frac{\partial^2{v}}{\partial{s}^2}	
\right) \,\, \Bigg|_{s=1}  = 0 ,
\end{align}
% 
%%%%%%%%%%%%%%%%%%%%%%%%%%%%%%%%%%%%%%%%%%
\section{Single Mode Decomposition}
\label{Sec: App. Single Mode Decomposition} 
%%%%%%%%%%%%%%%%%%%%%%%%%%%%%%%%%%%%%%%%%
%
To show the single mode decomposition satisfies the weak form solution and the boundary conditions, we consider the case of no lumped mass at tip where $M=J=0$ and check if the proposed approximate solution will solve the weak form and its subject boundary conditions. First we substitute the boundary conditions in the weak formulation and then we use approximation to recover the equation \eqref{Eq: weak form - discrete 2} .\\
We use the weak formulation of the problem in \eqref{Eq: weak form - 1} and do the integration by part to transfer the spatial derivative load to the test function as follow 
\begin{align}
\label{Eq: weak form single mode - 1}
& \int_{0}^{1} \frac{\partial^2{v}}{\partial{t}^2}  \tilde{v}  ds  
+ 
\frac{\partial}{\partial{s}}\left(  
\frac{\partial^2{v}}{\partial{s}^2} + \frac{\partial^2{v}}{\partial{s}^2} (\frac{\partial{v}}{\partial{s}})^2 
+ E_r
\prescript{RL}{0}{\mathcal D}_{t}^{\alpha}  \Big[\frac{\partial^2{v}}{\partial{s}^2}(1 + \frac{1}{2} (\frac{\partial{v}}{\partial{s}})^2) \Big]
+ \frac{1}{2}  E_r 
(\frac{\partial{v}}{\partial{s}})^2  \prescript{RL}{0}{\mathcal D}_{t}^{\alpha}  \frac{\partial^2{v}}{\partial{s}^2}
\right)
\tilde{v}  \Bigg|_{0}^{1}
\\ \nonumber
&
-\left(  
\frac{\partial^2{v}}{\partial{s}^2} + \frac{\partial^2{v}}{\partial{s}^2} (\frac{\partial{v}}{\partial{s}})^2 
+ E_r
\prescript{RL}{0}{\mathcal D}_{t}^{\alpha}  \Big[\frac{\partial^2{v}}{\partial{s}^2}(1 + \frac{1}{2} (\frac{\partial{v}}{\partial{s}})^2) \Big]
+ \frac{1}{2}  E_r 
(\frac{\partial{v}}{\partial{s}})^2  \prescript{RL}{0}{\mathcal D}_{t}^{\alpha}  \frac{\partial^2{v}}{\partial{s}^2}
\right)
\frac{\partial\tilde{v}}{\partial{s}}  \Bigg|_{0}^{1}
\\ \nonumber
&
+ \int_{0}^{1} \left(  
\frac{\partial^2{v}}{\partial{s}^2} + \frac{\partial^2{v}}{\partial{s}^2} (\frac{\partial{v}}{\partial{s}})^2 
+ E_r
\prescript{RL}{0}{\mathcal D}_{t}^{\alpha}  \Big[\frac{\partial^2{v}}{\partial{s}^2}(1 + \frac{1}{2} (\frac{\partial{v}}{\partial{s}})^2) \Big]
+ \frac{1}{2}  E_r 
(\frac{\partial{v}}{\partial{s}})^2  \prescript{RL}{0}{\mathcal D}_{t}^{\alpha}  \frac{\partial^2{v}}{\partial{s}^2}
\right)
\frac{\partial^2{\tilde{v}}}{\partial{s}^2}  ds
\\ \nonumber
&
- \left( 
\frac{\partial{v}}{\partial{s}}   (\frac{\partial^2{v}}{\partial{s}^2})^2  
+ E_r
\frac{\partial{v}}{\partial{s}}  \frac{\partial^2{v}}{\partial{s}^2}   \prescript{RL}{0}{\mathcal D}_{t}^{\alpha}  \frac{\partial^2{v}}{\partial{s}^2}
\right)
\tilde{v}  \Bigg|_{0}^{1}
\\ \nonumber
&
+\int_{0}^{1}\left( 
\frac{\partial{v}}{\partial{s}}   (\frac{\partial^2{v}}{\partial{s}^2})^2  
+ E_r
\frac{\partial{v}}{\partial{s}}  \frac{\partial^2{v}}{\partial{s}^2}   \prescript{RL}{0}{\mathcal D}_{t}^{\alpha}  \frac{\partial^2{v}}{\partial{s}^2}
\right)
\frac{\partial{\tilde{v}}}{\partial{s}}  ds
= f(t) ,
\end{align}
When $M=J=0$ our boundary conditions in equation \eqref{Eq: bc - dimless - body} will be
\begin{align}
\label{Eq: weak form single mode - 2}
& v \, \Big|_{s=0} = \frac{\partial{v}}{\partial{s}} \, \Big|_{s=0} = 0 ,
\\ \nonumber
& 
\left( 
\frac{\partial^2{v}}{\partial{s}^2}+  \frac{\partial^2{v}}{\partial{s}^2} ( \frac{\partial{v}}{\partial{s}})^2
+ E_r \, \prescript{RL}{0}{\mathcal D}_{t}^{\alpha} \,  \Big[\frac{\partial^2{v}}{\partial{s}^2}(1 + \frac{1}{2}  (\frac{\partial{v}}{\partial{s}})^2) \Big]
+ \frac{1}{2}  E_r \,  (\frac{\partial{v}}{\partial{s}})^2 \, \prescript{RL}{0}{\mathcal D}_{t}^{\alpha} \,  \frac{\partial^2{v}}{\partial^2{s}}
\right) \,\, \Bigg|_{s=1} = 0 ,
\\ \nonumber
& 
- \frac{\partial{v}}{\partial{s}}\left(
\frac{\partial^2{v}}{\partial{s}^2}+ \frac{\partial^2{v}}{\partial{s}^2} (\frac{\partial{v}}{\partial{s}})^2
+ E_r \, \prescript{RL}{0}{\mathcal D}_{t}^{\alpha} \, \Big[\frac{\partial^2{v}}{\partial{s}^2}(1 + \frac{1}{2} (\frac{\partial{v}}{\partial{s}})^2)\Big]
+ \frac{1}{2}  E_r \, (\frac{\partial{v}}{\partial{s}})^2 \, \prescript{RL}{0}{\mathcal D}_{t}^{\alpha} \, \frac{\partial^2{v}}{\partial{s}^2}
\right)
\\ \nonumber
& 
+ \left( 
\frac{\partial{v}}{\partial{s}} \, (\frac{\partial^2{v}}{\partial{s}^2})^2  + E_r \, \frac{\partial{v}}{\partial{s}} \, \frac{\partial^2{v}}{\partial{s}^2}  \, \prescript{RL}{0}{\mathcal D}_{t}^{\alpha} \, \frac{\partial^2{v}}{\partial{s}^2} 
\right) \,\, \Bigg|_{s=1}  = 0 ,
\end{align}
By substituting \eqref{Eq: weak form single mode - 2} in \eqref{Eq: weak form single mode - 1} we get
\begin{align}
\label{Eq: weak form single mode - 3}
&
\int_{0}^{1} \frac{\partial^2{v}}{\partial{t}^2} \, \tilde{v} \, ds  
+ 
\
\int_{0}^{1} \left(  
\frac{\partial^2{v}}{\partial{s}^2}+ \frac{\partial^2{v}}{\partial{s}^2} (\frac{\partial{v}}{\partial{s}})^2 
+ E_r
\prescript{RL}{0}{\mathcal D}_{t}^{\alpha} \, \Big[\frac{\partial^2{v}}{\partial{s}^2}(1 + \frac{1}{2} (\frac{\partial{v}}{\partial{s}})^2)\Big] \right.
\\ \nonumber
&
\left. + \frac{1}{2}  E_r 
(\frac{\partial{v}}{\partial{s}})^2 \, \prescript{RL}{0}{\mathcal D}_{t}^{\alpha} \, \frac{\partial^2{v}}{\partial{s}^2}
\right)
\, \frac{\partial^2{\tilde{v}}}{\partial{s}^2} \, ds
+\int_{0}^{1}\left( 
\frac{\partial{v}}{\partial{s}} \,  (\frac{\partial^2{v}}{\partial{s}^2})^2  
+ E_r
\frac{\partial{v}}{\partial{s}} \, \frac{\partial^2{v}}{\partial{s}^2} \,  \prescript{RL}{0}{\mathcal D}_{t}^{\alpha} \, \frac{\partial^2{v}}{\partial{s}^2}
\right)
\, \frac{\partial^2{\tilde{v}}}{\partial{s}^2}\tilde{v}^{\prime} \, ds
= f(t) ,
\end{align}
The modal discretization we used in \eqref{Eq: assumed mode} can be simplified as  $v(s,t)=q(t)\phi(s)$ where we choose $\phi(0)=\phi^{\prime}(0)=0$, and $\phi^{\prime\prime}(0)=\phi^{\prime\prime\prime}(0)=0$. Using the approximation, and having $ \tilde{v}=\phi(s)$, the evaluation of \eqref{Eq: weak form single mode - 3}  is
\begin{align}
\label{Eq: weak form single mode - 4}
&
\ddot{q}\int_{0}^{1} \phi \phi  ds  
\!+\! \int_{0}^{1} \left(  
q \phi^{\prime\prime} \!+\! q^{3}\phi^{\prime\prime}
{\phi^{\prime}}^2 \right)\phi^{\prime\prime}  ds
\!+\! \int_{0}^{1} E_r \phi^{\prime\prime} \prescript{RL}{0}{\mathcal D}_{t}^{\alpha}  q\phi^{\prime\prime} ds
\!+\!\int_{0}^{1}
\frac{1}{2}  E_r 
\phi^{\prime\prime}(\phi^{\prime})^2  \prescript{RL}{0}{\mathcal D}_{t}^{\alpha}  q^3 \phi^{\prime\prime}
ds
\\ \nonumber
&
\!+\! \int_{0}^{1}
\frac{1}{2}  E_r 
{\phi^{\prime}}^2 \phi^{\prime\prime} q^{2}  \prescript{RL}{0}{\mathcal D}_{t}^{\alpha}  q \phi^{\prime\prime}
ds
\!+\!\int_{0}^{1}\left( \phi^{\prime}   {\phi^{\prime\prime}}^2q^3 \!+\!  E_r
\phi^{\prime}\phi^{\prime\prime}\phi^{\prime\prime} q^2    \prescript{RL}{0}{\mathcal D}_{t}^{\alpha}  q 
\right)\phi^{\prime}
ds
\!=\! f(t) ,
\end{align}
considering \eqref{Eq: coeff unimodal} for the no lumped mass case as
\begin{align}
\label{Eq: weak form single mode - 5}
& \mathcal{M} = \int_{0}^{1} \phi^2 \, ds \,  \quad\mathcal{K}_l = \mathcal{C}_l = \int_{0}^{1}  {\phi^{\prime\prime}}^2 \,\, ds, \quad
\\ \nonumber
& \mathcal{K}_{nl} = \mathcal{C}_{nl} = \int_{0}^{1}  {\phi^{\prime}}^2 \, {\phi^{\prime\prime}}^2 \,\, ds,
\quad  \mathcal{M}_b = \int_{0}^{1} \phi \, ds  .
\end{align}
We are able to recover \eqref{Eq: weak form - discrete 2}
\begin{align}
\label{Eq: weak form single mode - 6}
\mathcal{M} \, \ddot{q} 
+ \mathcal{K}_l \, q + E_r \, \mathcal{C}_l \, \prescript{RL}{0}{\mathcal D}_{t}^{\alpha} q  
+ 2 \mathcal{K}_{nl} \, q^3 + \frac{E_r \, \mathcal{C}_{nl}}{2} \, \left( \prescript{RL}{0}{\mathcal D}_{t}^{\alpha} q^3 + 3 \, q^2 \, \prescript{RL}{0}{\mathcal D}_{t}^{\alpha} q  \right)
= -\mathcal{M}_b \, \ddot{V_b},
\end{align}
%%%%%%%%%%%%%%%%%%%%%%%%%%%%%%%%%%%%%%%%%
\section{Deriving the Linearized Equation}
\label{Sec: App. Linearization}
%%%%%%%%%%%%%%%%%%%%%%%%%%%%%%%%%%%%%%%%%
%

%
%%%%%%%%%%%%%%%%%%%%%%%%%%%
\subsection{Linearized Equation of Motion}
%%%%%%%%%%%%%%%%%%%%%%%%%%%
The source of nonlinearity in our problem is coming from geometry, so for the linear case we approximate the rotation angle as $\psi \simeq \frac {\partial{v}}{\partial{s}}$, and approximate the angular velocity and curvature as $ \frac {\partial{\psi}}{\partial{t}} \simeq \frac {\partial^2{v}}{\partial{t}\partial{s}}$, and $ \frac {\partial{\psi}}{\partial{s}} \simeq \frac {\partial^2{v}}{\partial{s}^2}$ respectively.
We recast the integral \eqref{Eq: work var} as $\delta W = \int_{0}^{L} \int_{A} \sigma \, \delta \varepsilon \, dA \,\, ds$ for the considered cantilever beam, in which the variation of strain is $\delta \varepsilon = -\eta \, \delta \frac{\partial{\psi}}{\partial{s}}$, using \eqref{Eq: strain-curvature relation}. 
Therefore, by assuming the constitutive equation \eqref{Eq: frac KV}, the variation of total work is expressed as
\begin{align}
\label{Eq: total work var linearized}
\delta w 
& = 
\int_{0}^{L} \int_{A} \left( -\eta \, E_{\infty} \, \frac{\partial{\psi}}{\partial{s}} -\eta \, E_{\alpha} \, \prescript{RL}{0}{\mathcal D}_{t}^{\alpha} \, \frac{\partial{\psi}}{\partial{s}} \right) \, (-\eta \, \delta \frac{\partial{\psi}}{\partial{s}}) \, dA  \,\, ds
\\ \nonumber
& = 
\int_{0}^{L} \left( E_{\infty} \, \left( \int_{A}  \, \eta^2 dA \right) \, \frac{\partial{\psi}}{\partial{s}} + E_{\alpha} \, \left( \int_{A}  \, \eta^2 dA \right) \, \prescript{RL}{0}{\mathcal D}_{t}^{\alpha} \, \frac{\partial{\psi}}{\partial{s}} \right) \, \delta \frac{\partial{\psi}}{\partial{s}} \,\, ds
\\ \nonumber
& = 
\int_{0}^{L} \left( E_{\infty} \, I \, \frac{\partial{\psi}}{\partial{s}} + E_{\alpha} \, I \, \prescript{RL}{0}{\mathcal D}_{t}^{\alpha} \, \frac{\partial{\psi}}{\partial{s}} \right) \, \delta \frac{\partial{\psi}}{\partial{s}} \,\, ds
\end{align}
where $I = \int_{A} \eta^2 \, dA$. By approximation \eqref{Eq: curvature}, we write the variation of curvature as
\begin{align}
\label{Eq: Curvature Variation Linearized}
\delta \frac{\partial{\psi}}{\partial{s}} 
=  \delta \frac{\partial^2{v}}{\partial{s}^2}.
\end{align}
Therefore, the variation of total energy becomes
\begin{align}
\label{Eq: Total Energy Linearized}
\delta w 
=
\int_{0}^{L}  
\left( E_{\infty} \, I \, \frac{\partial^2{v}}{\partial{s}^2}  + E_{\alpha} \, I \, \prescript{RL}{0}{\mathcal D}_{t}^{\alpha} \, \Big[\frac{\partial^2{v}}{\partial{s}^2}\Big]\right)  \, \delta \frac{\partial^2{v}}{\partial{s}^2} \,\, ds 
\end{align}
By expanding the terms and integrating by parts, we have
\begin{align}
\label{Eq: Total Energy Linearized Expantion}
\delta w 
& =
\int_{0}^{L}  
\frac{\partial^2}{\partial{s}^2}
\left( E_{\infty} \, I \, \frac{\partial^2{v}}{\partial{s}^2} + E_{\alpha} \, I \, \prescript{RL}{0}{\mathcal D}_{t}^{\alpha} \, \Big[\frac{\partial^2{v}}{\partial{s}^2}\Big]\right) \, 
\delta v \,\, ds 
\\ \nonumber
&  +
\left( E_{\infty} \, I \, \frac{\partial^2{v}}{\partial{s}^2}   + E_{\alpha} \, I \, \prescript{RL}{0}{\mathcal D}_{t}^{\alpha} \, \Big[\frac{\partial^2{v}}{\partial{s}^2}\Big]\right) 
\, \delta \frac{\partial{v}}{\partial{s}} \Bigg|_{0}^{L}
\\ \nonumber
& -
\frac{\partial}{\partial{s}}
\left( E_{\infty} \, I \, \frac{\partial^2{v}}{\partial{s}^2} + E_{\alpha} \, I \, \prescript{RL}{0}{\mathcal D}_{t}^{\alpha} \, \Big[\frac{\partial^2{v}}{\partial{s}^2}\Big]\right)\, \delta v \Bigg|_{0}^{L}
\end{align}
The prescribed geometry boundary conditions at the base of the beam, $s=0$, allow the variation of deflection and its first derivative to be zero at $s=0$, i.e. $\delta v(0,t) = \delta \frac{\partial{v}}{\partial{s}}(0,t) = 0$. Therefore,
\begin{align}
\label{Eq: total work var 4}
\delta w 
& =
\int_{0}^{L}  
\frac{\partial^2}{\partial{s}^2}
\left( E_{\infty} \, I \, \frac{\partial^2{v}}{\partial{s}^2} + E_{\alpha} \, I \, \prescript{RL}{0}{\mathcal D}_{t}^{\alpha} \, \Big[\frac{\partial^2{v}}{\partial{s}^2}\Big]\right)\, 
\delta v \,\, ds 
\\ \nonumber
&  +
\left( E_{\infty} \, I \, \frac{\partial^2{v}}{\partial{s}^2} + E_{\alpha} \, I \, \prescript{RL}{0}{\mathcal D}_{t}^{\alpha} \, \Big[\frac{\partial^2{v}}{\partial{s}^2}\Big]\right)\Bigg|_{s=L} \, \delta \frac{\partial{v}}{\partial{s}}(L,t) 
\\ \nonumber
& -
\frac{\partial}{\partial{s}}\left(
\left( E_{\infty} \, I \, \frac{\partial^2{v}}{\partial{s}^2} + E_{\alpha} \, I \, \prescript{RL}{0}{\mathcal D}_{t}^{\alpha} \, \Big[\frac{\partial^2{v}}{\partial{s}^2}\Big]\right)\right) \Bigg|_{s=L} \, \delta v(L,t) 
\end{align}
Let $\varrho$ be mass per unit volume of the beam, $M$ and $J$ be the mass and rotatory inertia of the lumped mass at the tip of beam. Using \eqref{Eq: displacement} and \eqref{Eq: velocity}, respectively, the kinetic energy is
\begin{align}
\label{Eq: Kinetic Linearized}
T 
& = \frac{1}{2} \, \int_{0}^{L} \int_{A} \varrho  \, (\frac{\partial{\textbf{r}}}{\partial{t}})^2 \, dA \, ds 
+ \frac{1}{2}M \left( (\frac{\partial{u}}{\partial{t}})^2  + (\frac{\partial{v}}{\partial{t}} + \dot V_b)^2 \right) \Big|_{s=L} + \frac{1}{2} J (\frac{\partial{\psi}}{\partial{t}})^2 \Big|_{s=L} , 
\\ \nonumber 
& = \frac{1}{2} \, \int_{0}^{L} \int_{A} \varrho  \, \Big\{ (\frac{\partial{u}}{\partial{t}} - \eta \, \frac{\partial{\psi}}{\partial{t}} \, \cos(\psi) )^2 + (\frac{\partial{v}}{\partial{t}}+ \dot V_b - \eta \, \frac{\partial{\psi}}{\partial{t}} \, \sin(\psi))^2 \Big\}  \,  dA \, ds 
\\ \nonumber
&
+ \frac{1}{2}M \left((\frac{\partial{u}}{\partial{t}})^2 + (\frac{\partial{v}}{\partial{t}} + \dot V_b)^2 \right) \Big|_{s=L} + \frac{1}{2} J (\frac{\partial{\psi}}{\partial{t}})^2 \Big|_{s=L} , 
\\ \nonumber
& = 
\frac{1}{2} \, \int_{0}^{L} \int_{A} \varrho  \,
\Big\{ (\frac{\partial{u}}{\partial{t}} )^2 - 2 \eta \, \frac{\partial{u}}{\partial{t}}  \, \frac{\partial{\psi}}{\partial{t}}  \, \cos(\psi) + \eta^2 \, (\frac{\partial{\psi}}{\partial{t}})^2 \, \cos^2(\psi) 
+ (\frac{\partial{v}}{\partial{t}} )^2 + {\dot V_b}^2 + 2 \frac{\partial{v}}{\partial{t}}  \, \dot V_b
\\ \nonumber
& 
- 2 \eta \, (\frac{\partial{v}}{\partial{t}} + \dot V_b) \frac{\partial{\psi}}{\partial{t}}  \, \sin(\psi) + \eta^2 (\frac{\partial{\psi}}{\partial{t}})^2 \, \sin^2(\psi) 
\Big\}  \,  dA \, ds + \frac{1}{2}M \left((\frac{\partial{u}}{\partial{t}})^2 + (\frac{\partial{v}}{\partial{t}} + \dot V_b)^2 \right) \Big|_{s=L}
\\ \nonumber
&
+ \frac{1}{2} J (\frac{\partial{\psi}}{\partial{t}})^2 \Big|_{s=L} , 
\\ \nonumber
& = 
\frac{1}{2} \, \int_{0}^{L} \int_{A} \varrho  \,
\Big\{ (\frac{\partial{u}}{\partial{t}})^2 + (\frac{\partial{v}}{\partial{t}})^2 + {\dot V_b}^2 + 2 \frac{\partial{v}}{\partial{t}} \, \dot V_b 
- 2 \eta \, \frac{\partial{u}}{\partial{t}}\, \frac{\partial{\psi}}{\partial{t}} \, \cos(\psi) + \eta^2 \, (\frac{\partial{\psi}}{\partial{t}})^2 
\\ \nonumber
&
- 2 \eta \, (\frac{\partial{v}}{\partial{t}} + \dot V_b) \frac{\partial{\psi}}{\partial{t}} \, \sin(\psi)  
\Big\}  \,  dA \, ds + \frac{1}{2}M \left((\frac{\partial{u}}{\partial{t}} )^2 + (\frac{\partial{v}}{\partial{t}}  + \dot V_b)^2 \right) \Big|_{s=L} + \frac{1}{2} J (\frac{\partial{\psi}}{\partial{t}})^2 \Big|_{s=L} .
\end{align}
Let
\begin{align*}
\rho = \int_{A} \varrho dA , \quad
\mathcal{J}_1 = \int_{A} \varrho \, \eta \,  dA , \quad
\mathcal{J}_2 = \int_{A} \varrho \, \eta^2 dA.
\end{align*}
$\rho$ is the mass per unit length of the beam, $\mathcal{J}_1$ is the first moment of inertia and is zero because the reference point of coordinate system attached to the cross section coincides with the mass centroid, and $\mathcal{J}_2$ is the second moment of inertia, which is very small for slender beam and can be ignored \cite{hamdan1997large}. Assuming that the velocity along the length of the beam, $\dot u$, is relatively small compared to the lateral velocity $\dot v + \dot V_b$, the kinetic energy of the beam can be reduced to
\begin{align}
\label{Eq: Kinetic Linearized - 2}
T = 
\frac{1}{2} \, \rho \int_{0}^{L} ( \frac{\partial{v}}{\partial{t}} + \dot V_b)^2  \,   ds 
+ \frac{1}{2}M (\frac{\partial{v}}{\partial{t}} + \dot V_b)^2  \Big|_{s=L} + \frac{1}{2} J (\frac{\partial{\psi}}{\partial{t}})^2 \Big|_{s=L}, 
\end{align}
where its variation can be taken as
\begin{align}
\label{Eq: Kinetic Variation Linearized}
\delta T 
& \!=\! \rho \int_{0}^{L} (\frac{\partial{v}}{\partial{t}} + \dot{V_b}) \, \delta \frac{\partial{v}}{\partial{t}} \, ds + M (\frac{\partial{v}}{\partial{t}} + \dot{V_b}) \, \delta \frac{\partial{v}}{\partial{t}} \Big|_{s=L} +  J \frac{\partial{\psi}}{\partial{t}} \, \delta \frac{\partial{\psi}}{\partial{t}} \Big|_{s=L} ,
\end{align}
in which $\frac{\partial{\psi}}{\partial{t}}$ is given in \eqref{Eq: angular velocity} and $\delta \frac{\partial{\psi}}{\partial{t}}$ can be obtained as $\delta \frac{\partial{\psi}}{\partial{t}} \simeq \delta \frac{\partial^2{v}}{\partial{t}^2} $. Therefore, 
\begin{align}
\label{Eq: Kinetic Variation Linearized - 2}
\delta T 
& \!\simeq\! \rho \!\int_{0}^{L}\! (\frac{\partial{v}}{\partial{t}} \!+\! \dot{V_b})  \delta \frac{\partial{v}}{\partial{t}}  ds 
\!+\! M (\frac{\partial{v}}{\partial{t}} \!+\! \dot{V_b})  \delta \frac{\partial{v}}{\partial{t}} \Big|_{s=L} 
\!+\! J \left( \frac{\partial^2{v}}{\partial{t}\partial{s}} ( 1 \!+\! (\frac{\partial{v}}{\partial{s}})^2) \delta \frac{\partial^2{v}}{\partial{t}\partial{s}} \!+\! \frac{\partial{v}}{\partial{s}} (\frac{\partial^2{v}}{\partial{t}\partial{s}})^2 \delta \frac{\partial{v}}{\partial{s}} \right) \Big|_{s=L} .
\end{align}
The time integration of $\delta T$ takes the following form through integration by parts
\begin{align}
\label{Eq: Kinetic Var Linearized- 2}
&\int_{t_1}^{t_2} \delta T \, dt 
\\ \nonumber
= & \int_{t_1}^{t_2} 
\Bigg\{
\rho \int_{0}^{L} (\frac{\partial{v}}{\partial{t}} + \dot{V_b}) \, \delta \frac{\partial{v}}{\partial{t}} \, ds 
+ M (\frac{\partial{v}}{\partial{t}} + \dot{V_b}) \, \delta \frac{\partial{v}}{\partial{t}} \Big|_{s=L} 
\\ \nonumber
& 
+ J \left( \frac{\partial^2{v}}{\partial{t}\partial{s}} ( 1 + (\frac{\partial{v}}{\partial{s}})^2) \delta \frac{\partial^2{v}}{\partial{t}\partial{s}} + \frac{\partial{v}}{\partial{s}} (\frac{\partial^2{v}}{\partial{t}\partial{s}})^2 \delta \frac{\partial{v}}{\partial{s}} \right) \Big|_{s=L}
\Bigg\} \, dt
\\ \nonumber
= &
\int_{t_1}^{t_2} \rho \int_{0}^{L} (\frac{\partial{v}}{\partial{t}} + \dot{V_b}) \, \delta \frac{\partial{v}}{\partial{t}} \, ds \, dt
+ M \int_{t_1}^{t_2} (\frac{\partial{v}}{\partial{t}} + \dot{V_b}) \, \delta \frac{\partial{v}}{\partial{t}} \Big|_{s=L} \, dt
\\ \nonumber
& 
+ J \int_{t_1}^{t_2} \left( \frac{\partial^2{v}}{\partial{t}\partial{s}} ( 1 + (\frac{\partial{v}}{\partial{s}})^2) \delta \frac{\partial^2{v}}{\partial{t}\partial{s}} + \frac{\partial{v}}{\partial{s}} (\frac{\partial^2{v}}{\partial{t}\partial{s}})^2 \delta \frac{\partial{v}}{\partial{s}}  \right) \Big|_{s=L} \, dt
\\ \nonumber
= &
\rho \int_{0}^{L} \int_{t_1}^{t_2}  (\frac{\partial{v}}{\partial{t}} + \dot{V_b}) \, \delta \frac{\partial{v}}{\partial{t}} \, dt \, ds
+ M \int_{t_1}^{t_2} (\frac{\partial{v}}{\partial{t}} + \dot{V_b}) \, \delta \frac{\partial{v}}{\partial{t}}  \, dt \, \Big|_{s=L}
\\ \nonumber
& 
+ J \int_{t_1}^{t_2} \left( \frac{\partial^2{v}}{\partial{t}\partial{s}} ( 1 + (\frac{\partial{v}}{\partial{s}})^2) \delta \frac{\partial^2{v}}{\partial{t}\partial{s}} + \frac{\partial{v}}{\partial{s}} (\frac{\partial^2{v}}{\partial{t}\partial{s}})^2 \delta \frac{\partial{v}}{\partial{s}} \right)  \, dt \, \Big|_{s=L}
\\ \nonumber
= &
\rho \int_{0}^{L} \left[ (\frac{\partial{v}}{\partial{t}} \!+\! \dot{V_b})  \delta{v} \Big|_{t_1}^{t_2} \!-\! \int_{t_1}^{t_2}  (\frac{\partial^2{v}}{\partial{t}^2} \!+\! \ddot{V_b})  \delta {v}  dt \right] ds
\!+\! M (\frac{\partial{v}}{\partial{t}} \!+\! \dot{V_b})  \delta v \Big|_{s=L} \Big|_{t_1}^{t_2} 
\!-\! M \int_{t_1}^{t_2} (\frac{\partial^2{v}}{\partial{t}^2} \!+\! \ddot{V_b})  \delta v   dt  \Big|_{s=L}
\\ \nonumber
& \quad
+ J  \frac{\partial^2{v}}{\partial{t}\partial{s}} ( 1 + (\frac{\partial{v}}{\partial{s}})^2) \delta \frac{\partial{v}}{\partial{s}}  \Big|_{s=L}  \Big|_{t_1}^{t_2}
- J \int_{t_1}^{t_2} 
\left(  \frac{\partial^3{v}}{\partial{t}^2\partial{s}} ( 1 + (\frac{\partial{v}}{\partial{s}})^2) + \frac{\partial{v}}{\partial{s}} (\frac{\partial^2{v}}{\partial{t}\partial{s}})^2 \right)  
\delta \frac{\partial{v}}{\partial{s}}   dt  \Big|_{s=L}
\\ \nonumber
\!=\! &
\!-\! \!\int_{t_1}^{t_2}\! \Bigg\{
\rho \!\int_{0}^{L}\!  (\frac{\partial^2{v}}{\partial{t}^2} \!+\! \ddot{V_b})  \delta {v}  ds 
\!+\! M (\frac{\partial^2{v}}{\partial{t}^2} \!+\! \ddot{V_b})  \delta v   \Big|_{s=L}
\!+\! J \left(  \frac{\partial^3{v}}{\partial{t}^2\partial{s}}( 1 \!+\! (\frac{\partial{v}}{\partial{s}})^2) \!+\! \frac{\partial{v}}{\partial{s}} (\frac{\partial^2{v}}{\partial{t}\partial{s}})^2  \right)  \delta \frac{\partial{v}}{\partial{s}}   \Big|_{s=L}
\Bigg\}  dt,
\end{align}
where we consider that $\delta v = \delta \frac{\partial{v}}{\partial{s}} = 0$ at $t=t_1$ and $t=t_2$. Therefore, the extended Hamilton's principle takes the form
\begin{align}
\label{Eq: Extended Hamilton Linearized}
&
\int_{t_1}^{t_2} \Bigg\{
\int_{0}^{L} \Bigg[
- \rho (\frac{\partial^2{v}}{\partial{t}^2}) -\frac{\partial^2}{\partial{s}^2}
\left( E_{\infty} \, I \, \frac{\partial^2{v}}{\partial{s}^2} + E_{\alpha} \, I \, \prescript{RL}{0}{\mathcal D}_{t}^{\alpha} \, \Big[\frac{\partial^2{v}}{\partial{s}^2}\Big]\right) \Bigg]\delta v \, ds
\\ \nonumber
&
- M (\frac{\partial^2{v}}{\partial{t}^2} + \ddot{V_b}) \Big|_{s=L} \, \delta v(L,t)  -
\left( E_{\infty} \, I \, \frac{\partial^2{v}}{\partial{s}^2} + E_{\alpha} \, I \, \prescript{RL}{0}{\mathcal D}_{t}^{\alpha} \, \Big[\frac{\partial^2{v}}{\partial{s}^2}\Big]\right) 
\Bigg|_{s=L} \, \delta \frac{\partial{v}}{\partial{s}}(L,t) 
\\ \nonumber
&  +
\frac{\partial}{\partial{s}}
\left( E_{\infty} \, I \, \frac{\partial^2{v}}{\partial{s}^2} + E_{\alpha} \, I \, \prescript{RL}{0}{\mathcal D}_{t}^{\alpha} \, \Big[\frac{\partial^2{v}}{\partial{s}^2}\Big]\right) \Bigg|_{s=L} \, \delta v(L,t) 
\,\, \Bigg\} \, dt = 0.
\end{align}
Invoking the arbitrariness of virtual displacement $\delta v$, we obtain the strong form of the equation of motion as:
\begin{align}
\label{Eq: Eqn of Motion Linearized}
& \rho \, \frac{\partial^2{v}}{\partial{t}^2}  
+ 
E_{\infty} \, I \, \frac{\partial^2}{\partial{s}^2}\left(  \frac{\partial^2{v}}{\partial{s}^2}  \right)
+ 
E_{\alpha} \, I \, \frac{\partial^2}{\partial{s}^2}\left(   \, \prescript{RL}{0}{\mathcal D}_{t}^{\alpha} \, \Big[\frac{\partial^2{v}}{\partial{s}^2} \Big]\right)=0
\end{align}
which is subject to the following natural boundary conditions:
\begin{align}
\label{Eq: Natural BC Linearized}
&E_{\infty} \, I \, \frac{\partial^2{v}}{\partial{s}^2}  + E_{\alpha} \, I \,  \, \prescript{RL}{0}{\mathcal D}_{t}^{\alpha} \, \Big[\frac{\partial^2{v}}{\partial{s}^2}(1 + \frac{1}{2} (\frac{\partial{v}}{\partial{s}})^2)\Big]
\,\, \Bigg|_{s=L} = 0 ,
\\ \nonumber
& 
M (\frac{\partial^2{v}}{\partial{t}^2} )
- \frac{\partial}{\partial{s}}\left(
E_{\infty} \, I \, \frac{\partial^2{v}}{\partial{s}^2}   + E_{\alpha} \, I \, \, \prescript{RL}{0}{\mathcal D}_{t}^{\alpha} \, \Big[\frac{\partial^2{v}}{\partial{s}^2}\Big] \right)\, \Bigg|_{s=L}  = 0 .
\end{align}
Therefore, the strong form reads as: find $v \in V$ such that
\begin{align}
\label{Eq: eqn of motion - 2 Linearized}
& m \, \frac{\partial^2{v}}{\partial{t}^2} 
+ 
\frac{\partial^2}{\partial{s}^2}\left(  
\frac{\partial^2{v}}{\partial{s}^2} 
+ E_r \prescript{RL}{0}{\mathcal D}_{t}^{\alpha} \, \Big[\frac{\partial^2{v}}{\partial{s}^2}\Big] 
\right)
= 0 ,
\end{align}
subject to the following boundary conditions:
\begin{align}
\label{Eq: bc Linearized}
& v \, \Big|_{s=0} = \frac{\partial{v}}{\partial{s}} \, \Big|_{s=0} = 0 ,
\\ \nonumber
& 
\left( 
\frac{\partial^2{v}}{\partial{s}^2}
+ E_r \, \prescript{RL}{0}{\mathcal D}_{t}^{\alpha} \, \Big[\frac{\partial^2{v}}{\partial{s}^2}\Big] 
\,
\right) \,\, \Bigg|_{s=L} = 0 ,
\\ \nonumber
& 
\frac{M m}{\rho} (\frac{\partial^2{v}}{\partial{t}^2})
- \frac{\partial}{\partial{s}} \left(
\frac{\partial^2{v}}{\partial{s}^2} + E_r \, \prescript{RL}{0}{\mathcal D}_{t}^{\alpha} \, \Big[\frac{\partial^2{v}}{\partial{s}^2} \Big]
\, \prescript{RL}{0}{\mathcal D}_{t}^{\alpha} \, \frac{\partial^2{v}}{\partial{s}^2}  
\right)\,\, \Bigg|_{s=L}  = 0 ,
\end{align}
where $m = \frac{\rho}{ E_{\infty} \, I}$ and $E_r = \frac{E_{\alpha}}{E_{\infty}}$. 

%%%%%%%%%%%%%%%%%%%%%%%%%%%
\subsection{Nondimensionalization of Linearized Equation of Motion}
%%%%%%%%%%%%%%%%%%%%%%%%%%%
%
Let the dimensionless variables
\begin{align}
s^{*} = \frac{s}{L}, \quad 
v^{*} = \frac{v}{L}, \quad 
t^{*} = t \left(\frac{1}{m L^4}\right)^{1/2}, \quad 
E_r^{*} = E_r \left(\frac{1}{m L^4}\right)^{\alpha/2}, \quad 
\end{align}
We obtain the following dimensionless equation by substituting the above dimensionless variables. 
\begin{align}
\label{Eq: eqn of motion - dimless - 0 Linearized}
& m \frac{L}{mL^4} \, \frac{\partial^2 v^*}{\partial {t^*}^2} + 
\frac{1}{L^2} \frac{\partial^2}{\partial {s^*}^2} 
\Bigg[  
\frac{L}{L^2} \frac{\partial^2 v^*}{\partial {s^*}^2} 
+ E_r^* (mL^4)^{\alpha/2} \frac{1}{(mL^4)^{\alpha/2}} \frac{L}{L^2} \prescript{RL}{0}{\mathcal D}_{t^*}^{\alpha} \frac{\partial^2 v^*}{\partial {s^*}^2}
\Bigg] = 0  ,
\end{align}
which can be simplified to
\begin{align}
\label{Eq: eqn of motion - dimless - 1 Linearized}
\frac{\partial^2 v^*}{\partial {t^*}^2}  
& + 
\frac{\partial^2}{\partial {s^*}^2} 
\Bigg[  
\frac{\partial^2 v^*}{\partial {s^*}^2} 
+ E_r^* \prescript{RL}{0}{\mathcal D}_{t^*}^{\alpha} \frac{\partial^2 v^*}{\partial {s^*}^2}
\Bigg]
= 0 ,
\end{align}
The dimensionless boundary conditions are also obtained by substituting dimensionless variables in \eqref{Eq: bc Linearized}. We can show similarly that they preserve their structure as: 
\begin{align*}
& v^{*} \, \Big|_{s^{*}=0} = \frac{\partial v^{*}}{\partial s^{*}} \, \Big|_{s^{*}=0} = 0 ,
\\ \nonumber
& 
\frac{1}{L} \Bigg[  
\frac{\partial^2 v^*}{\partial {s^*}^2} 
+ E_r^* \prescript{RL}{0}{\mathcal D}_{t^*}^{\alpha} \frac{\partial^2 v^*}{\partial {s^*}^2}
\Bigg]
\,\, \Bigg|_{s^{*}=1} = 0 ,
\\ \nonumber
& 
\frac{M^* \rho L m}{\rho} \frac{L}{mL^4} \left( \frac{\partial^2 v^{*}}{\partial^2 t^{*}} \right)
-\frac{1}{L^2} \frac{\partial v^{*}}{\partial s^{*}} \Bigg[  
\frac{\partial^2 v^*}{\partial {s^*}^2} 
+ E_r^* \prescript{RL}{0}{\mathcal D}_{t^*}^{\alpha} \frac{\partial^2 v^*}{\partial {s^*}^2}
\Bigg]
\,\, \Bigg|_{s^{*}=1}  = 0 ,
\end{align*}
\noindent Therefore, the dimensionless equation of motion becomes (after dropping $^*$ for the sake of simplicity)
\begin{align}
\label{Eq: eqn of motion - dimless - 2 Linearized}
\frac{\partial^2{v}}{\partial{t}^2} 
& + 
\frac{\partial^2}{\partial{s}^2} \left(  
\frac{\partial^2{v}}{\partial{s}^2}
+ E_r \prescript{RL}{0}{\mathcal D}_{t}^{\alpha} \, \Big[\frac{\partial^2{v}}{\partial{s}^2}\Big] 
\right)
= 0 ,
\end{align}
which is subject to the following dimensionless boundary conditions
\begin{align}
\label{Eq: bc - dimless Linearized}
& v \, \Big|_{s=0} = \frac{\partial{v}}{\partial{s}} \, \Big|_{s=0} = 0 ,
\\ \nonumber
& 
\left( 
\frac{\partial^2{v}}{\partial{s}^2}
+ E_r \, \prescript{RL}{0}{\mathcal D}_{t}^{\alpha} \, \Big[\frac{\partial^2{v}}{\partial{s}^2}\Big]  
\right) \,\, \Bigg|_{s=1} = 0 ,
\\ \nonumber
& 
- \frac{\partial}{\partial{s}} \left(
\frac{\partial^2{v}}{\partial{s}^2}
+ E_r \, \prescript{RL}{0}{\mathcal D}_{t}^{\alpha} \, \Big[\frac{\partial^2{v}}{\partial{s}^2} \Big]
\right)
\,\, \Bigg|_{s=1}  = 0 ,
\end{align}
We obtain the weak form of the problem by multiplying the strong form \ref{Eq: eqn of motion - dimless - 2 Linearized} with proper test functions $\tilde{v}(s) \in \tilde{V} $ and integrating over the dimensionless spatial computational domain $\Omega_s = [0,1]$. Therefore, by changing the order of integral and temporal derivatives, and through integration by parts, the weak form of problem can be written as:
\begin{align}
\label{Eq: weak form - 1 Linearized}
& \int_{0}^{1} \frac{\partial^2{v}}{\partial{t}^2} \, \tilde{v} \, ds  
+ 
\int_{0}^{1} 
\frac{\partial^2}{\partial{s}^2}\left(  
\frac{\partial^2{v}}{\partial{s}^2}
+ E_r
\prescript{RL}{0}{\mathcal D}_{t}^{\alpha} \, \Big[\frac{\partial^2{v}}{\partial{s}^2}\Big] 
\right)
\, \tilde{v} \, ds
= 0
\end{align}
where we transfer the spatial derivative load to the test function through integration by parts as
\begin{align}
\label{Eq: weak form - 3 linearized }
& \frac{\partial^2}{\partial {t}^2} \, \int_{0}^{1} v \, \tilde{v} \, ds 
+
\int_{0}^{1} \frac{\partial^2{v}}{\partial{s}^2}\, \frac{\partial^2{\tilde{v}}}{\partial{s}^2} \, ds
+ E_r 
\int_{0}^{1} \prescript{RL}{0}{\mathcal D}_{t}^{\alpha} \, \Big[\frac{\partial^2{v}}{\partial{s}^2}\Big]\,\, \frac{\partial^2{\tilde{v}}}{\partial{s}^2} \, ds = 0.
\end{align}
Using \eqref{Eq: assumed mode} and \eqref{Eq: Solution/Test Space}, the problem  \eqref{Eq: weak form - 3 linearized } read as: find $v_N \in V_N$ such that
\begin{align}
\label{Eq: weak form - discrete Linearized}
& \frac{\partial^2}{\partial {t}^2} \,  \int_{0}^{1} v_N \, \tilde{v}_N \, ds 
+ 
\int_{0}^{1} \frac{\partial^2{{v}_N}}{\partial{s}^2} \, \frac{\partial^2{\tilde{v}_N}}{\partial{s}^2} \, ds
+ E_r 
\int_{0}^{1} \prescript{RL}{0}{\mathcal D}_{t}^{\alpha} \, \Big[\frac{\partial^2{{v}_N}}{\partial{s}^2}\Big] \,\,  \frac{\partial^2{\tilde{v}_N}}{\partial{s}^2} \, ds
= 0
\end{align}
for all $ \tilde{v}_N \in \tilde{V}_N$. \\
By substituting unimodal discretization in \eqref{Eq: weak form - discrete Linearized} with the exact same method we used in section \ref{subsection}, we obtain the unimodal governing equation of motion as
\begin{align}
\label{Eq: weak form - discrete 2 Linearized}
\mathcal{M} \, \ddot{q}
+ \mathcal{K}_l \, q + E_r \, \mathcal{C}_l \, \prescript{RL}{0}{\mathcal D}_{t}^{\alpha} q  
= 0
\end{align}
in which
\begin{align}
\label{Eq: coeff unimodal Linearized}
& \mathcal{M} = \int_{0}^{1} \phi^2 \, ds \,\quad \quad\mathcal{K}_l = \mathcal{C}_l = \int_{0}^{1}  {\phi^{\prime\prime}}^2 \,\, ds, 
\end{align}
Therefore in absence of base excitation \eqref{Eq: weak form - discrete 2 Linearized} takes the form 
\begin{align}
\label{Eq: weak form linear - discrete 2 Linearized}
\ddot{q} + E_r \, c_l \, \prescript{RL}{0}{\mathcal D}_{t}^{\alpha} q  + k_l \, q  = 0
\end{align}
in which the coefficients $c_l = \frac{\mathcal{C}_l}{\mathcal{M}}$ and $k_l = \frac{\mathcal{K}_l}{\mathcal{M}}$ are given in \eqref{Eq: coeff unimodal}
%%%%%%%%%%%%%%%%%%%%%%%%%%%%%%%%%%%%%%%%%
\section{Eigenvalue Problem of Linear Model}
\label{Sec: App. Eigenvalue Problem of Linear Model} 
\begin{figure}[t]
	\centering
	\begin{subfigure}{0.4\textwidth}
		\centering
		\includegraphics[width=1\linewidth]{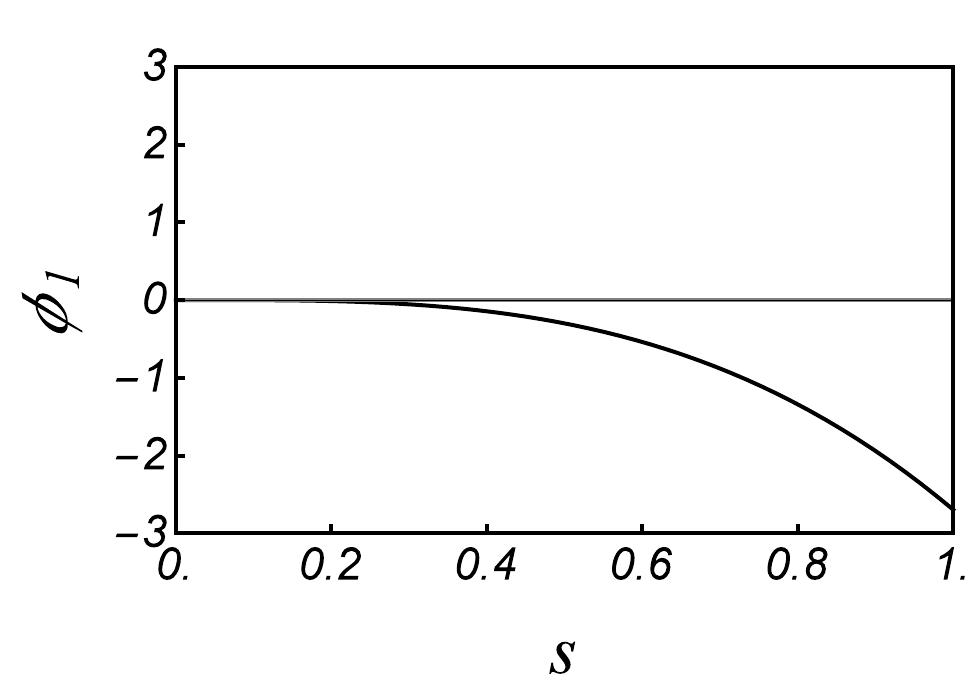}
	\end{subfigure}
	\begin{subfigure}{0.4\textwidth}
		\centering
		\includegraphics[width=1\linewidth]{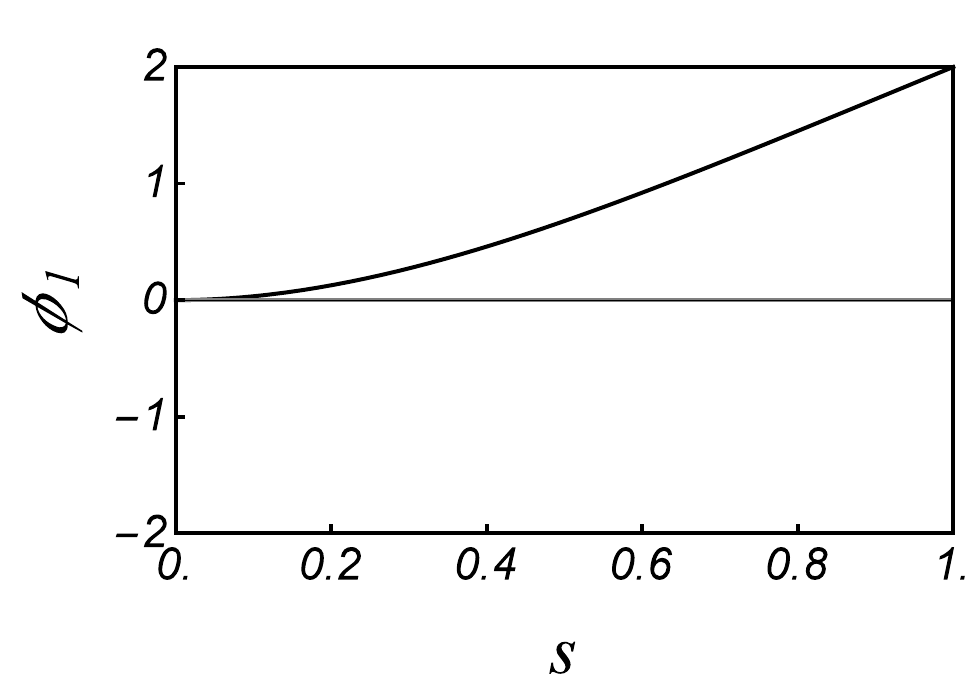}
	\end{subfigure}
	\caption{Left: The first eigenfunctions, $X_1(s)$, of the undamped linear counterpart of our model. It is used as the spatial functions in the single mode approximation. Right: The first eigenfunctions, $X_1(s)$, of the undamped linear counterpart of our model with no lumped mass at the tip. It is used as the spatial functions in the single mode approximation.}
	\label{Fig: Mode Shape Tip Mass}
\end{figure}

%%%%%%%%%%%%%%%%%%%%%%%%%%%%%%%%%%%%%%%%%
%
The assumed modes $\phi_i(s)$ in discretization \eqref{Eq: assumed mode} are obtained by solving the corresponding eigenvalue problem of free vibration of undamped linear counterparts to our model. Thus, the dimensionless linearized undamped equation of motion takes the form
\begin{align}
\label{Eq: Linear Euler Bernouli Beam}
& \frac{\partial^2}{\partial t^2}  v(s,t)
+ \frac{\partial^4}{\partial s^4} v(s,t)
= 0.
\end{align}
subject to linearized boundary conditions:
\begin{align}
\label{Eq: Linear BC}
& v(0,t) = 0 , &&  \frac{\partial^2{v}}{\partial{s}^2}(1,t) = -J \, \ddot{v}^{\prime}(1,t), 
\\ \nonumber
& v^{\prime}(0,t) = 0,
&& v^{\prime\prime\prime}(1,t) = M \, \ddot{v}(1,t) ,
\end{align}
where $ \dot{( \,\,\,)} = \frac{d }{dt}$ and $ ( \,\,\,)^{'} = \frac{d }{ds}$. We derive the corresponding eigenvalue problem by applying the separation of variables, i.e. $v(x,t) = X(s) T(t)$ to \eqref{Eq: Linear Euler Bernouli Beam}. Therefore, 
\begin{align}
\label{Eq: separ of vara}
& \ddot{T}(t) X(s) + T(t) X^{''''}(s) = 0,
\\
\nonumber
& \frac{\ddot{T}(t)}{T(t)}  + \frac{X^{''''}(s)}{X(s)} = 0,
\\
\nonumber
& \frac{\ddot{T}(t)}{T(t)} = -  \frac{X^{''''}(s)}{X(s)} = \lambda,
\end{align}
which gives the following equations
\begin{align}
\label{Eq: eigenProb 1}
&\ddot{T}(t) + \omega^2 T(t) = 0 , 
\\
\label{Eq: eigenProb 2}
&X^{\prime\prime\prime\prime}(s) - \beta^4 X(s) = 0, 
\end{align}
where $\beta^4 = \omega^2$ and the boundary conditions are
\begin{align*}
%\label{Eq: Linear BC - 2}
%
&X(0) = 0, && X^{\prime\prime}(1) =  J \, \omega^2 \, X^{\prime}(1), \\ \nonumber
&X^{\prime}(0) = 0,
&&X^{\prime\prime\prime}(1) = - M \, \omega^2 \, X(1).
\end{align*}
the solution to \eqref{Eq: eigenProb 2} is of the form $X(s) = A \sin(\beta s) + B \cos(\beta s) + C \sinh(\beta s) + D \cosh(\beta s)$, where $C = - A$ and $D = -B$, using the boundary conditions at $s=0$. Therefore, $$X(s) = A \left( \sin(\beta s) - \sinh(\beta s) \right) + B \left( \cos(\beta s) - \cosh(\beta s) \right).$$ Applying the first bondary condition at $s = 1$, i.e. $X^{\prime\prime}(1) =  J \, \omega^2 \, X^{\prime}(1)$ gives 
$$B = - \frac{ \sin(\beta) + \sinh(\beta) + J \beta^3 (\cos(\beta) - \cosh(\beta))}{ \cos(\beta) + \cosh(\beta)- J \beta^3 (\sin(\beta) - \sinh(\beta))} A, $$
that results in
$$X(s) = A \left[ \left( \sin(\beta s) - \sinh(\beta s) \right) - \frac{ \sin(\beta) + \sinh(\beta) + J \beta^3 (\cos(\beta) - \cosh(\beta))}{ \cos(\beta) + \cosh(\beta)- J \beta^3 (\sin(\beta) - \sinh(\beta))} \left( \cos(\beta s) - \cosh(\beta s) \right) \right].$$
Finally, using the second boundary condition at $s=1$ gives the following transcendental equation to solve for $\beta$'s for the case that $M=J=1$ :
\begin{align}
\label{Eq: freq}
&-\left(1 + \beta ^4 + \cos(\beta ) \cosh(\beta ) \right)
+\beta \, \left( \sin (\beta ) \cosh (\beta )-\cos (\beta ) \sinh (\beta ) \right) 
\\ 
&+\beta ^3 \, \left( \sin (\beta ) \cosh (\beta )-\sinh (\beta ) \cosh (\beta ) \right)
+ \beta ^4 \, \left( \sin (\beta ) \sinh (\beta ) + \cos (\beta ) \cosh (\beta ) \right) 
= 0 .
\end{align}
The first eigenvalue is computed as $\beta_1^2 = \omega_1 = 1.38569$, which results to the following first normalized eigenfunction, given in Fig. \ref{Fig: Mode Shape Tip Mass} (left). 

$$X_1(s) = 5.50054 \sin (\beta_1 s)-0.215842 \cos (\beta_1 s)-5.50054 \sinh (\beta_1 s)+0.215842 \cosh (\beta_1 s),$$

We note that \eqref{Eq: freq} reduces to $1 + \cos(\beta) \cosh(\beta) = 0 $ for the case that there is no lumped mass at the tip of beam; this in fact gives the natural frequencies of a linear cantilever beam. In this case, the first eigenvalue is computed as $\beta_1^2 = \omega_1 = 3.51602$, which results to the following first normalized eigenfunction, given in Fig. \ref{Fig: Mode Shape Tip Mass} (right). 

$$X_1(s) = 0.734096 \sin (\beta_1 s) - \cos (\beta_1 s)-0.734096 \sinh (\beta_1 s)+ \cosh (\beta_1 s).$$

\bibliographystyle{siamplain}
\bibliography{Structural_Health_arXiv}

\end{document}